\documentclass[aos,preprint]{imsart}
\usepackage{graphicx}
\usepackage{amssymb}
\usepackage{epstopdf}
\RequirePackage[authoryear]{natbib}
\usepackage{amsmath, mathabx}
\usepackage{mathrsfs}
\usepackage{hyperref}
\usepackage{comment}
\usepackage{xcolor}
\usepackage{subfigure}
\usepackage{enumitem}
\usepackage{algorithm}
\usepackage{algpseudocode}
\usepackage{stmaryrd}
\usepackage[T1]{fontenc}

\usepackage{blkarray}

\usepackage{url}
\usepackage{longtable}
\usepackage{mathtools}
\usepackage{multirow}
\usepackage{array}

\startlocaldefs

\renewcommand{\b}{{\mathbf{b}}}

\newcommand{\e}{{\mathbf{e}}}
\newcommand{\x}{{\mathbf{x}}}
\newcommand{\z}{{\mathbf{z}}}
\renewcommand{\j}{{\mathbf{j}}}
\renewcommand{\u}{{\mathbf{u}}}
\newcommand{\w}{\mathbf{w}}
\newcommand{\scA}{\mathscr{A}}
\renewcommand{\a}{{\mathbf{a}}}
\newcommand{\y}{{\mathbf{y}}}
\newcommand{\X}{{\mathbf{X}}}
\newcommand{\A}{{\mathbf{A}}}
\newcommand{\Y}{{\mathbf{Y}}}
\newcommand{\B}{{\mathbf{B}}}

\newcommand{\D}{{\mathbf{D}}}
\newcommand{\R}{{\mathbf{R}}}
\newcommand{\U}{{\mathbf{U}}}
\newcommand{\V}{{\mathbf{V}}}
\newcommand{\W}{{\mathbf{W}}}
\newcommand{\Q}{{\mathbf{Q}}}
\newcommand{\I}{{\mathbf{I}}}
\newcommand{\Z}{{\mathbf{Z}}}

\renewcommand{\O}{{\mathbf{O}}}

\newcommand{\scT}{{\mathscr{T}}}

\newcommand{\bcX}{{\mathbfcal{X}}}
\newcommand{\bcA}{{\mathbfcal{A}}}
\newcommand{\bcY}{{\mathbfcal{Y}}}

\newcommand{\bcZ}{{\mathbfcal{Z}}}
\newcommand{\bcB}{{\mathbfcal{B}}}
\newcommand{\bcS}{{\mathbfcal{S}}}
\newcommand{\bcW}{{\mathbfcal{W}}}

\newcommand{\bcT}{{\mathbfcal{T}}}

\newcommand{\bcE}{{\mathbfcal{E}}}

\newcommand{\br}{\mathbf{r}}
\newcommand{\bp}{\mathbf{p}}

\newcommand{\cH}{\mathcal{H}}

\newcommand{\cM}{\mathcal{M}}
\newcommand{\cN}{\mathcal{N}}
\newcommand{\cP}{\mathcal{P}}

\newcommand{\cT}{\mathcal{T}}

\newcommand{\bbN}{\mathbb{N}}
\newcommand{\bbP}{\mathbb{P}}
\newcommand{\bbE}{\mathbb{E}}
\newcommand{\bbR}{\mathbb{R}}

\newcommand{\bbO}{\mathbb{O}}
\newcommand{\bbM}{\mathbb{M}}
\newcommand{\bbQ}{\mathbb{Q}}

\newcommand{\balpha}{\boldsymbol{\alpha}}
\newcommand{\bbeta}{\boldsymbol{\beta}}

\newcommand{\bSigma}{\boldsymbol{\Sigma}}
\newcommand{\bvarepsilon}{\boldsymbol{\varepsilon}}

\newcommand{\rmvec}{{\rm{vec}}}
\renewcommand{\exp}{{\rm{exp}}}

\newcommand{\KL}{{\rm{KL}}}

\newcommand{\SVD}{{\rm{SVD}}}

\newcommand{\LR}{{\rm{LR}}}

\newcommand{\argmin}{\mathop{\rm arg\min}}
\newcommand{\argmax}{\mathop{\rm arg\max}}

\newcommand{\tuckerrank}{{\rm Tucrank}}
\newcommand{\F}{{\rm F}}
\newcommand{\grad}{{\rm grad \,}}

\newcommand{\QR}{{\rm QR}}

\newcommand{\RGN}{{\rm RGN}}
\newcommand{\var}{{\textrm Var}}

\newtheorem{Definition}{Definition}
\newtheorem{Theorem}{Theorem}

\newtheorem{Lemma}{Lemma}
\newtheorem{Remark}{Remark}
\newtheorem{Corollary}{Corollary}

\newtheorem{Proposition}{Proposition}

\DeclareMathAlphabet\mathbfcal{OMS}{cmsy}{b}{n}

\endlocaldefs
\usepackage{xr}
\begin{document}
	
	\begin{frontmatter}
		\title{Tensor-on-Tensor Regression: Riemannian Optimization, Over-parameterization, Statistical-computational Gap, and Their Interplay}
		\runtitle{Tensor-on-Tensor Regression}
		
		\begin{aug}
			\author[A]{\fnms{Yuetian} \snm{Luo}\ead[label=e1]{yuetian@uchicago.edu}}
			\and
			\author[B]{\fnms{Anru R.} \snm{Zhang}\ead[label=e2]{anru.zhang@duke.edu}}
			\address[A]{Data Science Institute, University of Chicago, \printead{e1}}
			\address[B]{Department of Biostatistics \& Bioinformatics and Department of Computer Science, Duke University, \printead{e2}}
		\end{aug}
		
		\begin{abstract}
			We study the tensor-on-tensor regression, where the goal is to connect tensor responses to tensor covariates with a low Tucker rank parameter tensor/matrix without prior knowledge of its intrinsic rank. We propose the Riemannian gradient descent (RGD) and Riemannian Gauss-Newton (RGN) methods and cope with the challenge of unknown rank by studying the effect of rank over-parameterization. We provide the first convergence guarantee for the general tensor-on-tensor regression by showing that RGD and RGN respectively converge linearly and quadratically to a statistically optimal estimate in both rank correctly-parameterized and over-parameterized settings. Our theory reveals an intriguing phenomenon: Riemannian optimization methods naturally adapt to over-parameterization without modifications to their implementation. We also prove the statistical-computational gap in scalar-on-tensor regression by a direct low-degree polynomial argument. Our theory demonstrates a ``blessing of statistical-computational gap" phenomenon: in a wide range of scenarios in tensor-on-tensor regression for tensors of order three or higher, the computationally required sample size matches what is needed by moderate rank over-parameterization when considering computationally feasible estimators, while there are no such benefits in the matrix settings. This shows moderate rank over-parameterization is essentially ``cost-free" in terms of sample size in tensor-on-tensor regression of order three or higher. Finally, we conduct simulation studies to show the advantages of our proposed methods and to corroborate our theoretical findings.
		\end{abstract}
		
		\begin{keyword}[class=MSC2020]
			\kwd[Primary ]{62H15}
			\kwd[; secondary ]{62C20}
		\end{keyword}
		
		\begin{keyword}
			\kwd{Tensor-on-tensor regression, over-parameterization, Riemannian optimization, statistical-computational gaps, low-degree polynomials}
		\end{keyword}
		
	\end{frontmatter}
	
	\begin{sloppypar}
		\section{Introduction} \label{sec: introduction}
		The analysis of tensor or multiway array data has emerged as a very active topic of research in statistics, applied mathematics, machine learning, and signal processing \citep{kolda2009tensor}, along with many important applications, such as neuroimaging analysis \citep{zhou2013tensor}, latent variable models \citep{anandkumar2014tensor}, and collaborative filtering \citep{bi2018multilayer}. This paper studies a general class of problems termed {\it tensor-on-tensor regression}, which aims to characterize the relationship between covariates and responses in the form of scalars, vectors, matrices, or high-order tensors:
		\begin{equation}\label{eq:model}
			\bcY_i = \langle \bcA_i, \bcX^* \rangle_* + \bcE_i, \quad i = 1,\ldots, n.
		\end{equation}
		Here, $\bcA_i \in \bbR^{p_1 \times \cdots \times p_d}, i = 1,\ldots,n$ are the known order-$d$ (or $d$-way) tensor covariates. $\bcY_i, \bcE_i \in\mathbb{R}^{p_{d+1} \times \cdots \times p_{d+m}}$ are both order-$m$ tensors and are observations and unknown noise, respectively. $\bcX^* \in \bbR^{p_1 \times \cdots \times p_d \times p_{d+1} \times \cdots \times p_{d+m}}$ is an order-$(d+m)$ tensor parameter of interest. $\langle \cdot, \cdot \rangle_*$ is the contracted tensor inner product defined as $\langle \bcA_i, \bcX^* \rangle_* \in \bbR^{p_{d+1} \times \cdots \times p_{d+m}}$,
		\begin{equation*}
			\left(\langle \bcA_i, \bcX^* \rangle_* \right)_{[j_1,\ldots, j_m]} = \sum_{\substack{k_l=1,\\l=1,\ldots,d}}^{p_l}  \bcA_{i[k_1,\ldots,k_d]}  \bcX^*_{[k_1,\ldots,k_d,j_1,\ldots,j_m]}.
		\end{equation*}
		Throughout the paper, we consider $d$ and $m$ to be fixed constants. We also stack all responses and errors to $\bcY$, $\bcE \in \bbR^{n \times p_{d+1} \times \cdots \times p_{d+m}}$, where $\bcY_{[i,:,\ldots,:]} = \bcY_i$ and $\bcE_{[i,:,\ldots,:]} = \bcE_i$. 
		Then the tensor-on-tensor regression model can be written succinctly as 
		$\bcY = \scA(\bcX^*) + \bcE$, where $\scA: \bbR^{p_1\times \cdots \times p_{d+m}} \to \bbR^{n \times p_{d+1} \times \cdots \times p_{d+m}}$ is a linear map such that
		\begin{equation}\label{eq: affine operator}
			\scA(\bcX^*)_{[i,:,\ldots,:]} = \langle \bcA_i, \bcX^* \rangle_* \quad \text{ for } \quad i = 1,\ldots, n.
		\end{equation}
		Our goal is to estimate $\bcX^*$ based on $(\bcY, \scA)$. 
		
		Tensor-on-tensor regression model was proposed and studied in \cite{raskutti2015convex,lock2018tensor}. The generic tensor-on-tensor regression covers many special tensor regression models in the literature, such as
		\begin{itemize}
			\item Scalar-on-tensor regression \citep{zhou2013tensor,mu2014square}: $m =0$; 
			\item Tensor-on-vector regression \citep{li2016parsimonious,sun2017store}: $d = 1$; 
			\item Scalar-on-matrix regression (or matrix trace regression) \citep{recht2010guaranteed}: $m = 0, d = 2$. 
		\end{itemize}
		There is a great surge of interest in tensor-on-tensor regression for its applications \citep{lock2018tensor,gahrooei2021multiple,llosa2022reduced}. Specific examples include:
		\begin{itemize}[leftmargin=*]
			\item {\it Neuroimaging Data Analysis}. 
			Studies in neuroscience are greatly facilitated by a variety of neuroimaging technologies. Tensor-on-tensor regression provides interpretable analysis of such datasets \citep{zhou2013tensor,li2016parsimonious}. For example, tensor-on-vector regression has been applied to compare MRI scans across different autism spectrum disorder groups \citep{sun2017store}, which has helped evaluate the effectiveness of a potential drug. Scalar-on-tensor regression has been used to predict neurological diseases, such as attention deficit hyperactivity disorder, and reveal regions of interest in the brain that affect the progression of diseases \citep{zhou2013tensor}.
			\item {\it Facial Image Data Analysis}. Attributes prediction from facial images is popular in social data analysis. Oftentimes, each facial image is labeled only with the name of the individual, often a celebrity, while people are interested in inferring more features from that. Tensor-on-tensor regression and tensor-variate analysis of variance have been proposed to predict describable attributes from a facial image \citep{lock2018tensor} and distinguish facial characteristics related to ethnic origin, age group, and gender \citep{llosa2022reduced}.
			\item {\it Longitudinal Relational Data Analysis}. Longitudinal relational data among a set of objects can be represented as a time series of matrices, where each entry of the matrices represents a directed relationship involving pairs of objects at a given time. The relation between one pair of objects may have an effect on the relation between members of another pair, an effective tensor-on-tensor regression model has been developed to estimate such effects \citep{hoff2015multilinear}.
		\end{itemize}
		
		Meanwhile, tensor datasets are often high-dimensional, i.e., the ambient data dimension is substantially bigger than the sample size. It is thus crucial to exploit the hidden low-dimensional structures from the datasets to facilitate the follow-up analyses. In tensor data analysis, low-rankness is among the most commonly considered structural assumptions. In this paper, we assume the target parameter $\bcX^*$ has an intrinsic low Tucker (or multilinear) rank $\br^* = (r^*_1,\ldots, r^*_d,r^*_{d+1},\ldots,r^*_{d+m})$, i.e., all fibers\footnote{Fibers are bar-shaped vectors and are counterpart of matrix columns and rows in a tensors \citep{kolda2009tensor}.} of $\bcX^*$ along mode-$k$ lie in a $r_k^*$ dimensional subspace of $\mathbb{R}^{p_k}$ for $k=1,\ldots,d+m$. 
		
		\subsection{Central Questions}
		A natural question on low-rank tensor-on-tensor regression is
		\vskip.05cm
		\noindent \fbox{ \parbox{0.98\textwidth}{
				{\it 1. Can we develop fast and statistically optimal solutions for the general low-rank tensor-on-tensor regression?
				}
		}}
		\vskip.05cm
		
		Various algorithms were proposed in the literature to solve specific instances of tensor-on-tensor regression with provable guarantees, such as variants of gradient descent methods \citep{rauhut2017low,yu2016learning,chen2016non,ahmed2020tensor,han2020optimal,hao2018sparse,tong2021scaling}, alternating minimization \citep{zhou2013tensor}, Bayesian Markov chain Monte Carlo \citep{guhaniyogi2015bayesian}, and Riemannian optimization methods \citep{kressner2016preconditioned,luo2021low} for scalar-on-tensor regression; regularized rank constrained least squares \citep{rabusseau2016low}, alternating minimization \citep{sun2017store} and envelope method \citep{li2016parsimonious} for tensor-on-vector regression. The theoretical guarantees of these methods were developed case-by-case under the assumption that the intrinsic tensor rank is known. In addition, \cite{hoff2015multilinear} proposed a Bayesian approach to solve the tensor-on-tensor regression when the mode numbers of the predictor and the response are equal. \cite{lock2018tensor,liu2020low} proposed alternating least squares procedures for solving the general tensor-on-tensor regression, and a numerical study on the effect of rank misspecification was performed in \cite{lock2018tensor} without theoretical exploration. Asymptotic analysis for the computationally intensive maximum likelihood estimator is provided in \cite{llosa2022reduced} for different low-rank tensor formats with known intrinsic ranks. The convex relaxation methods for tensor-on-tensor regression, including the computationally infeasible tensor nuclear norm relaxation, were studied in \cite{raskutti2015convex}. In summary, despite a great amount of effort in the literature, a general, fast, and statistically optimal framework for tensor-on-tensor regression is still underdeveloped.
		
		Moreover, the intrinsic rank $\br^*$ is usually unknown in practice, while tuning rank is even more challenging for tensors than matrices as $(d+m)$ parameter values need to be tuned simultaneously. Thus, an important question is:
		\vskip.05cm
		\noindent \fbox{ \parbox{0.98\textwidth}{
				{\it 2. Can we solve tensor-on-tensor regression robustly without knowing the intrinsic rank?}
		}}
		\vskip.05cm
		
		To this end, we adopt a rank over-parameterization scheme: we introduce a conservative guess of rank $\br:= (r_1,\ldots,r_{d+m})\geq (r_1^*,\ldots, r_{d+m}^*)$ and solve the following tensor-on-tensor regression under the possibly over-parameterized regime:
		\begin{equation}\label{eq:minimization}
			\begin{split}
				\widehat{\bcX}_{\textrm{opt}} =  &\argmin_{\bcX \in \mathbb{R}^{p_1\times \cdots \times p_{d+m}}} f(\bcX) := \frac{1}{2} \left\|\bcY - \scA(\bcX)\right\|_\F^2,\\
				&\text{subject to}\quad \tuckerrank(\bcX)\leq  \br.
			\end{split}
		\end{equation}
		Here, $\tuckerrank(\bcX)$ is the Tucker rank of $\bcX$ (see formal definition in the {\bf Notation and Preliminaries} Section). In most of the aforementioned literature, the ranks were assumed to be correctly specified and the results do not directly apply to the possibly over-parameterized scenario in \eqref{eq:minimization}. We will illustrate later that Riemannian optimization is an ideal scheme to treat rank-constrained optimization like \eqref{eq:minimization}. However, under the over-parameterized regime, the classic convergence theory of Riemannian optimization does not apply since the true parameter $\bcX^*$ is merely a boundary point of the Riemannian manifold consisting of tensors with incorrectly specified rank. 
		
		In addition, tensor problems often exhibit statistical-computational gaps \citep{hillar2013most,richard2014statistical}. For example, in scalar-on-tensor regression, i.e., $m =0$, and suppose $p_1 = \cdots = p_{d} = p$ and $r_1^* = \cdots = r_{d}^* = r^*$ is known and the design is Gaussian ensemble (to be formally introduced in Section \ref{sec:theory}), it has been shown that rank minimization recovers $\bcX^*$ with $\Omega(pr^* + r^{*d})$ samples \citep{mu2014square}; but the rank minimization is generally NP-hard to compute \citep{hillar2013most}. On the other hand, all existing polynomial-time algorithms require at least $\Omega(p^{d/2} r^* + r^{*{d}})$ samples to guarantee recovery \citep{han2020optimal}. So when $d \geq 3$, there exists a significant gap on the sample complexities between what can be achieved information theoretically and by existing polynomial-time algorithms. \cite{xia2022inference} leveraged this hypothetical gap to claim there is no need to debias in scalar-on-tensor regression inference. Intriguingly, this gap seems to close when $d = 2$, i.e., in the matrix case, since $p^{d/2} r^* + r^{*{d}} = pr^* + r^{*{d}}$. So we ask:
		\vskip.05cm
		\noindent \fbox{ \parbox{0.98\textwidth}{
				{\it 3. Is there a statistical-computational gap in tensor-on-tensor regression? What is the difference between tensor and matrix settings?}
		}}
		\vskip.05cm
		
		In the era of big data, Riemannian optimization and over-parameterization have become a common remedy for nonconvexity in high-dimensional statistics and machine learning, where the statistical-computational gap is a prevalent phenomenon. As these ingredients nicely gather in tensor-on-tensor regression, a more open-ended question is
		\vskip.05cm
		\noindent \fbox{ \parbox{0.98\textwidth}{
				{\it 4. Is there any interplay among Riemannian optimization, over-parameterization, and statistical-computational gap?}
		}}

		\subsection{Our Contributions} \label{sec: contributions}
		
		We aim to answer the four questions above. Our specific contributions include:
		
		\vskip.05cm
		
		\noindent{\bf (Over-parameterization, algorithms, convergence theory, and statistical optimality)} We address the unknown intrinsic rank through the rank over-parameterization scheme in \eqref{eq:minimization}. We introduce the Riemannian gradient descent (RGD) and Riemannian Gauss-Newton (RGN) algorithms for tensor-on-tensor regression and develop the corresponding convergence guarantees. We specifically show with proper initialization, RGD and RGN respectively converge linearly and quadratically to the true parameter $\bcX^*$ up to some statistical error. Especially in the noiseless setting, i.e., $\bcE = 0$, RGD and RGN respectively converge linearly and quadratically to the exact parameter $\bcX^*$. Our convergence theory for over-parameterized Riemannian optimization algorithms is novel, covers the rank under-parameterized cases as well, and cannot be inferred from the standard convergence theories in the Riemannian optimization literature, since the true parameter $\bcX^*$ only lies on the boundary of the working Riemannian manifold consisting of tensors with incorrectly specified rank. We further show the estimation error achieved by RGD and RGN matches the minimax risk lower bound under the Gaussian ensemble design. To our best knowledge, this is the first algorithmic convergence result for tensor-on-tensor regression with optimal statistical error guarantees. In the specific over-parameterized matrix trace regression setting, our results yield the first linear/quadratic convergence guarantee for RGD/RGN. Compared to the existing results on factorized GD in the over-parameterized matrix trace regression \citep{zhuo2021computational,zhang2021preconditioned}, our second-order algorithm RGN and the corresponding theory are novel, which improve the results in literature in many ways. 
		
		Our convergence theory reveals an intriguing phenomenon: in tensor-on-tensor regression, {\bf Riemannian optimization algorithms adapt to over-parameterized scenarios without modifications.} This is significantly different from the classic factorized gradient descent algorithm where preconditioning is needed. Table \ref{tab: comparison table} compares our results with the existing ones on over-parameterized matrix trace regression. 
		
		\begin{table}[h]
			\centering
			\begin{tabular}{c |c | c | c | c }
				\hline
				\multicolumn{5}{c}{Over-parameterized Matrix Trace Regression} \\
				\hline
				\multirow{2}{3em}{Algorithm} & statistical & convergence & require  & parameter   \\
				& error rate & rate &  tuning & matrix type  \\
				\hline
				RGD &  \multirow{2}{3em}{optimal} & \multirow{2}{3em}{linear} & \multirow{2}{2em}{no}  & \multirow{2}{3em}{general} \\
				(this work) & & & &\\
				\hline
				RGN  &  \multirow{2}{3em}{optimal} & \multirow{2}{3em}{quadratic} & \multirow{2}{2em}{no}   & \multirow{2}{3em}{general} \\
				(this work) & & & &\\
				\hline
				Factorized GD &  \multirow{2}{3em}{optimal} & \multirow{2}{3em}{sublinear} & \multirow{2}{2em}{yes}  & \multirow{2}{2em}{PSD} \\
				\citep{zhuo2021computational} & & & & \\
				\hline
				Preconditioned &  \multirow{2}{3.9em}{suboptimal} & \multirow{2}{3em}{linear} & \multirow{2}{2em}{yes} & \multirow{2}{2em}{PSD}  \\
				Factorized GD \citep{zhang2021preconditioned} & & & &\\
				\hline
			\end{tabular}
			\caption{Riemannian gradient descent (RGD), Riemannian Gauss-Newton (RGN) versus factorized gradient descent (Factorized GD), preconditioned factorized GD for over-parameterized matrix trace regression.  
			}\label{tab: comparison table}
		\end{table}
		
		Although developing proper initialization for all cases of tensor-on-tensor regression is difficult, we introduce spectral methods that yield adequate initializations for both RGD and RGN in four prominent instances, {\it scalar-on-tensor regression}, {\it tensor-on-vector regression}, {\it matrix trace regression} and {\it rank-$1$ tensor-on-tensor regression} under Gaussian ensemble design. 
		
		\vskip.05cm
		
		\noindent{\bf (Statistical-computational gap and sample size requirement)} In this paper, we establish rigorous evidence on the statistical-computational gap in scalar-on-tensor regression via low-degree polynomials methods. Our argument shows $n = \Omega(p^{d/2})$ samples are necessary for any polynomial-time method to succeed. Existing hardness evidence from low-degree polynomials is often established for statistical problems with the simple ``signal+noise" structure. Such a structure enables the decoupling of signal and noise that simplifies the analysis. To our best knowledge, our low-degree hardness evidence is the first one for problems with complex correlated structures.
		
\begin{figure}[h!]
			\centering
			\includegraphics[width=0.8\textwidth]{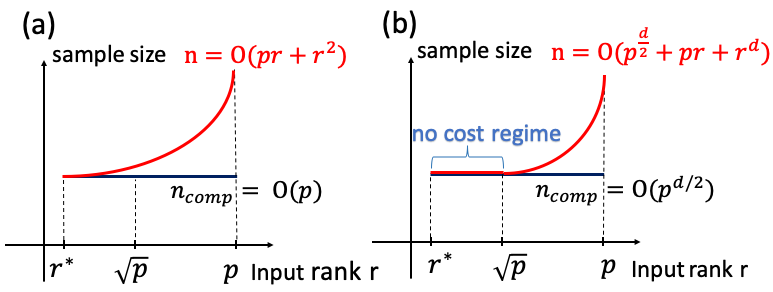}
			\caption{Comparison of sample size requirements in over-parameterized matrix trace (Panel (a)) and scalar-on-tensor regressions (Panel (b)) under Gaussian ensemble design. Here the red line denotes the sample size ($n$) requirements for the RGD and RGN to succeed with input rank $r$ and spectral initialization and the black line ($n_{\text{comp}}$) is the sample complexity of the computational limit, i.e. the minimum sample size requirement for any efficient algorithms. For simplicity, we assume $p_1 = \ldots = p_d = p$, $r_1 = \ldots = r_d = r$, $r^*_1 = \ldots = r^*_d = r^*$, $d$ and $r^*$ are some fixed constants, $\bcE = 0$ and $\bcX^*$ is well-conditioned. }\label{fig: tensor-matrix-comparison}
		\end{figure}

		Based on the computational lower bounds and algorithmic upper bounds developed in this paper, we draw Figure \ref{fig: tensor-matrix-comparison} to illustrate the sample size requirements in over-parameterized matrix trace regression with $d=2$ (Panel (a)) and scalar-on-tensor regression, a prominent instance of tensor-on-tensor regression, with $d\geq 3$ (Panel (b)). When the input rank $r$ is greater than $\sqrt{p}$, i.e., in the heavily over-parameterized regime, we show that an extra sample complexity is needed for RGD and RGN to converge in both regressions. When the input rank $r$ is between $r^*$ and $\sqrt{p}$, i.e., in the moderately over-parameterized regime, extra sample complexity is still required in matrix trace regression (Figure \ref{fig: tensor-matrix-comparison}(a)). On the other hand, in scalar-on-tensor regression (Figure \ref{fig: tensor-matrix-comparison}(b)), no larger sample size is required to account for the inflated input rank, as the red line is flat in the ``no extra cost" regime in Figure \ref{fig: tensor-matrix-comparison}(b). 
		
		This alludes to an important message, {\bf moderate rank over-parameterization is cost-free in terms of sample size for a computationally feasible optimal estimator in scalar-on-tensor regression.} The computational barrier, although being a tough scenario and is often referred to as the ``curse of computability," becomes a {\it ``blessing"} to over-parameterization here, as no extra samples are required if this large but essential sample size condition is met to guarantee that the computationally feasible estimator is achievable!
		
		\vskip.05cm
		
		\noindent{\bf (New technical tools)} We introduce a series of technical tools for theory development in this paper, including a tangent space projection error bound, a tensor decomposition perturbation bound under the over-parameterized setting, and a simple formula for computing expected values of Hermite polynomials on correlated multivariate Gaussian random variables while developing low-degree polynomials lower bounds. 
		See Section \ref{sec:technical contribution} for a summary of our technical contributions.
		
		\vskip.05cm
		
		\noindent{\bf (Implementation details and numerical experiments)} Finally, we discuss the implementation details of RGD and RGN for tensor-on-tensor regression in Section \ref{sec:implementation}. We specifically find a reduction from computing RGN update to solving $(m+1)$ separate least squares. This reduction yields a fast implementation of RGN. We conduct numerical studies to show the convergence and required sample size of our proposed algorithms match our theoretical findings. We also compare the numerical performance of our algorithms with existing ones. The results show the proposed algorithms have significant advantages in both rank correctly-specified and overspecified tensor-on-tensor regression.
		
		\subsection{Related Prior Work} \label{sec: literature review}
		
		This work is related to several lines of research on over-parameterization, Riemannian optimization, and computational barriers in tensor problems. 
		
		First, over-parameterization has attracted much attention in modern data science due to the great success of deep learning. The concept of over-parameterization generally refers to the scenario when learning problems include more model parameters than necessary. Recent studies show that over-parameterization brings both computational and statistical benefits when solving complex problems \citep{soltanolkotabi2018theoretical,bartlett2020benign,belkin2019reconciling}. There is a vast amount of literature on studying the role of over-parameterization to demystify deep learning \citep{bartlett2021deep,belkin2021fit}. This paper focuses on the effect of over-parameterization specifically in the rank-constrained tensor-on-tensor regression problem. In particular, we consider a special type of over-parameterization where the input rank to the model is overspecified.
		
		Second, Riemannian manifold optimization methods have been powerful in solving optimization problems with geometric constraints \citep{absil2009optimization}. Many progress in this topic were made for the low-rank matrix estimation \citep{keshavan2009matrix,boumal2011rtrmc,wei2016guarantees,meyer2011linear,mishra2014fixed,vandereycken2013low,huang2018blind,luo2020recursive,hou2020fast}. Moreover, Riemannian manifold optimization methods under various Riemannian geometries have been explored in many tensor problems, such as tensor decomposition \citep{elden2009newton,savas2010quasi,ishteva2009differential,breiding2018riemannian}, scalar-on-tensor regression \citep{kressner2016preconditioned,luo2021low}, tensor completion \citep{kasai2016low,dong2021new,kressner2014low,heidel2018riemannian,xia2017polynomial,steinlechner2016riemannian,wang2021entrywise,cai2021provable}, and robust tensor PCA \citep{cai2021generalized}.
		
		Third, many high-dimensional tensor problems exhibit the statistical-computational gaps, i.e. the gap between different signal-to-noise ratio thresholds that make the problem information-theoretically solvable versus polynomial-time solvable. Rigorous evidence for such gaps has been provided to tensor completion \citep{barak2016noisy}, tensor PCA/SVD \citep{zhang2018tensor,brennan2020reducibility,dudeja2020statistical,choo2021complexity}, tensor clustering \citep{luo2020tensor,han2020exact} and tensor-on-tensor association detection \citep{diakonikolas2022statistical}. This work provides a rigorous piece of evidence for the statistical-computational gap in scalar-on-tensor regression under the low-degree polynomials framework. 
		
		Finally, a special case of our setting, over-parameterized matrix trace regression, has attracted much attention recently. The results along this line include two categories: (1) $r \geq r^*$ and $n = O((p_1 + p_2)r)$: the problem is over-parameterized and identifiable \cite{zhuo2021computational,zhang2021preconditioned,ding2021rank}; (2) $r \geq r^*$ and $n = O((p_1 + p_2)r^*)$: as the sample size is smaller than the number of free parameters in the model, there can be infinitely many solutions to \eqref{eq:minimization} and the model is unidentifiable. One important finding in Category (2) is that with small magnitude initialization, vanilla gradient descent under the factorization formulation tends to implicitly bias towards a low-rank solution \citep{gunasekar2017implicit,li2018algorithmic,li2020towards,fan2022understanding,stoger2021small,ma2022global,jiang2022algorithmic}. Our work provides a unified simple Riemannian optimization framework to solve the general tensor-on-tensor regression problem under the setting in Category (1). The implication of our results in over-parameterized matrix trace regression is further discussed in Remarks \ref{rm: comparision-in-matrix} and \ref{rm: matrix-tensor-sample-compare}. 
		
		\subsection{Organization of the Paper} \label{sec: organization}
		After a brief introduction of notation and preliminaries in Section \ref{sec:notations}, we introduce our main algorithms, Riemannian gradient descent and Riemannian Gauss-Newton in Section \ref{sec: Algorithm}. The convergence results of RGD and RGN in the general tensor-on-tensor regression and applications in specific examples are discussed in Sections \ref{sec:theory} and \ref{sec: applications-initializations}, respectively. Computational limits are discussed in Section \ref{sec: computational-limit}. Technical contributions are summarized in Section \ref{sec:technical contribution}. 
		Implementation details of RGD/RGN and numerical studies are presented in Sections \ref{sec:implementation} and \ref{sec:numerics}, respectively.
		Conclusion and future work are given in Section \ref{sec: conclusion}. Additional algorithms, numerical studies and all technical proofs are collected in Supplements \ref{sec: HOSVD, STHOSVD}-\ref{sec: additional lemmas}.

		\subsection{Notation and Preliminaries}\label{sec:notations}
		
		Let $[r] = \{1,\ldots,r\}$ for any positive integer $r$. Lowercase letters (e.g., $a$), lowercase boldface letters (e.g., $\u$), uppercase boldface letters (e.g., $\U$), and boldface calligraphic letters (e.g., $\bcA$) denote scalars, vectors, matrices, and order-3-or-higher tensors, respectively. 
		We use bracket subscripts to denote sub-vectors, sub-matrices, and sub-tensors. 
		For any matrix $\D \in \mathbb{R}^{p_1\times p_2}$, let $\sigma_k(\D)$ be the $k$th largest singular value of $\D$. We also denote $\SVD_r(\D) = [\u_1 ~ \cdots \u_r]$ and QR($\D$) as the subspace composed of the leading $r$ left singular vectors and the Q part of the QR decomposition of $\D$, respectively. 
		$\I_r$ represents the $r$-by-$r$ identity matrix. Let $\mathbb{O}_{p, r} = \{\U \in \bbR^{p \times r}: \U^\top \U=\I_r\}$ and for any $\U\in \mathbb{O}_{p, r}$, denote $P_{\U} = \U\U^\top$. 
		The matricization operation $\mathcal{M}_k(\cdot)$ unfolds an order-$d$ tensor along mode $k$ to a matrix, say $\bcA\in\mathbb{R}^{p_1\times \cdots \times p_d}$ to $\mathcal{M}_k(\bcA)\in \mathbb{R}^{p_k\times p_{-k}}$, where $p_{-k} = \prod_{j\neq k}p_j$ and its detailed definition is provided in Supplement \ref{sec: additional-notation}. The Frobenius norm of tensor $\bcA$ is defined as $\|\bcA\|_{\F} = \left(\sum_{i_1,\ldots, i_d} \bcA_{[i_1,\ldots, i_d]}^2\right)^{1/2}$. The Tucker rank of an order-$d$ tensor $\bcA$, denoted by $\tuckerrank(\bcA)$, is defined as a $d$-tuple $\br := (r_1, \ldots, r_d)$, where $r_k = \text{rank}(\mathcal{M}_k(\bcA))$. Any Tucker rank-$(r_1,\ldots,r_d)$ tensor $\bcA$ admits the following Tucker decomposition \citep{tucker1966some}: $	\bcA = \llbracket \bcS; \U_1, \ldots, \U_d\rrbracket := \bcS\times_1 \U_1 \times \cdots \times_d \U_d,$
		where $\bcS \in \bbR^{r_1 \times \cdots \times r_d}$ is the core tensor and $\U_k = \SVD_{r_k}(\cM_k(\bcA))$ is the mode-$k$ top $r_k$ left singular vectors. Here, the mode-$k$ product of $\bcA \in \mathbb{R}^{p_1 \times \cdots \times p_d}$ with a matrix $\B\in \mathbb{R}^{r_k\times p_k}$, denoted by $\bcA \times_k \B$, is a $p_1 \times \cdots \times p_{k-1}\times r_k \times p_{k+1}\times \cdots \times p_d$-dimensional tensor, and its definition is provided in Supplement \ref{sec: additional-notation}. The following abbreviations are used to denote the tensor-matrix product along multiple modes: $\bcA \times_{k=1}^d \U_k := \bcA \times_{1} \U_1 \times \cdots \times_{d} \U_d$; $\bcA \times_{l\neq k}\U_l := \bcA \times_{1} \U_{1} \times \cdots \times_{k-1} \U_{k-1} \times_{k+1} \U_{k+1} \times \cdots \times_{d} \U_{d}$. For any order-$d$ tensor $\bcZ \in \bbR^{p_1 \times \cdots \times p_d}$ and a $d$-tuple $\br = (r_1,\ldots,r_d)$, let $\bcZ_{\max(\br)} := \bcZ \times_{k=1}^d P_{\widehat{\U}_k}$ be the best Tucker rank $\br$ approximation of $\bcZ$ in terms of Frobenius norm, where $(\widehat{\U}_1, \ldots, \widehat{\U}_d)$ is the solution to $\argmax_{\U_k \in \bbO_{p_k, r_k}, k=1,\ldots,d} \|\bcZ \times_{k=1}^d P_{\U_k} \|_{\F}$ \cite[Theorem 4.2]{de2000best}. Throughout the paper, let 
		$c(d)$ be a constant that depends on $d$ only, whose actual value varies from line to line; $c_1(m), c_2(d,m)$ are noted similarly. Finally, we denote $\scA^*$ as the adjoint of the linear operator $\scA$. 
		
		\section{Riemannian Optimization for Tensor-on-Tensor Regression}\label{sec: Algorithm}
\begin{figure}
			\centering
			\includegraphics[width=0.5\textwidth]{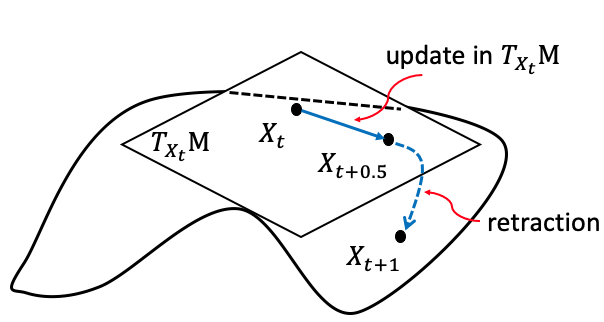}
			\caption{Pictorial illustration of steps in Riemannian optimization}
			\label{fig: Riemannian-opt-ill}
		\end{figure}

		Riemannian optimization concerns optimizing a real-valued function $f$ whose domain is a Riemannian manifold $\bbM$ \citep{absil2009optimization}. The continuous optimization on the Riemannian manifold often requires calculations on the tangent space due to its common non-linearity. A typical procedure of a Riemannian optimization method includes three steps per iteration: 1. find the tangent space of $\bbM$; 2. update the point on the tangent space; 3. map the point from the tangent space back to the manifold, i.e., retraction. A pictorial illustration for the three steps in  Riemannian optimization is presented in Figure \ref{fig: Riemannian-opt-ill}. The readers are also referred to \cite{absil2009optimization} and \cite{boumal2020introduction} for more discussions on Riemannian optimization.  

		\subsection{Geometry of Low Tucker Rank Tensor Manifolds}\label{sec: tensor manifold geometry}

		Denote the collection of $(p_1,\ldots, p_d, p_{d+1},\ldots, p_{d+m})$-dimensional tensors of Tucker rank $\br:= (r_1,\ldots,r_d,r_{d+1},\ldots,r_{d+m})$ by $\bbM_{\br}=\{\bcX \in \bbR^{p_1 \times \cdots \times p_{d+m}}, \tuckerrank (\bcX) = \br \}$. Then $\bbM_{\br}$ forms a $\left\{\prod_{j=1}^{d+m} r_j + \sum_{j=1}^{d+m} r_j(p_j - r_j)\right\}$-dimensional smooth submanifold embedded in $\bbR^{p_1 \times \cdots \times p_{d+m}}$ \citep{uschmajew2013geometry}.  Recall in the general over-parameterized scenario, $\br$ may be different from $\br^*$, the actual rank of the tensor of interest. Suppose $\bcX \in \bbM_{\br}$ has Tucker decomposition $\llbracket \bcS; \U_1,\ldots, \U_d, \U_{d+1}, \ldots, \U_{d+m} \rrbracket$. Define $\V_k=\QR(\cM_k(\bcS)^\top)$, which corresponds to the row space of $\cM_k(\bcS)$, and for $k = 1,\ldots, d+m$, define
		\begin{equation} \label{def: Wk}
			\W_k := \left(\U_{d+m} \otimes \cdots \otimes \U_{k+1} \otimes \U_{k-1} \otimes\cdots \otimes \U_1\right) \V_k \in \bbO_{p_{-k}, r_{k}}, 
		\end{equation}
		where $p_{-k} = \prod^{d+m}_{j=1,j\neq k} p_j$. By the tensor matricization formula provided in Supplement \ref{sec: additional-notation}, 
		$\U_k,\W_k$ correspond to the subspaces of the column and row spans of $\cM_k(\bcX)$, respectively. \cite{koch2010dynamical} provided the explicit formulas for the tangent space of $\bbM_{\br}$ at $\bcX$, denoted by $T_{\bcX} \bbM_{\br}$ (see Supplement \ref{sec: additional-notation} for the expression). 
		We equip $\bbM_{\br}$ with the Riemannian metric induced by the natural Euclidean inner product $\langle \cdot,\cdot \rangle$. Under this metric, the following operator $P_{T_{\bcX}}$ projects any tensor $\bcZ\in \bbR^{p_1 \times \cdots \times p_{d+m}}$ onto the tangent space $T_{\bcX} \bbM_{\br}$,
		\begin{equation} \label{eq: tangent space projector}
			P_{T_{\bcX}} (\bcZ) := \bcZ \times_{k=1}^{d+m} P_{\U_k} + \sum_{k=1}^{d+m} \cT_k(P_{\U_{k\perp}} \cM_k(\bcZ) P_{\W_k}  ),
		\end{equation} where $\cT_k(\cdot)$ denotes the mode-$k$ tensorization, i.e., the reverse operator of $\cM_k(\cdot)$.

		\subsection{Riemannian Gradient Descent and Gauss-Newton for Tensor-on-Tensor Regression}\label{sec: RGD-and-RGN}

		The Riemannian gradient of a smooth function $f:\bbM_{\br} \to \bbR$ at $\bcX\in \bbM_{\br}$ is defined as the unique tangent vector ${\rm grad}\, f(\bcX) \in T_\bcX \bbM_{\br}$ such that $\langle {\rm grad}\, f(\bcX), \bcZ \rangle = {\rm D}\, f(\bcX)[\bcZ], \forall\, \bcZ \in T_\bcX \bbM_{\br},$ where ${\rm D}f(\bcX)[\bcZ]$ denotes the directional derivative of $f$ at point $\bcX$ along direction $\bcZ$. We can calculate the Riemannian gradient for the tensor-on-tensor regression as follows.
		\begin{Lemma}[Riemannian gradient]\label{lm: gradient} For $f(\bcX)$ in \eqref{eq:minimization}, ${\rm grad}\, f(\bcX) = P_{T_\bcX}(\scA^*(\scA(\bcX) - \bcY)),$ where $\scA^*$ is the adjoint operator of $\scA$. 
		\end{Lemma}
		By Lemma \ref{lm: gradient}, a natural idea of RGD update is
		$\bcX^{t+0.5} = \bcX^t - \alpha_t P_{T_{\bcX^t}} \scA^*( \scA(\bcX^t) - \bcY)$,
		where the stepsize $\alpha_t$ is chosen as the local steepest descent direction with a closed form as
		\begin{equation} \label{eq: stepsize-choice}
			\begin{split}
				\alpha_t := &  \argmin_{\alpha \in \bbR} \frac{1}{2}  \left\|\bcY - \scA\left( \bcX^t - \alpha P_{T_{\bcX^t}} \scA^*( \scA(\bcX^t) - \bcY )  \right) \right\|_\F^2\\
				= & \frac{ \|P_{T_{\bcX^t}} (\scA^*( \scA(\bcX^t) - \bcY ))\|_\F^2}{\|\scA P_{T_{\bcX^t}} (\scA^*( \scA(\bcX^t) - \bcY ))\|_{\F}^2}.
			\end{split}
		\end{equation}
		As illustrated in Figure \ref{fig: Riemannian-opt-ill}, the updated iterate $\bcX^{t+0.5}$ may not be on the Riemannian manifold $\mathbb{M}_{\br}$. We can apply two types of computationally efficient retractions to bring $\bcX^{t+0.5}$ back to $\mathbb{M}_{\br}$: truncated high-order singular value decomposition (T-HOSVD) \citep{de2000multilinear} or sequentially truncated high-order singular value decomposition (ST-HOSVD) \citep{vannieuwenhoven2012new}. The pseudocode of T-HOSVD and ST-HOSVD are given in Algorithms \ref{alg: t-HOSVD} and \ref{alg: st-HOSVD} in Supplement \ref{sec: HOSVD, STHOSVD}, respectively.

		Moreover, the first-order methods, such as RGD described above, can suffer from slow convergence and low precision in large-scale settings. A natural remedy is to apply second-order methods, such as the Newton algorithm. For tensor-on-tensor regression, the Riemannian Newton relies on the construction and inversion of Riemannian Hessian, which is analytically difficult to develop and computationally intensive. Alternatively, the following Riemannian Gauss-Newton update is a nice approximation of the Riemannian Newton for the nonlinear least squares objective \citep[Section 8.4.1]{absil2009optimization}: 
		\begin{equation}\label{eq:RGN-original}
			-\grad f(\bcX^t)= P_{T_{\bcX^t}} \left( \scA^* (\scA (\eta)) \right), \text{ where } \eta \in T_{\bcX^t} \bbM_{\br}.
		\end{equation} 
		Gauss-Newton has a similar per-iteration complexity as first-order methods but requires much fewer iterations to converge in several other tensor decomposition problems \citep{sorber2013optimization}. The direct calculation of \eqref{eq:RGN-original} is still complicated. Surprisingly, we can show the Gauss-Newton equation \eqref{eq:RGN-original} for tensor-on-tenor regression is equivalent to the following least squares equation.
		\begin{Lemma}\label{lm: gauss-newton}
			For $f(\bcX)$ in \eqref{eq:minimization}, suppose the current iterate is $\bcX^t$. Then the Riemannian Gauss-Newton update is $\eta^{\RGN} = \argmin_{\eta \in T_{\bcX^t} \bbM_{\br} } \frac{1}{2}\|\bcY - \scA P_{T_{\bcX^t}}(\bcX^t + \eta)\|_\F^2$.
		\end{Lemma}
		As we will discuss in Section \ref{sec:implementation} that under some mild condition on $\scA$, the least squares problem in \eqref{eq: alg least square} has a unique solution and can be implemented and solved efficiently via solving $(m+1)$ separate least squares based on Lemma \ref{lm: gauss-newton}. The pseudocode of the overall RGD and RGN procedures are summarized in Algorithm \ref{alg:RGD-RGN}.

		\begin{algorithm}[h]
			\caption{Riemannian Gradient Descent/Gauss-Newton for (Over-parameterized) Tensor-on-Tensor Regression}
			\begin{algorithmic}[1]
			\State \noindent {\bf Input}: $\bcY\in \mathbb{R}^{n \times p_{d+1} \times \cdots \times p_{d+m}}, \bcA_1,\ldots, \bcA_n \in \mathbb{R}^{p_1\times \cdots \times p_d}$, $t_{\max}$, input Tucker rank $\br$, and initialization $\bcX^0$ of Tucker rank $\br$.
				\For{$t=0, 1, \ldots, t_{\max}-1$}
				\State (RGD Update) Compute $\bcX^{t+0.5} = \bcX^t - \alpha_t P_{T_{\bcX^t}} \scA^*( \scA(\bcX^t) - \bcY )$, where $\alpha_t$ is given in \eqref{eq: stepsize-choice}. 
				
				(RGN Update) Solve the least squares problem
				\begin{equation} \label{eq: alg least square}
					\bcX^{t+0.5} = \argmin_{\bcX \in T_{\bcX^t} \bbM_{\br} } \frac{1}{2}\|\bcY - \scA P_{T_{\bcX^t}}(\bcX)\|_2^2.
				\end{equation}
				\State Update $\bcX^{t+1} =\cH_{\br} \left( \bcX^{t+0.5} \right)$. Here $\cH_{\br} (\cdot)$ is the retraction map onto $\bbM_{\br}$, e.g., ST-HOSVD and T-HOSVD.
				\EndFor
				\State {\bf Output}: $\bcX^{t_{\max}}$.
			\end{algorithmic}
			\label{alg:RGD-RGN}
		\end{algorithm}

		\begin{Remark}[Riemannian Optimization for Bounded Rank Constraint]\label{rmk: bounded-rank}
			The classic RGD/RGN methods are designed to optimize on smooth manifolds. This corresponds to minimizing the objective function in \eqref{eq:minimization} with the fixed Tucker rank constraint $\tuckerrank (\bcX) = \br$ since $\bbM_\br$ is a smooth manifold. Note that $\{\bcX:\tuckerrank(\bcX) \leq \br\}$ is not a smooth manifold while such bounded rank constraint is essential to handle over-parameterization, the classic theory no longer applies. Regardless, we propose to continue using Algorithm \ref{alg:RGD-RGN} even with the bounded rank constraint. 
		\end{Remark}

		\section{Theory of RGD/RGN in Tensor-on-Tensor Regression}\label{sec:theory}
		
		For technical convenience in the convergence analysis of RGD and RGN, we first introduce the Tensor Restricted Isometry Property (TRIP).
		\begin{Definition}[Tensor Restricted Isometry Property (TRIP)]\label{def: RIP}
			Let $\scA: \bbR^{p_1 \times \cdots \times p_{d+m}} \to \bbR^{n \times p_{d+1} \times \cdots \times p_{d+m}}$ be a linear map. For a fixed $(d+m)$-tuple $\br = (r_1,\ldots,r_{d+m})$ with $1 \leq r_k \leq p_k$, define the $\br$-tensor restricted isometry constant to be the smallest number $R_{\br}$ such that $(1-R_{\br}) \|\bcZ\|^2_{\F} \leq \|\scA(\bcZ)\|_\F^2 \leq (1+R_{\br}) \|\bcZ\|_{\F}^2$ holds for all $\bcZ$ of Tucker rank at most $\br$. If $0\leq R_{\br} < 1$, we say $\scA$ satisfies $\br$-tensor restricted isometry property ($\br-$TRIP).
		\end{Definition}
		TRIP can be seen as a tensor generalization of the popular restricted isometry property (RIP) \citep{candes2011tight}. TRIP was used in various tensor inverse problems \citep{rauhut2017low}. The next Proposition \ref{prop: TRIP-under-Gaussian-design} shows $\scA$ satisfies TRIP with high probability when $\scA$ is generated from a sufficient number of sub-Gaussian measurements.
		
		\begin{Proposition}[TRIP Under sub-Gaussian] \label{prop: TRIP-under-Gaussian-design} Suppose $\scA$ is defined as \eqref{eq: affine operator} and each entry of $\bcA_i$ is independently drawn from mean zero variance $1/n$ sub-Gaussian distributions. There exists universal constants $C, c > 0$ such that for any Tucker rank $\br = (r_1,\ldots,r_{d+m})$ and $0\leq R_{\br} < 1$, as long as $n \geq C(\sum_{i=1}^d (p_i-r_i)r_i + \prod_{i=1}^d r_i )\log(d)/R^2_{\br}$, $\scA$ satisfies the TRIP with $\br$-TRIP constant $R_{\br}$ with probability at least $1 - \exp(-c(\sum_{i=1}^d p_i))$. \end{Proposition}

		Now, we are ready to present the convergence theories for RGD and RGN.
		\begin{Theorem}[Convergence of RGD]\label{th: local contraction general setting-RGD}
			Assume the tensor rank of $\bcX^*$ is $\br^*$ and the input rank to Algorithm \ref{alg:RGD-RGN} is $\br \geq \br^*$. Suppose $\scA$ satisfies {$2\br $}-TRIP, and the initialization $\bcX^0$ satisfies $	\|\bcX^0 - \bcX^*\|_{\F} \leq \frac{ R_{2\br} }{(d+m)(1+R_{2\br +\br^*}-R_{2\br})}\underline{\lambda}$,
			where $\underline{\lambda}:= \min_{k=1,\ldots,d+m}\sigma_{r^*_k}(\cM_k(\bcX^*))$ is the minimum of least singular values at each matricization of $\bcX^*$. In addition, we assume $R_{2\br} \leq \frac{1}{8(\sqrt{d+m} + 1) + 1} $ and $ \underline{\lambda} \geq \frac{2(1 + R_{2\br + \br^*} - R_{2\br} ) (\sqrt{d+m} + 1)(d+m)}{R_{2\br} (1 - R_{2\br})} \|(\scA^*(\bcE ))_{\max(2\br)}\|_{\F} $. Then for all $t\ge 0$,
			\begin{equation} \label{ineq: RGD-error-bound}
				\begin{split}
					& \|\bcX^{t} - \bcX^* \|_{\F} \\
					\leq & 2^{-t} \|\bcX^0 - \bcX^*\|_{\F} + \frac{2(\sqrt{d+m} + 1)}{1 - R_{2\br}} \|(\scA^*(\bcE ))_{\max(2\br)}\|_{\F}.    
				\end{split}
			\end{equation} 
			Recall $(\scA^*(\bcE ))_{\max(2\br)}$ denotes the best Tucker rank $2\br$ approximation of the tensor $\scA^*(\bcE )$.
			
			Especially if $\bcE = 0$, $\{\bcX^t\}$ converges linearly to $\bcX^*$:
			\begin{equation*}
				\|\bcX^{t} - \bcX^* \|_{\F} \leq 2^{-t} \|\bcX^0 - \bcX^*\|_{\F}, \quad \forall\, t \ge 0.
			\end{equation*}
		\end{Theorem}
		
		\begin{Theorem}[Convergence of RGN]\label{th: local contraction general setting-RGN}
			Assume the tensor rank of $\bcX^*$ is $\br^*$ and the input rank to Algorithm \ref{alg:RGD-RGN} is $\br \geq \br^*$. Suppose $\scA$ satisfies {$2\br$}-TRIP and the initialization $\bcX^0$ satisfies $\|\bcX^0 - \bcX^*\|_{\F} \leq \frac{1-R_{2\br} }{4(d+m) (\sqrt{d+m} + 1)(1+R_{2\br + \br^*}-R_{2\br})}\underline{\lambda}$. Then for all $t\ge 0$,
			\begin{equation*}
				\begin{split}
					& \|\bcX^t - \bcX^*\|_{\F} \\
					\leq & 2^{-2^{t}} \|\bcX^0 - \bcX^*\|_{\F} + \frac{2(\sqrt{d+m} + 1)}{1 - R_{2\br}} \|(\scA^*(\bcE ))_{\max(2\br)}\|_{\F}.
				\end{split}
			\end{equation*} 
			Especially if $\bcE = 0$, $\{\bcX^t\}$ converges quadratically to $\bcX^*$:
			\begin{equation*}
				\|\bcX^t - \bcX^*\|_{\F} \leq 2^{-2^{t}} \|\bcX^0 - \bcX^*\|_{\F}, \quad \forall\, t \ge 0.
			\end{equation*}
		\end{Theorem}
		Theorems \ref{th: local contraction general setting-RGD} and \ref{th: local contraction general setting-RGN} show that with proper assumptions on $\scA$ and initialization, iterates of RGD and RGN converge linearly and quadratically to the ball of center $\bcX^*$ and radius $O(\|(\scA^*(\bcE ))_{\max(2\br)}\|_{\F})$, respectively. If $\bcE = 0$, i.e., in the noiseless case, $\bcX^t$ generated by RGD/RGN converges linearly/quadratically to the exact $\bcX^*$. These results show the convergence of RGD and RGN are both robust against rank over-parameterization. We note that the error bound $O(\|(\scA^*(\bcE ))_{\max(2\br)}\|_{\F})$, which is achievable by RGD and RGN, depends on the input rank $\br$ and will increase as $\br$ increases. This is confirmed by the simulation study in Section \ref{sec: numerics-input-rank-effect}, indicating that selecting an appropriate input rank $\br$ remains crucial for the accuracy of the estimators.

	One challenge in establishing Theorems \ref{th: local contraction general setting-RGD} and \ref{th: local contraction general setting-RGN} is to show the contraction of the iterates in the rank overspecified scenario. Standard analysis will result in a condition which requires $ \underline{\lambda}':= \min_{k = 1,\ldots, d+m} \sigma_{r_k}(\cM_k(\bcX^*))$ to be larger than some positive threshold. However, it can never be satisfied since $\underline{\lambda}'$ is zero in the rank overspecified scenario. Instead, we show via a refined analysis that lower bounding $\underline{\lambda}$ is still enough. One such example is Lemma \ref{lm: orthogonal projection} in Section \ref{sec:technical contribution}, where we obtain a projection error bound proportional to $\underline{\lambda}$ rather than $\underline{\lambda}'$ even in the rank overspecified scenario.
		
		\begin{Remark}[General Input Rank and Under-parameterization]\label{rm: under-parameterization}
			Suppose $\br$ is a general input rank (possibly under-parameterized, e.g., $r_k < r_k^*$ for some $k$), we can rewrite \eqref{eq:model} into $\bcY_i = \langle \bcA_i, \bcX' \rangle_* + \bcE_i'$, where $\bcX'$ is the best rank $\br$ approximation of $\bcX^*$ and $\bcE_i' = \bcE_i + \langle \bcA_i, \bcX^* - \bcX' \rangle_* $. Similar results to Theorems \ref{th: local contraction general setting-RGD} and \ref{th: local contraction general setting-RGN} hold if $\bcE_i$ is replaced by $\bcE_i'$. We have the following contraction error bounds for RGD and RGN for general input rank and under-parameterized cases:
			$$ \|\bcX^{t} - \bcX' \|_{\F} \leq  \frac{\|\bcX^0 - \bcX'\|_{\F}}{2^{t}} + \frac{2(\sqrt{d+m} + 1)}{1 - R_{2\br}} \|(\scA^*(\bcE' ))_{\max(2\br)}\|_{\F},$$ 
			$$\|\bcX^t - \bcX'\|_{\F}  \leq  \frac{\|\bcX^0 - \bcX'\|_{\F}}{2^{2^{t}}} + \frac{2(\sqrt{d+m} + 1)}{1 - R_{2\br}} \|(\scA^*(\bcE' ))_{\max(2\br)}\|_{\F}.$$
		\end{Remark}	
		
		\begin{Remark}\label{rmk: convergence-remark-bounded-rank}
			{\bf (Convergence Guarantees Under Over-parameterized Scenario Compared with Literature)} When $\bcE = 0$, the convergent point $\bcX^*$ of RGD and RGN has Tucker rank $\br^*$, which falls out of the manifold $\bbM_{\br}$ when $\br > \br^*$, i.e., the over-parameterized scenario. Because of this, the standard convergence theory of RGD/RGN does not imply the convergence results in Theorems \ref{th: local contraction general setting-RGD} and \ref{th: local contraction general setting-RGN} to our best knowledge. Especially in the low-rank matrix trace regression setting, \cite[Theorem 4.1]{barber2018gradient} established a local convergence result of RGD with a bounded rank constraint for a general objective $f$ satisfying restricted strong convexity and smoothness. However, the local convergence radius implied by their theory shrinks to $0$ in our setting and does not directly apply. Also see more discussions on the convergence of various Riemannian optimization algorithms with bounded rank constraints in \cite{schneider2015convergence,levin2021finding,olikier2022apocalypse}.
		\end{Remark}

		\begin{Remark}[Conditions]\label{rm: conditions}
			We impose the mild condition $\underline{\lambda} \geq \Omega(\|(\scA^*(\bcE ))_{\max(2\br)}\|_{\F})$ while analyzing RGD. Since the forthcoming Theorem \ref{th: lower-bound} shows $\Omega(\|(\scA^*(\bcE ))_{\max(2\br)}\|_{\F})$ is the essential statistical error, $\underline{\lambda} \leq O(\|(\scA^*(\bcE ))_{\max(2\br)}\|_{\F})$ can be a trivial case from a statistical perspective because the initialization $\bcX^0$ is already optimal and no further refinement is needed in such the scenario. Another key condition on initialization will be discussed in Section \ref{sec: applications-initializations}.
		\end{Remark}

		Next, we show in two ways that the statistical error $O\left(\|(\scA^*(\bcE ))_{\max(2\br)}\|_{\F}\right)$ achieved by RGD and RGN is essential. First, in Theorem \ref{th: MLE-guarantee}, we show the estimators with small loss, such as the global minimizer of the loss function \eqref{eq:minimization}, achieve the same error rate. 
		\begin{Theorem}[Upper Bound for Estimators with Small Loss and Global Minimizers] \label{th: MLE-guarantee}
			Suppose $\scA$ satisfies $2\br$-TRIP with TRIP constant $R_{2\br}$ (Definition \ref{def: RIP}). Let $\widehat{\bcX}$ be any estimator such that $\tuckerrank(\widehat{\bcX}) \leq \br$ and $\|\bcY - \scA(\widehat{\bcX})\|_\F^2 \leq \|\bcY - \scA(\bcX^*)\|_\F^2$, i.e., the loss function value of $\widehat{\bcX}$ is no bigger than $\bcX^*$. Then $\|\widehat{\bcX} - \bcX^*\|_\F \leq \frac{2}{1-R_{2\br}} \|(\scA^*(\bcE ))_{\max(2\br)}\|_{\F}.$ 
		\end{Theorem}
		
		Second, we focus on the Gaussian ensemble design, which has been widely considered as a benchmark-setting in the literature on compressed sensing, and matrix/tensor regression \citep{candes2011tight,raskutti2015convex}. In Theorem \ref{th: lower-bound}, we establish the minimax estimation error rate under Gaussian ensemble design, which demonstrates the statistical optimality of RGD and RGN when $d$ and $m$ are constants.
		\begin{Definition}[Tensor-on-tensor Regression Under Gaussian Ensemble Design] \label{def: gaussian-design} We say the tensor-on-tensor regression \eqref{eq:model} is generated from the Gaussian ensemble design if $\{\bcA_i\}_{i=1}^n$ and $\{ \bcE_i \}_{i=1}^n$ are generated independently, $\bcA_i$ has i.i.d. $N(0,1/n)$ entries, and $\bcE_i$ has i.i.d. $ N(0, \sigma^2/n)$ entries.
		\end{Definition}
		
		\begin{Theorem}\label{th: lower-bound}
			{\bf (Error Bound Under Gaussian Ensemble and Minimax Risk Upper and Lower Bounds)} Consider the tensor-on-tensor regression problem \eqref{eq:model} under Gaussian ensemble design (Definition \ref{def: gaussian-design}) and let $df = \sum_{i=1}^{d+m} r_i(p_i - r_i) + \prod_{i=1}^{d+m} r_i$. 
			\begin{itemize}[leftmargin=*]
				\item (Upper bound) When $n \geq C(\sum_{i=1}^d (p_i-r_i)r_i + \prod_{i=1}^d r_i )\log(d)$ for some large positive constant $C$, with probability at least $1 - \exp(-c_1(d,m) \underline{p})$, $\|(\scA^*(\bcE ))_{\max(2\br)}\|_{\F} \leq c_2(d,m)  \sigma\sqrt{\frac{df}{n}}$ for some $c_1(d,m),c_2(d,m) > 0$, where $\underline{p}:= \min_j p_j$. Furthermore, for $\widehat{\bcX}$ in Theorem \ref{th: MLE-guarantee}, we have
				$\mathbb{E}\|\widehat{\bcX} - \bcX\|_{\F} \leq C_2(d,m)\sigma \sqrt{\frac{df}{n}}$. 
				\item (Lower bound) Consider the parameter space of all $p_1\times \cdots \times p_{d+m}$-dimensional tensors of Tucker rank at most $\br=(r_1,\ldots, r_{d+m})$:
				\begin{equation*}
					\mathcal{F}_{\bp, \br} := \left\{  \bcX \in \bbR^{p_1 \times \cdots \times p_{d+m}}, \tuckerrank(\bcX) \leq \br \right\}.
				\end{equation*} 
				Suppose $\min_k r_k \geq C'$ for some absolute constant $C'$. Then there exists a absolute constant $c > 0$ that does not depend on $\br$ and $\bp$ such that $\inf_{\widehat{\bcX}} \sup_{ \bcX \in \mathcal{F}_{\bp,\br} } \bbE \| \widehat{\bcX} - \bcX \|_\F \geq c\sigma \sqrt{\frac{df}{n}}.$
			\end{itemize}
		\end{Theorem}

		\section{Applications, Initialization, and Guarantees in Specific Scenarios}\label{sec: applications-initializations}
		
		The convergence theory in Theorems \ref{th: local contraction general setting-RGD} and \ref{th: local contraction general setting-RGN} rely on a good initialization. As it is challenging to develop a universal initialization algorithm that handles all settings of tensor-on-tensor regression with provable guarantees, we focus on the four most representative cases appearing in applications and literature, {\it scalar-on-tensor regression}, {\it tensor-on-vector regression}, {\it matrix trace regression}, and {\it rank-$1$ tensor-on-tensor regression} to show various spectral methods yield adequate initializations. 
		
		\subsection{Scalar-on-tensor Regression}\label{sec: init-scalar-tensor}
		
		The scalar-on-tensor regression corresponds to the general tensor-on-tensor regression model \eqref{eq:model} with $m = 0$. It can be written as
		\begin{equation} \label{eq: scalar-on-tensor-model}
			\y = \scA(\bcX^*) + \bvarepsilon, \text{ or } \y_i = \langle \bcA_i, \bcX^* \rangle + \bvarepsilon_i,\, i\in [n].
		\end{equation}
		Here, $\y \in \bbR^n$ are observations, $\bvarepsilon \in\mathbb{R}^n$ are unknown noise, and $\bcX^* \in \bbR^{p_1 \times \cdots \times p_d}$ is an order-$d$ Tucker rank $\br^*$ tensor that links response $\y_i$ to tensor covariates $\bcA_i$, which is the parameter of interest. $\scA(\bcX^*) = (\langle \bcA_1, \bcX^*\rangle, \ldots, \langle \bcA_n, \bcX^*\rangle)^\top$. We propose the following Algorithm \ref{alg: initial-scalar-on-tensor} on initialization. 
		
		\begin{algorithm}[!h] \caption{Initialization for (Over-parameterized) Scalar-on-tensor Regression}\label{alg: initial-scalar-on-tensor}
			\begin{algorithmic}[1]
				\State \textbf{Input:} $\mathbf{y}_i \in \mathbb{R}$, $\bcA_i \in \bbR^{p_{1} \times \cdots \times p_{d}}$ for $i =1,\ldots,n$ and input Tucker rank $\br = (r_1,\ldots, r_d)$.
				\State Calculate $\widetilde{\U}_k^{0} = \SVD_{r_k}(\mathcal{M}_k(\scA^*(\y))), k=1,\ldots, d$.
				\State
				For $k = 1$ to $d$, apply one-iteration HOOI, i.e., calculate $$\widetilde{\U}_k^{1} =\SVD_{r_k}   \big(\cM_{k}( \scA^*(\y) \times_{j < k} ( \widetilde{\U}_{j}^{0})^\top \times_{j > k} (\widetilde{\U}_{j}^{0})^\top)\big).$$
				Recall $\SVD_r(\cdot)$ returns the matrix composed of the leading $r$ left singular vectors of matrix ``$\cdot$".
				\State \textbf{Output}: $\bcX^0 = \scA^*(\y)\times_{k = 1}^d \widetilde{\U}_k^{1}(\widetilde{\U}_k^{1})^\top$.
			\end{algorithmic}
		\end{algorithm}

		\begin{Theorem}[Initialization and Overall Guarantees in Scalar-on-tensor Regression] \label{th: initialization-scalar-on-tensor} Consider the over-parameterized scalar-on-tensor regression under Gaussian ensemble design. Denote $df = \sum_{i=1}^d (p_i -r_i) r_i + \prod_{i=1}^d r_i $ and suppose $n \geq c(d)  \left( \frac{(\|\bcX^*\|_{\F}^2 + \sigma^2)}{\underline{\lambda}^2} \left((\prod_{i=1}^d p_i)^{1/2} +  df \right) \right)$ for some constant $c(d)$. 
			Then with probability at least $1 - \underline{p}^{-C}$ for some $C > 0$,
			\begin{itemize}[leftmargin=*]
				\item 
				$\bcX^0$ returned from Algorithm \ref{alg: initial-scalar-on-tensor} satisfies the initialization conditions in Theorems \ref{th: local contraction general setting-RGD} and \ref{th: local contraction general setting-RGN};
				\item consider RGD and RGN initialized with $\bcX^0$, then as long as $t_{\max} \geq \log\left( \frac{\underline{\lambda} \sqrt{n/df} }{c_1(d) \sigma }   \right) \vee 0 $ for RGD or $t_{\max} \geq\log \log \left ( \frac{\underline{\lambda} \sqrt{n/df} }{ c_2(d)  \sigma } \right) \vee 0 $ for RGN, we have the output of RGD or RGN satisfies $\left\| \bcX^{t_{\max} } - \bcX^* \right\|_\F \leq c_3(d) \sigma \sqrt{\frac{df}{n}}.$
			\end{itemize}
		\end{Theorem}
		In establishing Theorem \ref{th: initialization-scalar-on-tensor}, we introduce a new perturbation bound for over-parameterized tensor decomposition. See Theorem \ref{th: perturbation-HOOI} in Section \ref{sec:technical contribution} for more details. Compared with RGD, RGN only requires a {\it double logarithmic} number of iterations to achieve the same $O(\sigma \sqrt{\frac{df}{n}})$ error rate.
		
		\begin{Remark}[Sample Complexity for Over-parameterized Scalar-on-tensor Regression]\label{rm: scalar-on-tenor-sample-complexity} Suppose $\|\bcX^*\|_{\F}^2 \geq C \sigma^2$ for some $C >0$, $\kappa := \bar{\lambda}/\underline{\lambda} = O(1)$ where $\bar{\lambda} = \max_{k=1,\ldots,d} \sigma_1(\cM_k(\bcX^*)) $ and $p_1 = p_2 = \ldots = p$, $r_1 = r_2 = \ldots = r$, $r^*_1 = r_2^* \ldots = r^*$, then the overall sample complexity for RGD/RGN in over-parameterized scalar-on-tensor regression with spectral initialization is $\Omega(r^*(p^{d/2}   + pr+ r^d ) )$. Compared to the sample complexity required for the global minimizer (see Theorem \ref{th: MLE-guarantee}) in this example, i.e., $\Omega(pr+r^d)$ proved in Theorem \ref{th: lower-bound}, there is a significant gap between what can be achieved by the inefficient global minimizer and efficient RGD/RGN algorithms. Rigorous evidence for this statistical-computational gap will be provided in Section \ref{sec: computational-limit}.
		\end{Remark}

		\subsection{Tensor-on-vector Regression} \label{sec: init-tensor-on-vector}
		In this section, we consider the tensor-on-vector regression model:
		\begin{equation} \label{eq: tensor-on-vector-regression}
			\bcY_i = \bcX^* \times_1 \a_i^\top + \bcE_i, \quad \text{ for } i = 1,\ldots,n,
		\end{equation}
		where $\bcY_i, \bcE_i \in \bbR^{p_2 \times \cdots \times p_{m+1}}$ are the observation and noise, $\bcX^* \in \bbR^{p_1 \times \cdots \times p_{m+1}}$ is the parameter tensor of interest with Tucker rank $\br^*$ and $\a_i \in \bbR^{p_1}$ is the covariate vector. We can also write the model compactly as $\bcY = \bcX^* \times_1 \A + \bcE$ where $\bcY,\bcE \in \bbR^{n  \times \cdots \times p_{m+1}}$, $\bcY_{[i,:,\ldots,:]} = \bcY_i$, $\bcE_{[i,:,\ldots,:]} = \bcE_i$ and $\A = [\a_1, \ldots, \a_n]^\top \in \bbR^{n \times p_1}$ is the collection of covariate vectors. We propose the following Algorithm \ref{alg: initial-tensor-on-vector} for initialization and its guarantee is provided in Theorem \ref{th: initialization-tensor-on-vector}.
		
		\begin{algorithm}[!h] \caption{Initialization for (Over-parameterized) Tensor-on-vector Regression}\label{alg: initial-tensor-on-vector}

			\begin{algorithmic}[1]
				\State \textbf{Input:} $\bcY_i \in \mathbb{R}^{p_2 \times \cdots \times p_{m+1}}$, $\a_i \in \bbR^{p_{1}}$ for $i =1,\ldots,n$ and input Tucker rank $\br = (r_1,\ldots, r_{m+1})$.
				\State Compute the QR decomposition of $\A$ and denote it by $\Q_\A \R_\A$.
				\State Calculate $\widetilde{\U}_k^{0} = \SVD_{r_k}(\cM_k( \bcY \times_1 \Q_\A^\top )), k=1,\ldots, m+1$.
				\State
				For $k = 1$ to $m+1$, apply one-iteration HOOI, i.e., calculate $$\widetilde{\U}_k^{1} =\SVD_{r_k}   \big(\cM_{k}( (\bcY \times_1 \Q_\A^\top) \times_{j < k} ( \widetilde{\U}_{j}^{0})^\top \times_{j > k} (\widetilde{\U}_{j}^{0})^\top)\big).$$
				\State Compute $\widebar{\bcX}^0 =(\bcY \times_1 \Q_\A^\top) \times_{k = 1}^{m+1} \widetilde{\U}_k^{1}(\widetilde{\U}_k^{1})^\top$.
				\State Return $\bcX^0 = \widebar{\bcX}^0 \times_1 \R_\A^{-1}$.
				\State \textbf{Output:}
			$\bcX^0$.
			\end{algorithmic}
		\end{algorithm}
		\begin{Theorem}[Initialization and Overall Guarantees in Tensor-on-vector Regression] \label{th: initialization-tensor-on-vector} Consider the over-parameterized tensor-on-vector regression under Gaussian ensemble design. Denote $df = \sum_{i=1}^{m+1} (p_i -r_i) r_i + \prod_{i=1}^{m+1} r_i $. Suppose $$n \geq c(m) \left( \left((\prod_{i=1}^{m+1} p_i)^{1/2} +  df \right) \sigma^2/\underline{\lambda}^2 + p_1 \right) $$ for some constant $c(m)$. Then with probability at least $1 - \exp(-c\underline{p})$ for some $c > 0$,
			\begin{itemize}[leftmargin=*]
				\item $\bcX^0 $ returned from Algorithm \ref{alg: initial-tensor-on-vector} satisfies the initialization conditions in Theorems \ref{th: local contraction general setting-RGD} and \ref{th: local contraction general setting-RGN};
				\item moreover, consider RGD and RGN initialized with $\bcX^0$, then as long as $t_{\max} \geq \log\left( \frac{\underline{\lambda} \sqrt{n/df} }{c_1(m) \sigma }   \right) \vee 0 $ for RGD or $t_{\max} \geq\log \log \left ( \frac{\underline{\lambda} \sqrt{n/df} }{ c_2(m)  \sigma } \right) \vee 0 $ for RGN, we have the output of RGD or RGN satisfies
				\begin{equation*}
					\| \bcX^{t_{\max} } - \bcX^* \|_\F \leq c_3(m) \sigma \sqrt{\frac{df}{n}}.   
				\end{equation*} 
			\end{itemize}
		\end{Theorem}

		\subsection{Matrix Trace Regression} \label{sec: init-matrix-regression}
		In this model, we observe
		\begin{equation} \label{eq: scalar-on-matrix-model}
			\y_i = \langle \A_i, \X^* \rangle + \bvarepsilon_i, i=1,\ldots, n; \text{ or } \y = \scA(\X^*) + \bvarepsilon,
		\end{equation} 
		where $\y, \bvarepsilon \in\mathbb{R}^n$ are observations and unknown noise and $\X^* \in \bbR^{p_1 \times p_2}$ is a rank $r^*$ parameter matrix of interest. 
		
		In matrix trace regression, we can take the retraction map $\cH_r$ in RGD and RGN as the best rank $r$ matrix projection operator: 
		$\cP_r(\B) = \U_{[:, 1:r]}\bSigma_{[1:r, 1:r]}\V_{[:,1:r]}^\top$, where  $\B=\U\bSigma\V^\top$ is the SVD.
		Different from the low-rank projection for tensor of order 3 or higher, $\cP_r$ can be computed efficiently by truncated SVD. Moreover, suppose $\X^t$ has economic SVD $\U^t \bSigma^t \V^{t\top}$, then the projection of $\Z \in \bbR^{p_1 \times p_2}$ onto the tangent space $T_{\X^t} \bbM_r$ can be written succinctly as $P_{T_{\X^t}}(\Z) = P_{\U^t} \Z P_{\V^t} + P_{\U_\perp^t} \Z P_{\V^t} + P_{\U^t} \Z P_{\V_\perp^t}.$
		
		We have the following corollary on the guarantees of RGD and RGN in over-parameterized matrix trace regression.
		\begin{Corollary}[Convergence of RGD/RGN in Matrix Trace Regression]\label{coro: local contraction general setting-RGD-RGN-matrix} Consider the (over-parameterized) matrix trace regression model in \eqref{eq: scalar-on-matrix-model} with $r \geq r^*$. Let $\cH_r$ be the rank $r$ truncated SVD. 
			Suppose $\scA$ satisfies {$2r$}-RIP.
			
			(RGD) Suppose the initialization $\X^0$ satisfies $	\|\X^0 - \X^*\|_{\F} \leq \frac{ R_{2r} }{(1+R_{2r +r^*}-R_{2r})}\sigma_{r^*}(\X^*)$. In addition, we assume $R_{2r} \leq \frac{1}{17} $ and $ \sigma_{r^*}(\X^*) \geq \frac{4(1 + R_{2r + r^*} - R_{2r} ) }{R_{2r} (1 - R_{2r})} \|(\scA^*(\bvarepsilon ))_{\max(2r)}\|_{\F} $. Then $\{\X^t\}$ generated by RGD satisfy for all $t\ge 0$,
			\begin{equation*}
				\|\X^{t} - \X^* \|_{\F} \leq 2^{-t} \|\X^0 - \X^*\|_{\F} + \frac{4}{1 - R_{2r}} \|(\scA^*(\bvarepsilon ))_{\max(2r)}\|_{\F}. 
			\end{equation*} 
			
			(RGN) If the initialization $\X^0$ satisfies $\|\X^0 - \X^*\|_{\F} \leq \frac{1-R_{2r} }{8(1+R_{2r + r^*}-R_{2r})}\sigma_{r^*}(\X^*)$. Then $\{\X^t\}$ generated by RGN satisfy for all $t\ge 0$,
			\begin{equation*}
				\|\X^t - \X^*\|_{\F} \leq 2^{-2^{t}} \|\X^0 - \X^*\|_{\F} + \frac{4}{1 - R_{2r}} \|(\scA^*(\bvarepsilon ))_{\max(2r)}\|_{\F}.
			\end{equation*} 
			
			Especially if $\bvarepsilon = 0$, 
			$\|\X^{t} - \X^* \|_{\F} \leq 2^{-t} \|\X^0 - \X^*\|_{\F}$ for RGD and $\|\X^t - \X^*\|_{\F} \leq 2^{-2^{t}} \|\X^0 - \X^*\|_{\F}$ for RGN.
		\end{Corollary}

		An efficient initialization for the matrix trace regression is $\X^0 = \cP_r(\scA^*(\y))$. The guarantee of $\X^0$ and overall performance of RGD and RGN in matrix trace regression are given in Theorem \ref{th: matrix-initialization}. 
		
		\begin{Theorem}\label{th: matrix-initialization}
			{\bf (Initialization and Overall Guarantees in Over-parameterized Matrix Trace Regression)} Consider the over-parameterized matrix trace regression under Gaussian ensemble design. Denote $df = (p_1 + p_2-r)r $ and suppose $n \geq   \frac{C(\sigma^2 + \|\X^*\|_\F^2)}{\sigma^2_{r^*}(\X^*)}    df$ for some $C > 0$. Then with probability at least $1 - \exp(-c \underline{p} )$,
			\begin{itemize}[leftmargin=*]
				\item $\X^0 = \cP_r(\scA^*(\y)) $ satisfies the initialization conditions in Corollary \ref{coro: local contraction general setting-RGD-RGN-matrix}; 
				\item moreover, consider RGD and RGN initialized with $\X^0$, then as long as $t_{\max} \geq \log\left( \frac{\sigma_{r^*}(\X^*) }{c_1 \sigma } \sqrt{\frac{n}{df}}  \right) \vee 0 $ for RGD or $t_{\max} \geq\log \log \left ( \frac{\sigma_{r^*}(\X^*)  }{ c_2  \sigma }\sqrt{\frac{n}{df}} \right) \vee 0 $ for RGN, we have the output of RGD or RGN satisfies
				\begin{equation*}
					\| \bcX^{t_{\max} } - \bcX^* \|_\F \leq c_3 \sigma \sqrt{\frac{df}{n}}.   
				\end{equation*} 
			\end{itemize}
		\end{Theorem}

		\begin{Remark}\label{rm: comparision-in-matrix} 
			{\bf (Comparison with Existing Results on Over-parameterized Matrix Trace Regression)} Recently, \cite{zhuo2021computational,zhang2021preconditioned} studied the local convergence of factorized gradient descent (GD) in the same setting as ours. In particular, \cite{zhuo2021computational} showed the convergence rate of the original factorized GD slows down to being sublinear when the input rank $r$ is greater than the actual rank $r^*$.  \cite{zhang2021preconditioned} proposed to overcome that by preconditioning the factorized GD; they showed that the convergence rate of preconditioned factorized GD can be boosted back to linear for all $r \geq r^*$. However, the preconditioning step in \cite{zhang2021preconditioned} requires a carefully chosen damping parameter in each iteration and such the choice depends on the unknown noise variance. In contrast, our proposed RGD and RGN algorithms are easy to implement, tuning-free, and are unified in both rank correctly-specified and overspecified settings. In addition, in terms of the theoretical guarantees, the estimation error bound in \cite{zhang2021preconditioned} is suboptimal in the noisy setting, while our bound is minimax optimal as shown in Theorem \ref{th: lower-bound}. Finally, our result is also more general since our $\X^*$ can be a general rank $r^*$ matrix while existing works only focus on positive-semidefinite $\X^*$. The readers are referred to Table \ref{tab: comparison table} for a summary of comparisons.
			
			Meanwhile, to satisfy $r$-RIP, we need $n = \Omega((p_1 + p_2)r)$, so our theory is still based on the ``sample size ($n$) $\geq$ parameter degree of freedom ($df$)" scenario. A follow-up question is whether the ``implicit regularization" phenomenon discussed in the {\bf Related Prior Work} Section appears in Riemannian formulated matrix trace regression in the highly over-parameterized regime, i.e., ``$df > n$," as such phenomenon was recently observed in factorized gradient descent \citep{gunasekar2017implicit,li2018algorithmic}. In fact, the direct application of RGD proposed in this paper does not enjoy implicit regularization in the highly over-parameterized regime because when the input rank $r$ is equal to $p_1 \wedge p_2$, RGD reduces to gradient descent in the whole $p_1$-by-$p_2$ matrix parameter space, which does not enjoy implicit regularization as it will converge to the minimum Frobenius norm solution in this over-parameterized setting with near origin initialization \citep{gunasekar2017implicit}. Our theory so far does not cover the highly over-parameterized regime and further investigation is left as future work.
		\end{Remark}
		
		\subsection{Rank-$1$ Tensor-on-tensor Regression} \label{sec: rank-1-initialization}
  For the general tensor-on-tensor regression model, although $\bbE(\scA^*(\bcY)) = \bcX^*$ is low-rank, the noise structure of $\scA^*(\bcY) - \bcX^*$ is complicated that significantly deviates from the commonly studied additive tensor PCA model in the literature. It is thus challenging to provide an optimal theoretical guarantee for the initialization schemes T-HOSVD and ST-HOSVD in general.

  In this section, we introduce a modified initialization scheme with theoretical guarantees for general $d$ and $m$ when $\bcX^*$ is a rank-$1$ tensor and input rank is also $1$. For simplicity, we assume $n$ is even. Suppose $\bcX^* = \lambda \u_1 \circ \u_2 \circ \cdots \circ \u_{d+m} \in \bbR^{p_1 \times \cdots \times p_{d+m}}$, where ``$\circ$'' denotes the outer product of vectors. Then in this special setting, the model \eqref{eq:model} can be rewritten as 
		\begin{equation} \label{eq: rank-1-model}
			\bcY_i = \lambda\langle \bcA_i, \u_1 \circ \cdots \circ \u_d \rangle \u_{d+1} \circ \cdots \circ \u_{d+m} + \bcE_i. 
		\end{equation} 
  Let $\bcY^1$ and $\bcY^2$ collect $\bcY_i$s in the first and second halves of the data: $\bcY^1_{[i,:,\ldots,:]} = \bcY_i$ and $\bcY^2_{[i,:,\ldots,:]} = \bcY_{n/2+i}$ for $i = 1,\ldots, n/2$. 
We propose an initialization procedure in Algorithm \ref{alg: initial-rank-1-tensor-on-tensor} and provide its theoretical guarantee in Theorem \ref{th: initialization-rank-1-tensor-tensor}. The high-level idea for Algorithm \ref{alg: initial-rank-1-tensor-on-tensor} is as follows: we use the first half of the data $\bcY^1$ to get estimates $\widehat{\u}_{k}$ for $k = d+1, \ldots, d+m$ and then use the second half of the data $\bcY^2$ to estimate $\u_{k}$ for $k = 1, \ldots, d$ after projecting the data to the subspace spanned by $\{\widehat{\u}_{k} \}_{k=d+1}^{d+m}$; finally, a one-iteration HOOI is applied to obtain the initialization. 
		\begin{algorithm}[!h] \caption{Initialization for Rank-$1$ Tensor-on-tensor Regression}\label{alg: initial-rank-1-tensor-on-tensor}
			\begin{algorithmic}[1]
				\State \textbf{Input:} $\bcY_i \in \mathbb{R}^{p_{d+1} \times\cdots \times p_{d+m} }$, $\bcA_i \in \bbR^{p_{1} \times \cdots \times p_{d}}$ for $i =1,\ldots,n$.
				\State Calculate $\widetilde{\u}_k^{0} = \SVD_{1}(\mathcal{M}_{k-d+1}(\bcY^1)), k=d+1,\ldots, d+m$.
				\State Compute $\y_i' = \langle \bcY_i, \widetilde{\u}_{d+1}^{0} \circ \cdots \circ \widetilde{\u}_{d+m}^{0} \rangle$ for $i = n/2+1, \ldots, n$. 
				\State Calculate $\widetilde{\u}_k^{0} = \SVD_{1}(\mathcal{M}_k( \sum_{i=n/2+1}^{n} \y_i' \bcA_i )), k=1,\ldots, d$.
				\State
				For $k = 1$ to $d+m$, apply one-iteration HOOI, i.e., calculate $$\widetilde{\u}_k^{1} =\SVD_{1}   \big(\cM_{k}( \scA^*(\bcY) \times_{j < k} ( \widetilde{\u}_{j}^{0})^\top \times_{j > k} (\widetilde{\u}_{j}^{0})^\top)\big).$$
				\State \textbf{Output}: $\bcX^0 = \scA^*(\bcY)\times_{k = 1}^{d+m} \widetilde{\u}_k^{1}(\widetilde{\u}_k^{1})^\top$.
			\end{algorithmic}
		\end{algorithm}
		\begin{Theorem}[Initialization and Overall Guarantees in Rank-$1$ Tensor-on-tensor Regression] \label{th: initialization-rank-1-tensor-tensor} Consider the rank-$1$ tensor-on-tensor regression under Gaussian ensemble design \eqref{eq: rank-1-model}. Denote $df = \sum_{i=1}^{d+m} p_i $ and suppose $\lambda > C' \sigma$ for some $C' > 0$. If $n \geq c(d,m)  \left( \frac{(\lambda^2 + \sigma^2)}{\lambda^2} \left((\prod_{i=1}^d p_i)^{1/2} +  \bar{p} \right) + \frac{\sigma^4}{\lambda^4}\left(\prod^{d+m}_{i=d+1} p_i + \bar{p} \right) \right)$ for some constant $c(d,m)$ depending on $d$ and $m$ only, where $\bar{p}= \max_{k=1, \ldots, d+m} p_i$. 
			Then with probability at least $1 - \underline{p}^{-C}$ for some $C > 0$,
			\begin{itemize}[leftmargin=*]
				\item 
				$\bcX^0$ returned from Algorithm \ref{alg: initial-rank-1-tensor-on-tensor} satisfies the initialization conditions in Theorems \ref{th: local contraction general setting-RGD} and \ref{th: local contraction general setting-RGN};
				\item Considering RGD and RGN initialized with $\bcX^0$, as long as $t_{\max} \geq \log\left( \frac{\lambda \sqrt{n/df} }{c_1(d,m) \sigma }   \right) \vee 0 $ for RGD or $t_{\max} \geq\log \log \left ( \frac{\lambda \sqrt{n/df} }{ c_2(d,m)  \sigma } \right) \vee 0 $ for RGN, we have the output of RGD or RGN satisfies $\left\| \bcX^{t_{\max} } - \bcX^* \right\|_\F \leq c_3(d,m) \sigma \sqrt{\frac{df}{n}}.$
			\end{itemize}
		\end{Theorem}

		\section{Computational Limits}\label{sec: computational-limit}
		In this section, we provide rigorous evidence for the computational barrier in scalar-on-tensor regression via the low-degree polynomials method. Without loss of generality, we assume $\bvarepsilon_i \overset{i.i.d.}\sim N(0,\sigma^2)$ with $0 \leq \sigma^2 < 1$, $\bcA_i \overset{i.i.d.}\sim N(0,1)$ and $\|\bcX^*\|_\F + \sigma^2 = 1$ in establishing the computational lower bound for scalar-on-tensor regression \eqref{eq: scalar-on-tensor-model} (see Supplement \ref{sec: equiva-scalar-on-tensor-reg} for a proof). We also consider the setting $p_1 = \ldots = p_d = p$ and $r_1^* = \ldots = r_d^* = r^*$ throughout this section.

		We consider a canonical hypothesis testing formulation of scalar-on-tensor regression:
		\begin{equation}\label{eq: hypo-test-scalar-on-tensor}
			\begin{split}
				& H_0: \{ (\y_i, \rmvec(\bcA_i)) \}_{i=1}^n \overset{i.i.d.}\sim N(0, \I_{1+p^d}), \\
				& H_1: \{ (\y_i, \rmvec(\bcA_i)) \}_{i=1}^n \text{: } \bcX^* = \sqrt{1-\sigma^2} \x^{*\otimes d},\\
				&\quad  \x^* = (x_1^*,\ldots,x_p^*), x_j^* \overset{i.i.d.}\sim \textrm{Uniform}(\{ p^{-1/2}, -p^{-1/2}\}); \\
				&\quad \text{ for } i \in [n], \bcA_i \overset{i.i.d.}\sim N(0,1), \y_i \text{ is  i.i.d. generated}\\
				&\quad \text{ via } \y_i = \langle \bcX^*, \bcA_i \rangle + \bvarepsilon_i, ~
				\bvarepsilon_i \overset{i.i.d.}\sim N(0,\sigma^2).
			\end{split}
		\end{equation}
        Since we aim to develop a lower bound, the hardness result for \eqref{eq: hypo-test-scalar-on-tensor} also implies the hardness result for a bigger class in the sense of minimax. The idea of using low-degree polynomials to predict the statistical-computational gaps is recently developed in a line of work \citep{hopkins2017bayesian,hopkins2018statistical}. In comparison to sum-of-squares (SOS) computational lower bounds, the low-degree polynomials method is simpler to establish and appears to always yield the same results for natural average-case hardness problems. Low-degree polynomials computational hardness results have been provided to a number of problems, such as the planted clique detection \citep{hopkins2018statistical, barak2019nearly}, community detection in stochastic block model \citep{hopkins2017bayesian, hopkins2018statistical}, the spiked tensor model \citep{hopkins2017power, hopkins2018statistical, kunisky2019notes}, the spiked Wishart model \citep{bandeira2020computational}, sparse PCA \citep{ding2019subexponential}, spiked Wigner model \citep{kunisky2019notes}, clustering \citep{loffler2020computationally,davis2021clustering,lyu2022optimal}, planted vector recovery \citep{mao2021optimal}, certifying RIP \citep{ding2020average} and random k-SAT \citep{bresler2022algorithmic}. It is gradually believed that the low-degree polynomials method is able to capture the essence of what makes sum-of-squares algorithms succeed or fail \citep{hopkins2018statistical, kunisky2019notes}. Our results on the computational hardness of distinguishing between $H_0$ and $H_1$ in scalar-on-tensor regression based on low-degree polynomials are given below.

		\begin{Theorem}[Low-degree Hardness for Scalar-on-tensor Regression]\label{th: lower-degree-polynomial-tensor-regression}
			Consider the hypothesis test \eqref{eq: hypo-test-scalar-on-tensor}. For any $0<\delta<1$, if $n\leq \frac{(p/dD)^{d/2}\delta}{2(1-\sigma^2)}$, we have
			\begin{equation} \label{ineq: truncated-likelihood-ratio-bound}
				\sup_{\text{polynomial } f: \substack{deg(f) \leq D \\
					\mathbb{E}_{H_0} f(\{\y_i, \bcA_i\}_{i=1}^n) = 0,\\ \text{Var}_{H_0} f(\{\y_i, \bcA_i\}_{i=1}^n) = 1}}\mathbb{E}_{H_1}f(\{\y_i, \bcA_i\}_{i=1}^n) \leq  \frac{\delta}{1-\delta}.
			\end{equation}
		\end{Theorem}
		
		It has been widely conjectured in the literature that for a broad class of hypothesis testing problems: $H_0$ versus $H_1$, there is a test with runtime $n^{\tilde{O}(D)}$ and Type I + II error tending to zero if and only if there is a successful $D$-simple statistic, i.e., a polynomial $f$ of degree at most $D$, such that $\bbE_{H_0} f(X) = 0$, $\text{Var}_{H_0} (f^2(X)) = 1$, and $\bbE_{H_1} f(X) \to \infty$ \citep{hopkins2018statistical,kunisky2019notes}. Therefore, by setting $D = C \log p$ for any $C > 0$, Theorem \ref{th: lower-degree-polynomial-tensor-regression} provides firm evidence for the statistical-computational gap when $n = O(p^{d/2-\varepsilon})$ for any $\epsilon > 0$. Compared to the sample size requirement in the upper bound mentioned in Remark \ref{rm: scalar-on-tenor-sample-complexity}, the computational lower bound established in Theorem \ref{th: lower-degree-polynomial-tensor-regression} is sharp when $r^*=O(1), r \leq \sqrt{p}$. Our Theorem \ref{th: lower-degree-polynomial-tensor-regression} answers the question raised by \cite{rauhut2017low} on the sample complexity requirement for efficient estimators in scalar-on-tensor regression. We note the first computational hardness evidence for scalar-on-tensor regression was provided recently in \cite{diakonikolas2022statistical} in the Statistical Query model. We complement their results by providing a direct low-degree polynomials argument and figuring out the explicit dependence of the sample complexity on the degrees tolerated in low-degree polynomials. Finally, we also show in the Supplement \ref{sec: supp-test-to-estimation} Proposition \ref{prop: testing-imply-estimation} that the hardness of testing $H_0$ versus $H_1$ implies the hardness of estimating $\bcX^*$. 
		
		\begin{Remark} \label{rem: comp-limit-proof-idea}
			{\bf (Proof Ideas and Comparison with Existing Arguments)} Here we briefly discuss the proof idea of Theorem \ref{th: lower-degree-polynomial-tensor-regression} and the key technical novelty therein. A detailed proof and preliminaries of low-degree polynomials are provided in Supplement \ref{sec: proof-comp-limit}. First, it has been established in \cite{hopkins2018statistical,kunisky2019notes} that the left-hand side of \eqref{ineq: truncated-likelihood-ratio-bound} is equal to the norm of the truncated likelihood ratio under the null:
			\begin{equation} \label{eq: low-degree-truncated-likelihood-ratio}
				\begin{split}
					&\sup_{ \substack{\text{polynomial } f: ~deg(f) \leq D  \\
						\mathbb{E}_{H_0} f(\{y_i, \bcA_i\}_{i=1}^n) = 0,\\ \text{Var}_{H_0} f(\{y_i, \bcA_i\}_{i=1}^n) = 1}}\mathbb{E}_{H_1}f(\{y_i, \bcA_i\}_{i=1}^n) = \sqrt{ \bbE_{H_0} \left( \left( \frac{p_{H_1}( \{y_i, \bcA_i\}_{i=1}^n )}{p_{H_0}(\{y_i, \bcA_i\}_{i=1}^n)} \right)^{\leq D} - 1 \right)^2 },
				\end{split}
			\end{equation} where $p_{H_0}$ and $p_{H_1}$ denote the likelihood under the null and alternative, respectively, and $f^{\leq D}$ is the projection of a function $f$ to the linear subspace of degree-$D$ polynomials, where the projection is orthonormal with respect to the inner product induced under $H_0$. A standard trick to bound the right hand of \eqref{eq: low-degree-truncated-likelihood-ratio} is to evaluate it separately under the orthogonal basis functions $\{f_j \}_{j\geq 1}$ under the null, and then the argument boils down to bound $\sum_{j=1}^D (\bbE_{H_1} f_j(\{y_i, \bcA_i\}_{i=1}^n))^2$, which is the sum of second moments of the orthogonal basis functions under the {\it alternative}. See \eqref{eq: evalutation-truncated-likelihood-ratio} in Supplement \ref{sec: preliminary-low-degree-hermite} for details. There have been many successes in bounding $\sum_{j=1}^D (\bbE_{H_1} f_j(\{y_i, \bcA_i\}_{i=1}^n))^2$ when the testing problem under $H_1$ has the ``signal + noise" structure \citep{hopkins2018statistical, kunisky2019notes}. Such a structure simplifies the analysis as the noise part and signal part are decoupled. In contrast, there is little low-degree polynomial hardness evidence when the problem under $H_1$ has correlated structures, such as the regression problem considered in this paper. One of our main technical contributions in tackling this challenge is a formula for computing the expectation of Hermite polynomials for correlated multivariate Gaussian random variables (Lemma \ref{lm:correlated-hermitian} in Section \ref{sec:technical contribution}). With this key technical tool, we can bound $\sum_{j=1}^D (\bbE_{H_1} f_j(\{y_i, \bcA_i\}_{i=1}^n))^2$ under the $H_1$ in \eqref{eq: hypo-test-scalar-on-tensor} to prove the result. See Supplement \ref{sec: low-degree-lower-bound-proof} for the detailed calculation.

		\end{Remark}

		\begin{Remark}\label{rm: matrix-tensor-sample-compare}
			{\bf (Comparing Rank Overspecification in Matrix Trace Regression and Scalar-on-tensor Regression)} Suppose $r^* = O(1)$. In matrix trace regression, the sample size requirement of the ``spectral initialization + local refinement" estimation scheme is $O(pr)$, where $r$ is the input rank. Thus, the sample complexity increases linearly as the input rank $r$ increases. Meanwhile, the sample complexity of the scalar-on-tensor regression under the same estimation scheme is $O(p^{d/2})$ when $r\leq \sqrt{p}$ (see Remark \ref{rm: scalar-on-tenor-sample-complexity}). Due to the computational lower bound of scalar-on-tensor regression in Theorem \ref{th: lower-degree-polynomial-tensor-regression}, the sample complexity $\Omega( p^{d/2})$ is essential for any polynomial-time algorithm to succeed under proper assumptions. Therefore, no extra samples are needed for efficient estimators in moderate over-parameterized scalar-on-tensor regression; while such a phenomenon does not exist in its matrix counterpart. See Figure \ref{fig: tensor-matrix-comparison} for a pictorial illustration of this distinction. 
			
			In addition to the ``spectral initialization + local refinement", random initialization + refinement by some simple local methods is another effective approach for solving matrix and tensor problems. Such a ``random initialization + local refinement" scheme has been shown to be effective in over-parameterized matrix trace regression, where only $O(pr^{*2})$ samples are needed \citep{li2018algorithmic}. However, initialization with a small enough magnitude and the factorization formulation seem to be critical there. Due to the space limit, we leave a thorough comparison of these two popular approaches for over-parameterized tensor-on-tensor regression problems as future work. 
		\end{Remark}

		\section{Technical Contributions}\label{sec:technical contribution}
		
		We develop several technical tools to establish the theoretical results in this paper. We summarize them in this section.
		
		\noindent{\bf Tackle Over-parameterization in the Convergence Analysis.} 
		In the proof of Theorems \ref{th: local contraction general setting-RGD} and \ref{th: local contraction general setting-RGN}, we first observe that for any $k \in [d+m]$, the mode-$k$ singular subspace of $\bcX^t$, denoted by $\U_k^t$, can be decomposed as $\U_k^t = [\widebar{\U}_k^t \quad \widecheck{\U}_k^t]$ where $\widebar{\U}_k^t$ is composed of the first $r_k^*$ columns of $\U_k^t$ and $\widecheck{\U}_k^t$ is composed of the rest of the $(r_k - r_k^*) $ columns of $\U_k^t$. Then the projection operator onto $\U_{k\perp}^t$, the orthogonal complement of $\U_k^t$, satisfies
		\begin{equation} \label{eq: subspace-decomposition}
			P_{\U_{k \perp}^t} = \I_{p_k} - P_{\U_k^t} = \I_{p_k} - P_{\widebar{\U}_k^t} -  P_{\widecheck{\U}_k^t} =  (\I_{p_k} - P_{\widecheck{\U}_k^t} ) (\I_{p_k} - P_{\widebar{\U}_k^t} ).
		\end{equation} 
		This implies $\|(\I_{p_k} - P_{\U_k^t}) \Z\| \leq \|(\I_{p_k} - P_{\widebar{\U}_k^t} ) \Z \|$ for any matrix $\Z$ with compatible dimension. Based on this property, we can focus on the first $r_k^*$ columns of $\U_k^{t}$ and establish the following lemma, which plays a key role in establishing the convergence of RGD and RGN. 
		\begin{Lemma}[An Over-parameterized Projection Error Bound]
			\label{lm: orthogonal projection}
			Suppose $\bcX^t \in \bbR^{p_1 \times \cdots \times p_{d+m}}$ is an order-$(d+m)$ Tucker rank $\br := (r_1, \ldots,r_{d+m})$ tensor and $\bcX^*\in \bbR^{p_1 \times \cdots \times p_{d+m}}$ is an order-$(d+m)$ Tucker rank $\br^* := (r^*_1, \ldots,r^*_{d+m})$ tensor with $\br^* \leq \br$. Then we have
			\begin{equation*}
				\| P_{(T_{\bcX^t})_\perp} \bcX^* \|_{\F} \leq \frac{2(d+m) \| \bcX^t - \bcX^* \|^2_{\F} }{\underline{\lambda}},
			\end{equation*} 
			where $P_{(T_{\bcX})_\perp} := \I- P_{T_{\bcX}}$ is the orthogonal complement of the projector $P_{T_{\bcX}}$ \eqref{eq: tangent space projector} and $\underline{\lambda} := \min_{k=1,\ldots,d+m} \sigma_{r^*_k}(\cM_k(\bcX^*))$. Especially in the matrix setting, i.e., $d+m=2$, a sharper upper bound holds: $\|P_{(T_{\X^t})_\perp} \X^*\|_\F \leq \frac{2\|\X^t-\X^*\|_\F^2}{\sigma_{r^*}(\X^*)}.$
		\end{Lemma}

		\noindent{\bf Initialization Guarantees for Scalar-on-tensor Regression and Tensor-on-vector Regression}. A key step of Algorithms \ref{alg: initial-scalar-on-tensor} and \ref{alg: initial-tensor-on-vector} is the one-iteration HOOI (OHOOI) algorithm (Algorithm \ref{alg: HOOI} below). Such one loop update improves the dependence of $r^*$ in sample complexity compared to the vanilla T-HOSVD based initialization in both scalar-on-tensor and tensor-on-vector regressions. In the proofs of Theorems \ref{th: initialization-scalar-on-tensor} and \ref{th: initialization-tensor-on-vector}, we develop the following deterministic tensor perturbation bound for OHOOI in the over-parameterized regime.
		\begin{algorithm}[!h] \caption{One-iteration Higher-Order Orthogonal Iteration (OHOOI)} \label{alg: HOOI}
			\begin{algorithmic}[1]
				\State \textbf{Input:} $\widetilde{\bcT}\in \bbR^{p_{1} \times \cdots \times p_{d}}$, initialization $\widetilde{\U}_k^{0}\in \bbO_{p_k,r_k}, k=1,\ldots, d$, input Tucker rank $\br = (r_1,\ldots, r_d)$.
				\State
				For $k = 1$ to $d$, update $\widetilde{\U}_k^{1} =\SVD_{r_k}   \big(\cM_{k}( \widetilde{\bcT} \times_{j < k} ( \widetilde{\U}_{j}^{0})^\top \times_{j > k} (\widetilde{\U}_{j}^{0})^\top )  \big).$
				\State {\bf Output:} $\widehat{\bcT} = \widetilde{\bcT}\times_{k = 1}^d P_{\widetilde{\U}_k^{1}}$.
			\end{algorithmic}
		\end{algorithm}

		\begin{Theorem}[Perturbation Bound for Over-parameterized Tensor Decomposition] \label{th: perturbation-HOOI}
			Suppose $\widetilde{\bcT}, \bcT \in \bbR^{p_1 \times \cdots \times p_d}$, $\bcT$ is of Tucker rank $\br^* = (r^*_1,\ldots,r^*_d)$ with Tucker decomposition $\bcB \times_1 \U_1 \times \cdots \times_d \U_d$, where $\bcB \in \bbR^{r_1^* \times \cdots \times r_d^*}$ and $\U_k \in \bbO_{p_k, r_k^*}$ for $k = 1, \ldots, d$. Let $\bcZ = \widetilde{\bcT} - \bcT$. Suppose the inputs of the OHOOI algorithm are $\widetilde{\bcT}$, Tucker rank $\br = (r_1,\ldots,r_d)$ with $\br \geq \br^*$ and initializations $\widetilde{\U}_k^0 \in \bbO_{p_k, r_k}$ for $k = 1,\ldots,d$. If the initialization error satisfies $\max_{k=1,\ldots,d}\| \widetilde{\U}^{0\top}_{k\perp} \U_k \| \leq \frac{\sqrt{2}}{2} $. Then the output of Algorithm \ref{alg: HOOI}, $\widehat{\bcT}$, satisfies $ \| \widehat{\bcT} - \bcT \|_{\F} \leq (2^{ \frac{d+1}{2} }\cdot d + 1) \|\bcZ_{\max(\br)}\|_{\F}$. 
		\end{Theorem}

		\noindent{\bf Low-degree Polynomials Evidence for Problems With Correlated Structures.} As we have mentioned in Remark \ref{rem: comp-limit-proof-idea}, the main task in the proof of Theorem \ref{th: lower-degree-polynomial-tensor-regression} is to compute the norm of the truncated likelihood ratio. See Supplement \ref{sec: preliminary-low-degree-hermite} for a preliminary of low-degree polynomials method. Since the data are i.i.d. Gaussian under the null hypothesis of \eqref{eq: hypo-test-scalar-on-tensor}, the main challenge boils down to computing the expected Hermite polynomials on correlated multivariate Gaussian. In the following Lemma \ref{lm:correlated-hermitian}, we provide a simple formula for that. This lemma can be useful in establishing low-degree polynomial hardness evidence for other problems with complex structures. Let $\{h_k \}_{k \in \bbN}$ be the normalized univariate Hermite polynomials $h_k = \frac{1}{\sqrt{k!}} H_k$ where $\{ H_k \}_{k \in \bbN}$ are univariate Hermite polynomials which are defined by the following recurrence: $H_0(x) = 1, H_1(x) = x$, $H_{k+1}(x) = x H_k(x) - k H_{k-1}(x)$ for $ k \geq 1$.

		\begin{Lemma}[Expected Hermitian Polynomials on Correlated Multivariate Gaussian]\label{lm:correlated-hermitian}
			Suppose $w$ is a positive integer, $ Y\in \mathbb{R}, \X = (X_1,\ldots,X_w) \in \mathbb{R}^w$ are random variable and random vectors, respectively, and $(Y, \X)\sim \mathcal{N}\left(0, \begin{bmatrix}
				1 & \u^\top\\
				\u & \I_w
			\end{bmatrix}\right)$ with $\u = (u_1,\ldots,u_w)$. For any integers $\alpha, \beta_1,\ldots, \beta_w\geq 0$, $\mathbb{E}\left( h_\alpha(Y)\prod_{j=1}^w h_{\beta_j}(X_j) \right) =  \sqrt{\frac{\alpha!}{\prod_{j=1}^w \beta_j!}}\cdot\prod_{j=1}^w u_j^{\beta_j} 1( \alpha= \sum_{j=1}^w\beta_j )$, where $1(\cdot)$ in the indicator function.
		\end{Lemma}

		\section{Implementation Details of RGD and RGN}\label{sec:implementation}
		
		In this section, we complement the implementation details of RGD and RGN proposed in Section \ref{sec: RGD-and-RGN}. 
		\vskip.3cm
		
		{\bf \noindent Implementation of RGD.} First, by the definition of the adjoint map, $\scA^*: \bbR^{n \times p_{d+1} \times \cdots \times p_{d+m}} \to \bbR^{p_1 \times \cdots \times p_{d+m}}$ satisfies $ \langle \scA(\bcZ_1), \bcZ_2 \rangle = \langle \bcZ_1, \scA^*(\bcZ_2) \rangle$ for any $\bcZ_1 \in \bbR^{p_1 \times \cdots \times p_{d+m}}, \bcZ_2 \in \bbR^{n \times p_{d+1} \times \cdots \times p_{d+m}}$. Simple manipulation yields:
		\begin{equation*}
			\scA^*(\bcZ_2)_{[k_1,\ldots,k_d, j_1,\ldots,j_m]} = \sum_{i=1}^n \bcZ_{2[i,j_1,\ldots,j_m]} \bcA_{i[k_1,\ldots,k_d]}.
		\end{equation*} 
		Combining this with the formula of  projection $P_{T_{\bcX^t}}$ in \eqref{eq: tangent space projector}, we can calculate $\bcX^{t+0.5} = \bcX^t - \alpha_t P_{T_{\bcX^t}} \scA^*( \scA(\bcX^t) - \bcY)$ and implement the RGD update.
		
		\vskip.3cm

		{\bf\noindent Implementation of RGN.} To illustrate the implementation details of RGN, we first introduce the following lemma.

		\begin{Lemma}[Spectrum of $P_{T_{\bcX^t}} \scA^* \scA P_{T_{\bcX^t}}$]\label{lm: spectral norm bound of Atop A}
			Suppose $\bcX^t$ is of Tucker rank at most $\br$ and the linear map $\scA$ satisfies the 2$\br$-TRIP. Then for any tensor $\bcZ\in T_{\bcX^t} \bbM_\br$,
			\begin{equation}\label{eq: spectrum of LAAL}
				(1 - R_{2\br}) \|\bcZ\|_{\F} \leq \|P_{T_{\bcX^t}} \scA^* \scA P_{T_{\bcX^t}}(\bcZ) \|_{\F} \leq (1+R_{2\br})\|\bcZ\|_{\F},
			\end{equation} and
			\begin{equation}\label{eq: spectrum of LAAL inverse}
				\frac{\|\bcZ\|_{\F}}{1+R_{2\br}} \leq \|(P_{T_{\bcX^t}} \scA^* \scA P_{T_{\bcX^t}})^{-1}(\bcZ) \|_{\F} \leq \frac{\|\bcZ\|_{\F}}{1 - R_{2\br}}.
			\end{equation}
		\end{Lemma}
		Lemma \ref{lm: spectral norm bound of Atop A} shows the linear operator $P_{T_{\bcX^t}} \scA^* \scA P_{T_{\bcX^t}}$, which is a mapping from $T_{\bcX^t} \bbM_\br$ to itself, is provably invertible under TRIP condition, which further implies the least squares in RGN update, $\bcX^{t+0.5} = \argmin_{\bcX \in T_{\bcX^t} \bbM_{\br} } \frac{1}{2}\|\bcY - \scA P_{T_{\bcX^t}}(\bcX)\|_\F^2$, has a unique solution. In the following Proposition \ref{prop: efficient-implementation-RGN}, we show that the RGN update can be reduced to solving $(m+1)$ least squares, which renders a fast implementation of RGN.
		
		\begin{Proposition}[Efficient Implementation of RGN Update] \label{prop: efficient-implementation-RGN}
			Suppose $\bcX^t$ has Tucker decomposition $\bcS^{t} \times_{k=1}^{d+m} \U_k^t$. Then the RGN update, i.e., $\bcX^{t+0.5} = \argmin_{\bcX \in T_{\bcX^t} \bbM_{\br} } \frac{1}{2}\|\bcY - \scA P_{T_{\bcX^t}}(\bcX)\|_\F^2$, is equal to $\bcX^{t+0.5} = \bcB^t \times_{k=1}^{d+m} \U_k^t + \sum_{k=1}^{d+m} \bcS^t \times_k \U_{k\perp}^t \D_k^{t} \times_{j \neq k} \U_j^t $, where 
			\begin{itemize}[leftmargin=*]
				\item $(\bcB^t,\{\D_k^{t}\}_{k=1}^d)$ is the solution of the following least squares with design matrix size $ n\prod_{l=d+1}^{d+m} r_l \times ( \prod_{k=1}^{d+m} r_k + \sum_{k=1}^d r_k(p_k - r_k) ) $:
				\begin{equation*}
					\begin{split}
						&(\bcB^t,\{\D_k^{t}\}_{k=1}^d) \\
						=& \argmin_{\substack{\bcB\in \bbR^{r_1 \times \cdots \times r_{d+m}}, \\ \D_k\in \bbR^{(p_k -r_k) \times r_k}, k=1,\ldots,d}} \sum_{i=1}^n \Big\| \bcY_i \times_{l=1}^m \U_{l+d}^{t\top} - \langle \bcA_i \times_{j=1}^d \U_j^{t\top} , \bcB \rangle_* \\
						&- \sum_{k=1}^d \langle \bcA \times_k \U_{k\perp}^{t\top} \times_{j\neq k} \U_j^{t\top} , \bcS^t \times_k \D_k \rangle_*  \Big\|_\F^2\\
						= & \argmin_{\substack{\bcB\in \bbR^{r_1 \times \cdots \times r_{d+m}}, \\ \D_k\in \bbR^{(p_k -r_k) \times r_k}, k=1,\ldots,d}} \sum_{i=1}^n \sum_{j_l \in [r_{d+l}], l = 1,\ldots,m } \\
						& ~~ \Big(  \left(\bcY_i \times_{l=1}^m \U_{l+d}^{t\top} \right)_{[j_1,\ldots,j_m]}  -  \langle \bcA_i \times_{j=1}^d \U_j^{t\top} , \bcB_{[:,\ldots,:,j_1,\ldots,j_m]} \rangle - \\
						& - \sum_{k=1}^d \langle \U_{k\perp}^{t\top} \cM_k\left( \bcA_i \times_{j\neq k} \U_j^{t\top} \right) \left( \cM_k( \bcS^t_{[:,\ldots,:,j_1,\ldots,j_m]} ) \right)^\top  , \D_k \rangle  \Big)^2
					\end{split}
				\end{equation*}
				\item for $k = d+1, \ldots,d+m$, 
				\begin{equation*}
					\D_k^{t \top} = \argmin_{\D_k^\top \in \bbR^{r_k \times (p_k - r_k)}} \| \Y_{ki} - \A_{ki} \D_k^\top \|_\F^2, 
				\end{equation*} where $\A_{ki} = \left( \cM_{k-d}( \langle \bcA_i \times_{j=1}^d \U_j^{t\top} , \bcS^t \rangle_*  ) \right)^\top \in \bbR^{ \prod_{l=d+1,l\neq k}^{d+m} r_{l} \times r_k }  $, $ \Y_{ki} = \left( \cM_{k-d} ( \bcY_i \times_{l \neq k-d} \U_{l+d}^{t\top}  ) \right)^\top \U_{k\perp}^t \in  \bbR^{ \prod_{l=d+1,l\neq k}^{d+m} r_{l} \times (p_k-r_k) }$.
			\end{itemize} 
		\end{Proposition}

		In the tensor-on-vector regression ($d=1$), the update of RGN has a cleaner and fully closed expression as follows.
		\begin{Lemma}[RGN Update in Tensor-on-vector Regression] \label{lm: RGN-tensor-on-vector} Consider the RGN for tensor-on-vector regression in \eqref{eq: tensor-on-vector-regression}. Suppose $\A^\top \A$ is invertible where $\A$ is the collection of covariate vectors and the iterate at iteration $t$ is $\bcX^t = \llbracket \bcS^t; \U^t_1,\U^t_{2}, \ldots, \U^t_{m+1} \rrbracket$. Then the solution $\bcX^{t+0.5}$ in \eqref{eq: alg least square} has a closed-form expression:
			\begin{equation*}
				\bcX^{t+0.5} = \bcB^t \times_{k=1}^{1+m} \U_k^t + \sum_{k=1}^{1+m} \bcS^t \times_{k} \U_{k\perp}^t \D^t_k \times_{j \neq k} \U_j^t, 
			\end{equation*} where 
			\begin{equation*}
				\begin{split}
					&\cM_1(\bcB^t) = ( \U_1^{t\top} \A^\top \A \U_1^t )^{-1} \U_1^{t\top} \A^\top\\
					& \qquad\qquad  \cdot \left( \cM_1(\bcY) \otimes_{j=(1+m)}^2 \U_j^t -  \A \U_{1\perp}^t \U_{1\perp}^{t\top} (\A^\top \A)^{-1} \A^\top \cM_1(\bcY) \W^t_1 \V^{t\top}_1  \right);\\
					&\cM_1(\bcS^t \times_{1} \U_{1\perp}^t \D^t_1 \times_{j \neq 1} \U_j^t) = \U_{1\perp}^t \U_{1\perp}^{t\top} (\A^\top \A)^{-1} \A^\top \cM_1(\bcY) \W^t_1 \W^{t\top}_1;
				\end{split}
			\end{equation*}
			and 
			\begin{equation*}
				\begin{split}
					&\cM_k(\bcS^t \times_{k} \U_{k\perp}^t \D^t_k \times_{j \neq k} \U_j^t)\\
					=& \U_{k\perp}^t \U_{k\perp}^{t\top}  \cM_k(\bcY \times_1 \A^\top) \W^t_k \left( \V_k^{t\top} \left( \otimes_{\substack{i \neq k,\\i\neq 1}} \I_{r_i} \otimes(\U_1^{t\top} \A^\top \A \U_1^t ) \right) \V_k^t \right)^{-1} \W^{t\top}_k
				\end{split}
			\end{equation*} for $k = 2,\ldots, (1+m)$. Recall $\V^t_k=\QR(\cM_k(\bcS^t)^\top)$ and $\W_k^t$ is defined in \eqref{def: Wk}.
		\end{Lemma}

		\section{Numerical Studies}\label{sec:numerics}
		We conduct simulation studies to investigate the numerical performance of RGD/RGN in tensor-on-tensor regression and to verify our theoretical findings. In each simulation, we generate $\bcE_i$ with i.i.d. $ N(0, \sigma^2)$ entries, $\bcA_i $ with i.i.d. $ N(0, 1)$ entries, $\{\U_k\}_{k=1}^{d+m} $ uniformly at random from $\mathbb{O}_{p, r^*}$ for some to-be-specified $p$ and $r^*$, and $\bcS \in \bbR^{r^* \times \cdots \times r^*}$ with i.i.d. $N(0,1)$ entries; then we form $\bcX^* = \bcS \times_1 \U_1 \times \cdots \times_{d+m} \U_{d+m}$ and generate $\bcY_i$ for $i = 1,\ldots,n$. The input rank of RGD and RGN is set to be $\br = (r,\ldots,r)$ and $r \geq r^*$. In the simulation study, we will experiment with various values of $r$. Additionally, $r$ can be chosen by a data-driven approach. See Supplement \ref{sec: rank-selection} for details. For simplicity, we mainly focus on two examples: scalar-on-tensor regression and tensor-on-vector regression. In the scalar-on-tensor regression, we consider $d=3$; in the tensor-on-vector regression, we consider $m = 3$. Spectral initializations discussed in Section \ref{sec: applications-initializations} are applied in both examples.
		
		Throughout the simulation studies, the error metric we consider is the relative root mean squared error (Relative RMSE) $\|\bcX^{t} - \bcX^*\|_{\F}/ \|\bcX^*\|_{\F}$. The algorithm is terminated when it reaches the maximum number of iterations $t_{\max} = 300$ or the corresponding error metric is less than $10^{-13}$. Unless otherwise noted, the reported results are based on averages of 100 simulations and on a computer with Intel Xeon E5-2680 2.5GHz CPU.

		\subsection{Numerical Performance of RGD and RGN} \label{sec: numerical RGD/RGN property}
		
\begin{figure}[ht]
			\centering
			\subfigure[Scalar-on-tensor regression: RGD]{\includegraphics[width=0.46\textwidth]{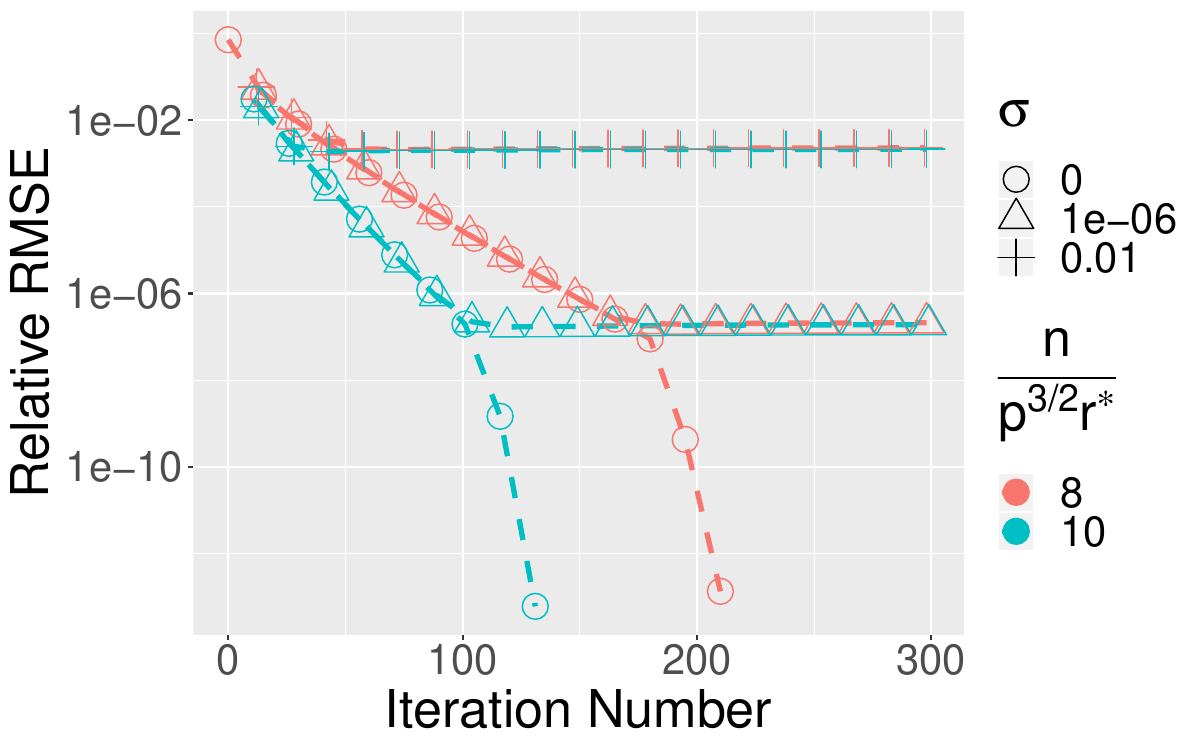}
			}\qquad 
			\subfigure[Scalar-on-tensor regression: RGN]{\includegraphics[width=0.46\textwidth]{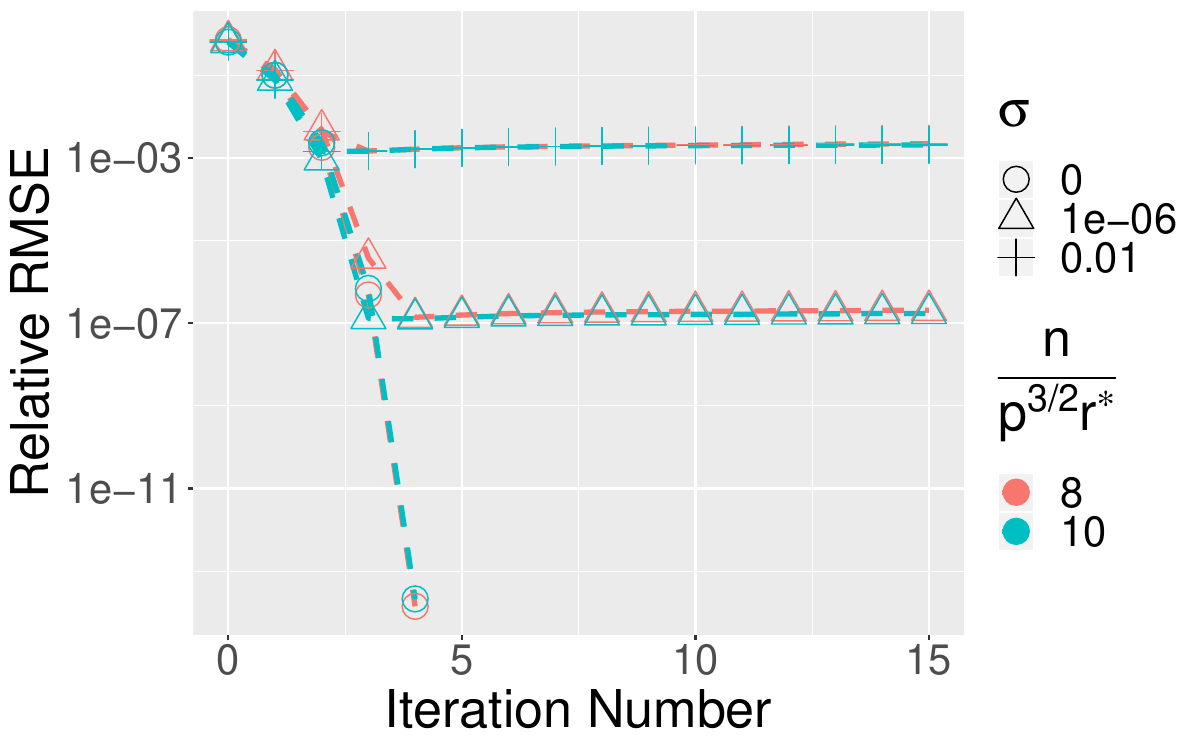}
			}
			\caption{Convergence performance of RGD/RGN in over-parameterized scalar-on-tensor regression with spectral initialization. Here, $p = 30, r^* = 3, r= 10$.}\label{fig: RGD/RGN-scalar-tensor-prop}
		\end{figure}

\begin{figure}[ht]
			\centering
			\subfigure[Tensor-on-vector regression: RGD]{\includegraphics[width=0.46\textwidth]{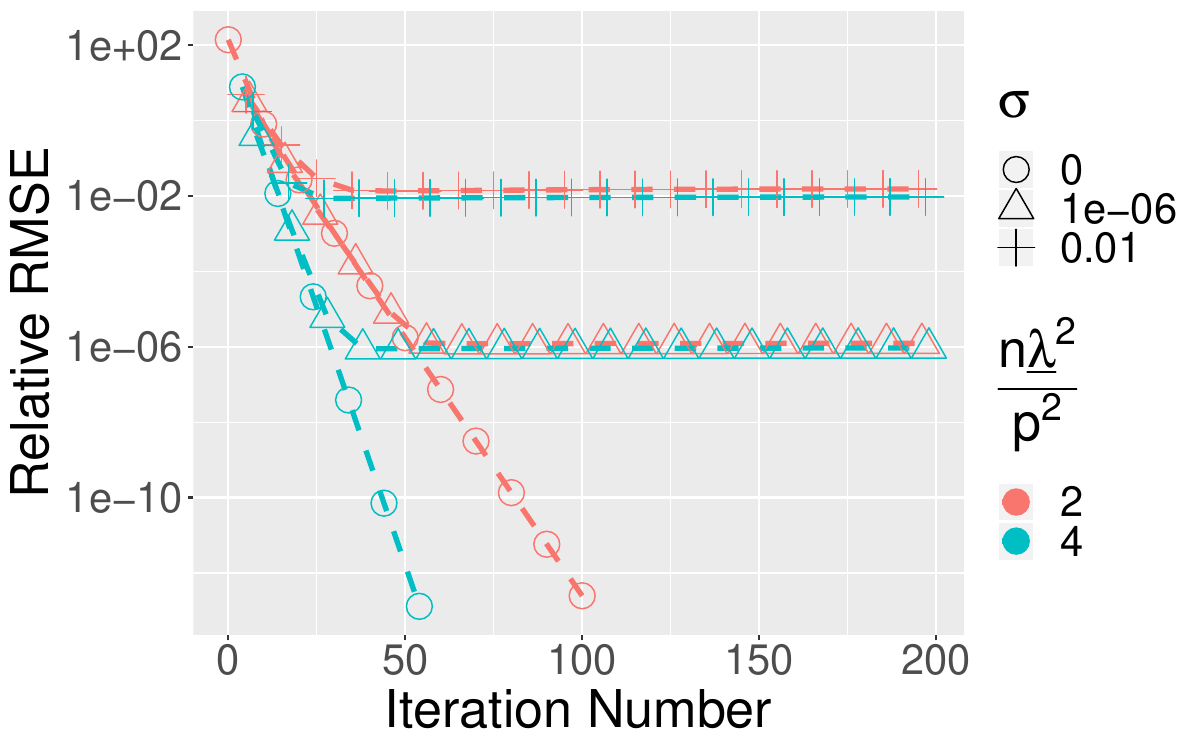}
			}\qquad 
			\subfigure[Tensor-on-vector regression: RGN]{\includegraphics[width=0.46\textwidth]{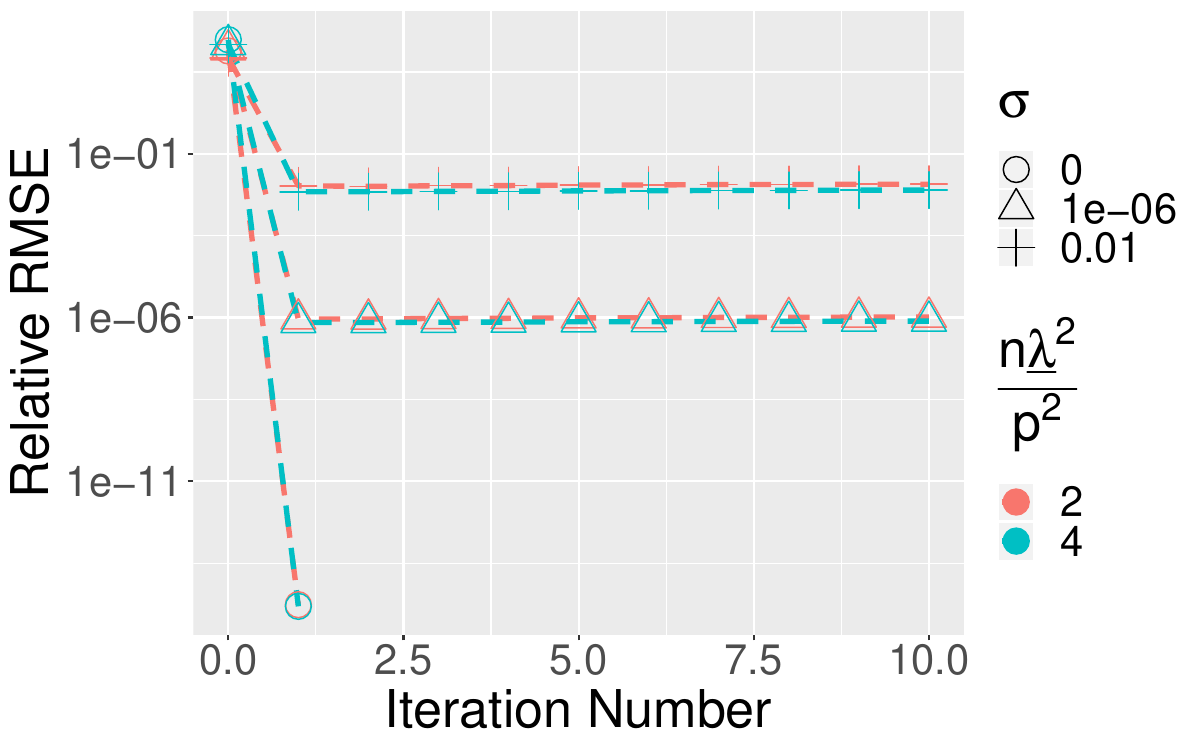}
			}
			\caption{Convergence performance of RGD/RGN in over-parameterized tensor-on-vector regression with spectral initialization. Here $p = 30, r^* = 3, r= 10$.}\label{fig: RGD/RGN-tensor-on-vector-prop}
		\end{figure}

		In this simulation, we examine the convergence rate of RGD/RGN in over-parameterized scalar-on-tensor regression and tensor-on-vector regression. We set $\sigma \in \{0,10^{-6}, 10^{-2}\}$, $p = 30$, $r^* = 3$, and $r = 10$. In scalar-on-tensor regression, we choose $n$ such that $ \frac{n}{p^{3/2} r^*} \in \{8,10 \}$; in tensor-on-vector regression, we let $ \frac{n\underline{\lambda}^2}{p^2 } \in \{2,4 \} $ where $\underline{\lambda} = \min_{k} \sigma_{r^*}(\cM_k( \bcS )) $. The convergence performance of RGD and RGN in scalar-on-tensor regression and tensor-on-vector regression are presented in Figures \ref{fig: RGD/RGN-scalar-tensor-prop} and \ref{fig: RGD/RGN-tensor-on-vector-prop}, respectively. In both examples, we find the estimation error of RGD converges linearly to the minimum precision in the noiseless setting and converges linearly to a limit determined by the noise level in the noisy setting. In scalar-on-tensor regression, we find RGN converges quadratically and in tensor-on-vector regression, we observe RGN converges with almost one iteration. We tried several other simulation settings and observed the similar phenomenon.
		
		\subsection{Effect of Input Rank and Sample Size on the Performance of RGD and RGN} \label{sec: numerics-input-rank-effect}
\begin{figure}[ht]
			\centering
			\subfigure[RGD]{\includegraphics[width=0.46\textwidth]{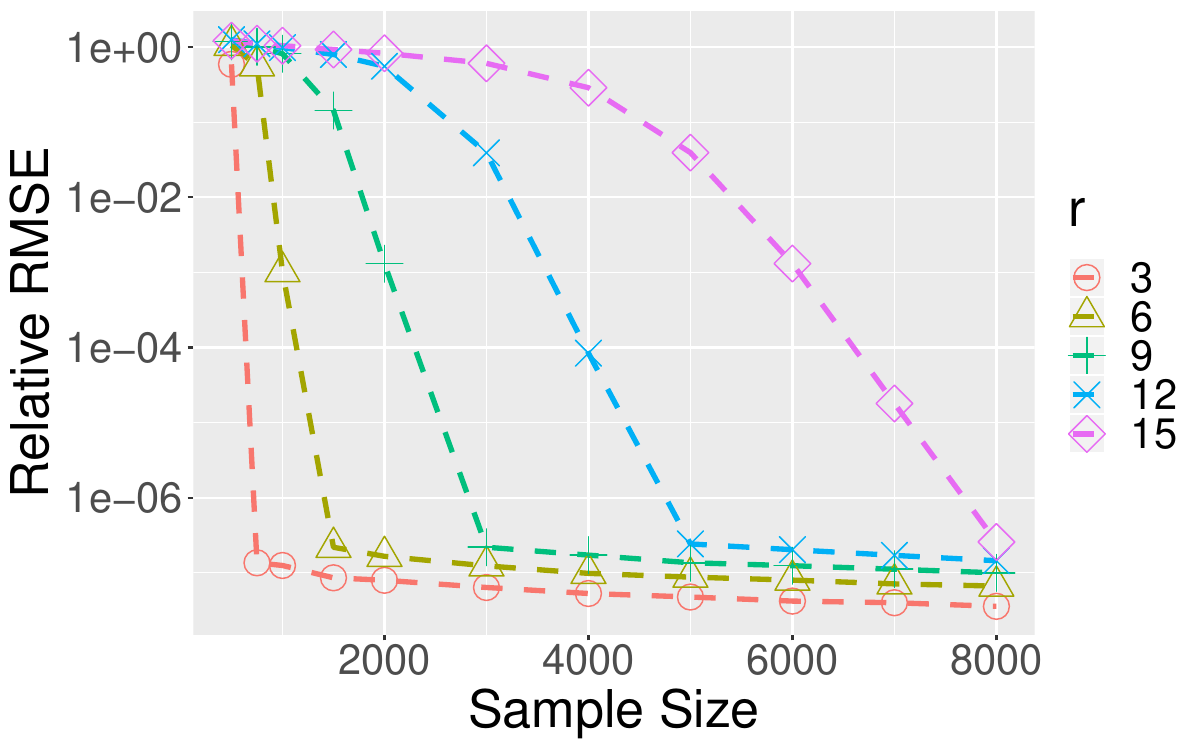}
			}\qquad 
			\subfigure[RGN]{\includegraphics[width=0.46\textwidth]{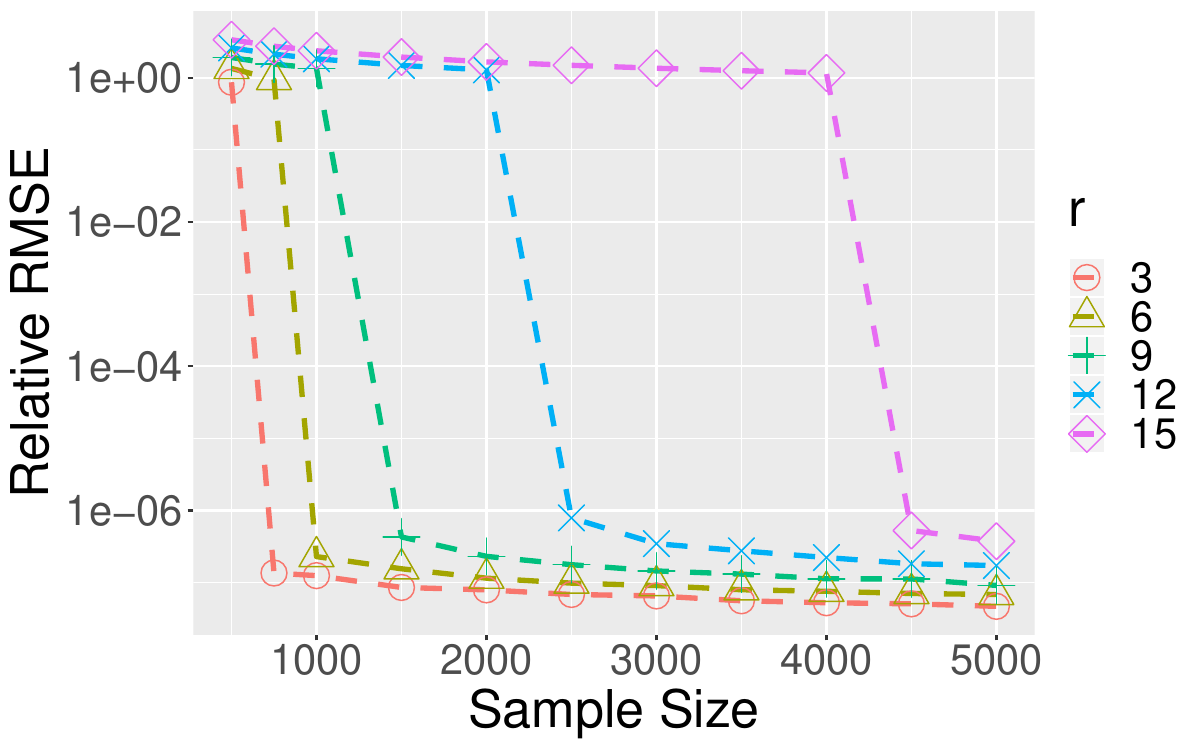}
			}
			\caption{Convergence performance of RGD/RGN in over-parameterized scalar-on-tensor regression with spectral initialization. Here $p = 30, r^* = 3,n \in [500,8000], r\in \{3,6,9,12,15 \}$.}\label{fig: RGD/RGN-scalar-tensor-diffrank}
		\end{figure}

\begin{figure}[ht]
			\centering
			\subfigure[RGD]{\includegraphics[width=0.46\textwidth]{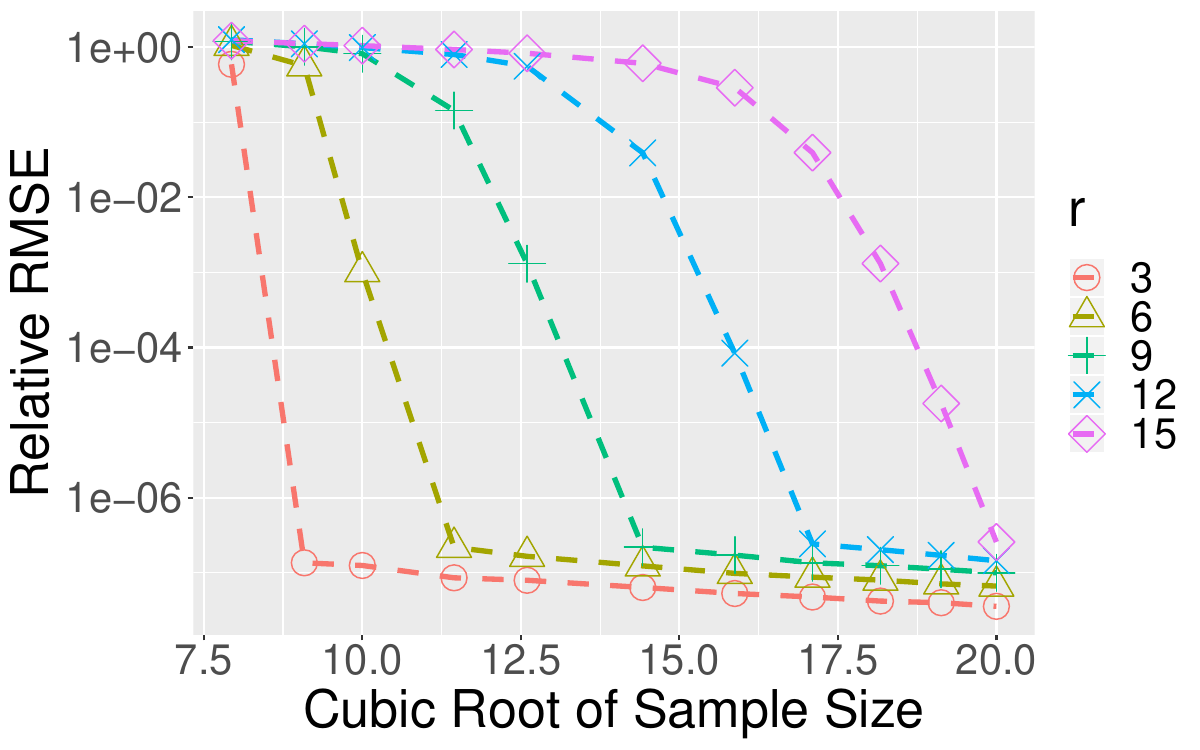}
			}\qquad 
			\subfigure[RGN]{\includegraphics[width=0.46\textwidth]{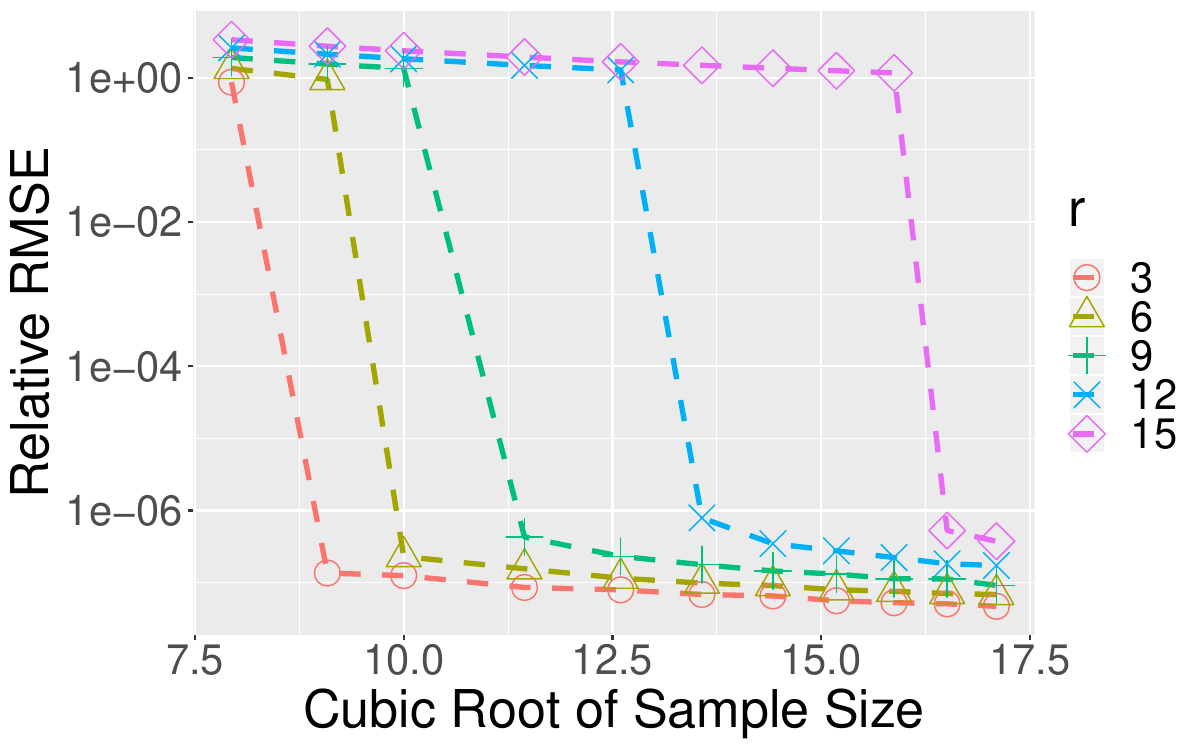}
			}
			\caption{Rescaled plot for the convergence performance of RGD/RGN in over-parameterized scalar-on-tensor regression with spectral initialization. Here $p = 30, r^* = 3,n \in [500,8000], r\in \{3,6,9,12,15 \}$.}\label{fig: RGD/RGN-scalar-tensor-diffrank_rescale}
		\end{figure}

		We also examine the effect of input rank $r$ and sample size $n$ on the convergence of RGD and RGN and we focus on the scalar-on-tensor regression example. We let $p = 30, r^* = 3$, $\sigma = 10^{-6}$, $n\in [500,8000]$ and input rank $r\in \{3,6,9,12,15\}$. The performance of RGD and RGN in this simulation study is given in Figure \ref{fig: RGD/RGN-scalar-tensor-diffrank}. We can see that for both RGD and RGN, the sample size requirement for convergence increases as the input rank $r$ increases. For a fixed $n$, the relative RMSE attainable by RGD and RGN increases as the input rank increases. In addition, the phase transition on the sample complexity for the failure/success in RGN is sharper than the one in RGD. This is because RGN enjoys a higher-order convergence compared to RGD and RGD converges slowly when the number of samples is around the threshold. This matches our main theoretical results in Sections \ref{sec:theory} and \ref{sec: applications-initializations}. Moreover, our results suggest that the number of samples needed for the convergence of RGD and RGN increases at the scale of $r^d$ for large $r$ (here $d=3$) and this is indeed suggested in Figure \ref{fig: RGD/RGN-scalar-tensor-diffrank_rescale} after we plot the cubic root of the sample size with respect to Relative RMSE.
		
		\subsection{Scalar-on-tensor Regression versus Matrix Trace Regression under Over-parameterization} \label{sec: bless-comp-limit}
		In this simulation, we compare the sample size requirements to ensure successful recovery in over-parameterized scalar-on-tensor regression and matrix trace regression with an increasing input rank via RGD. We focus on the noiseless setting, i.e., $\sigma = 0$. We say an algorithm achieves successful recovery if the averaged relative root mean squared error (Relative RMSE) $\|\bcX^{t} - \bcX^*\|_{\F}/ \|\bcX^*\|_{\F}$ is smaller than $0.01$. In scalar-on-tensor regression, we set $p = 90, r^*=1$, $r \in [1,\ldots,8]$, $n = [800,900,\ldots,3500]$ and in the matrix trace regression, we set $p = 100, r^{*} = 1$, $r \in [1,\ldots,8]$ and $n = [200,\ldots,3000]$. For every input rank $r$, we increase the sample size by $100$ at each time from the one that ensures the successful recovery with input rank $r-1$ until RGD succeeds. 
		
\begin{figure}[ht]
			\centering
			\includegraphics[width=\textwidth]{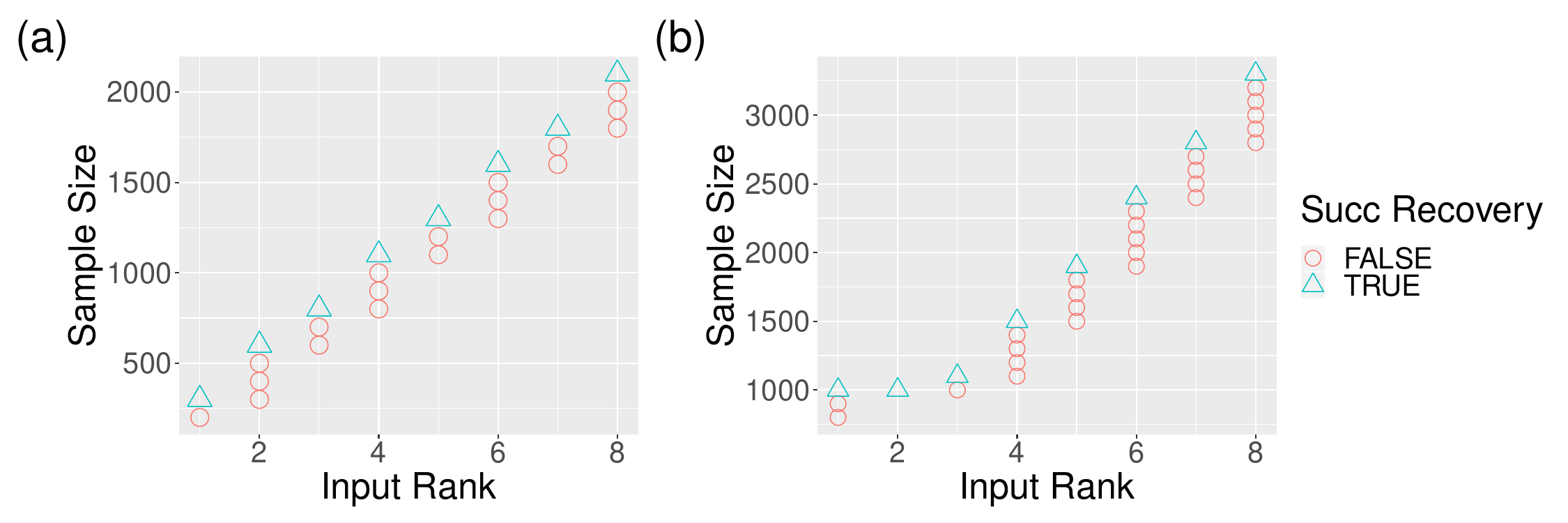}
			\caption{Comparison of successful recovery of RGD under over-parameterized matrix trace regression (Panel (a)) and scalar-on-tensor regression (Panel (b)).
			}
			\label{fig: sample-compare-matrix-tensor_succ}
		\end{figure}

		Figure \ref{fig: sample-compare-matrix-tensor_succ} shows as the input rank increases, the line of triangles for the sample size requirement of successful recovery in scalar-on-tensor regression is flat at the beginning stage while increases for large input $r$. In contrast, the sample size requirement for successful recovery of RGD in the matrix trace regression always increases linearly as input rank increases. This matches our theoretical findings in Section \ref{sec: computational-limit} that a ``free lunch'' on the sample complexity appears in over-parameterized scalar-on-tensor regression, but not in the matrix trace regression. Meanwhile, Figure \ref{fig: sample-compare-matrix-tensor_succ} shows when the input rank is equal to $r^*$, the phase transitions on sample complexity for the failure/success of RGD in matrix trace regression and scalar-on-tensor regression appear around $n = 300 \approx 2pr^*$ and $n = 1000 \approx p^{3/2}r^*$, respectively. This matches the results in Section \ref{sec: computational-limit} that there is a statistical-computational gap in scalar-on-tensor regression and $\Omega(p^{d/2})$ (here $d = 3$) samples are needed for any polynomial-time algorithm to succeed. 
		
		\subsection{Comparison of Riemannian Optimization Methods with Existing Algorithms} \label{sec: numerical comparison}
		
\begin{figure}[h]
			\centering
			\includegraphics[width=\textwidth]{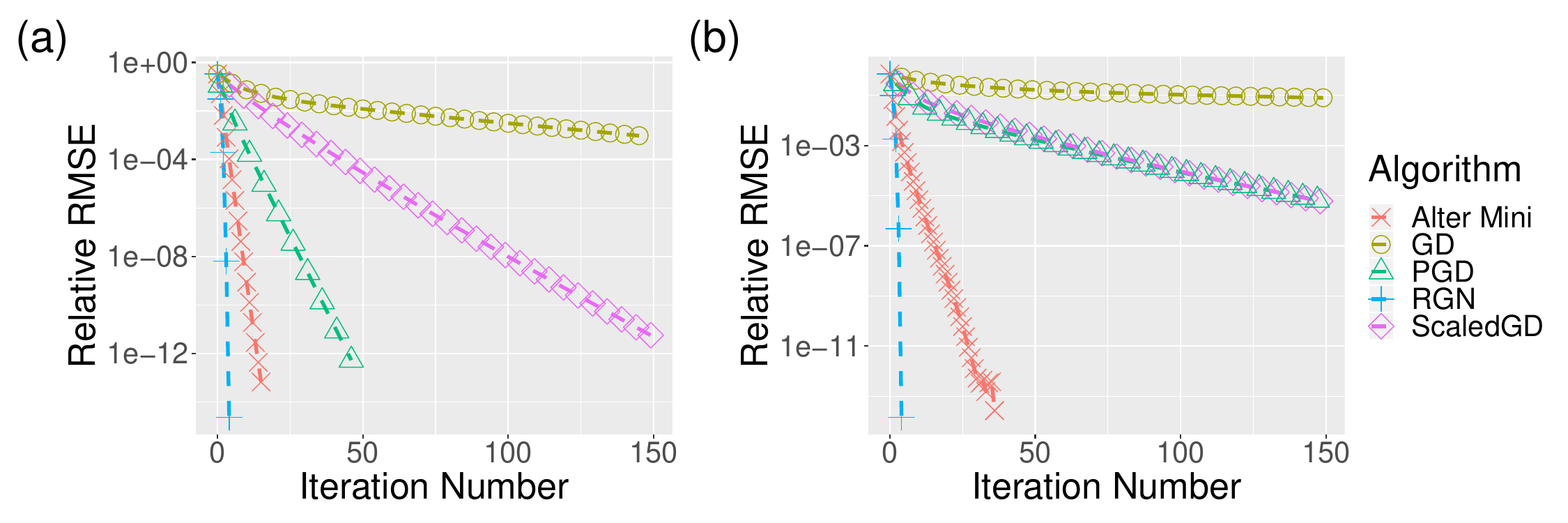}
			\caption{Panel (a): $r = 3$; Panel (b): $r = 10$. Relative RMSE of RGN (this work), alternating minimization (Alter Mini), projected gradient descent (PGD), gradient descent (GD), and scaled gradient descent (ScaledGD) in noiseless scalar-on-tensor regression.
			} \label{fig: tenreg comparison noiseless}
		\end{figure}

		In the second simulation, we compare RGN with other existing algorithms, including alternating minimization (Alter Mini) \citep{zhou2013tensor}, projected gradient descent (PGD)\citep{rauhut2017low}, gradient descent (GD) \citep{han2020optimal} and scaled gradient descent \citep{tong2021scaling}, in both exact and over-parameterized scalar-on-tensor regression. While implementing PGD, GD, and scaled GD, we evaluate five choices of step size, $\frac{1}{n} \cdot \{0.1,0.25, 0.5,0.75, 1 \}$, then choose the best one following \cite{zheng2015convergent}. We set $p = 30, r^* = 3, r \in \{3,10\}, n = 8p^{3/2}r^*$ and consider the noiseless case ($\sigma = 0$). Figure \ref{fig: tenreg comparison noiseless} shows RGN converges quadratically in both settings, while the other baseline algorithms converge at a  much slower linear rate. Moreover, when we go from exact-parameterization (Panel (a)) to over-parameterization (Panel (b)), the convergence rate of all baseline algorithms slows down significantly while RGN maintains its robust and fast second-order convergence performance. 
		
		\section{Conclusion and Discussions} \label{sec: conclusion}
		In this work, we propose Riemannian gradient descent and Riemannian Gauss-Newton methods for solving the general tensor-on-tensor regression. We provide optimal statistical and computational guarantees for these algorithms in both rank correctly-specified and overspecified settings and discover an intriguing blessing of the statistical-computational gap in the over-parameterized scalar-on-tensor regression. Our current initialization and computational results are established for several representative examples. It is of great interest to see whether these results can be extended to the general tensor-on-tensor regression problem. Moreover, the rank overspecification studied in this paper falls in the moderate over-parameterized regime in the sense that the model still includes more samples than the degree of freedom of parameters. It is interesting to consider the highly over-parameterized regime and study the analogy of implicit regularization effect \citep{gunasekar2017implicit,li2018algorithmic} in factorization formulated tensor problems. Some progress has been made recently in the tensor decomposition setting \citep{razin2021implicit,ge2021understanding}.
		
		\section*{Acknowledgements}
		
		The authors would like to thank Ilias Diakonikolas and Daniel Kane for helpful discussions. Diakonikolas and Kane developed a computational lower bound in the Statistical Query model (which further yields a low-degree polynomial computational lower bound) for low-rank scalar-on-tensor rank-one regression before this work; and the proof was later incorporated into a full paper in \cite{diakonikolas2022statistical}. However, the low-degree polynomial computational lower bound in Theorem \ref{th: lower-degree-polynomial-tensor-regression} of this paper is tighter and its proof is direct and arguably simpler. We also thank the Editor, the Associated Editor, and two anonymous referees for their helpful suggestions, which helped improve the presentation and quality of this paper.
		
		\vskip 0.2in
		\bibliographystyle{apalike}
		\bibliography{reference}

		\newpage
		\appendix
		
			\begin{center}
			{\LARGE Supplement to "Tensor-on-Tensor Regression: Riemannian Optimization, Over-parameterization, Statistical-computational Gap, and Their Interplay"	
				
			}		
			\medskip
			
			{\large Yuetian Luo \quad and\quad Anru R. Zhang}
			\medskip
		\end{center}
		
		In this supplement, we provide a table of contents, detailed algorithms, and all technical proofs.
		\tableofcontents
		
		\section{T-HOSVD and ST-HOSVD} \label{sec: HOSVD, STHOSVD}
		
		In this section, we present the procedures of truncated HOSVD (T-HOSVD) \citep{de2000multilinear} and sequentially truncated HOSVD (ST-HOSVD) \citep{vannieuwenhoven2012new}. For simplicity, we present the sequentially truncated HOSVD with the truncation order from mode $1$ to mode $(d+m)$.
		\begin{algorithm}[htbp]
			\caption{Truncated High-order Singular Value Decomposition (T-HOSVD)}
			\begin{algorithmic}[1]
			\State \noindent {\bf Input}: $\bcY \in \bbR^{p_{1} \times \cdots \times p_{d+m}}$, Tucker rank $\br = (r_1,\ldots, r_{d+m})$.
				\State Compute $\U_k^0 = \SVD_{r_k}(\cM_{k}(\bcY))$ for $k=1, \ldots,d+m$.
			\State \noindent {\bf Output}: $\widehat{\bcY} = \bcY \times_{k=1}^{d+m} P_{\U_k^0}$.
			\end{algorithmic}
			\label{alg: t-HOSVD}
		\end{algorithm}
		
		\begin{algorithm}[htbp]
			\caption{Sequentially Truncated High-order Singular Value Decomposition (ST-HOSVD) }
			\begin{algorithmic}[1]
				\State \noindent {\bf Input}: $\bcY \in \bbR^{p_{1} \times \cdots \times p_{d+m}}$, Tucker rank $\br = (r_1,\ldots, r_{d+m})$.
				\State Compute $\U_1^0 = \SVD_{r_1}(\cM_{1}(\bcY))$.
				\For{$k = 2, \ldots, d+m$}
				\State Compute  $\U_k^0 = \SVD_{r_k}(\cM_{k}(\bcY \times^{k-1}_{l=1} P_{\U_l^0} ))$.
				\EndFor
				\State \noindent {\bf Output}: $\widehat{\bcY} = \bcY \times_{k=1}^{d+m} P_{\U_k^0}$.
			\end{algorithmic}
			\label{alg: st-HOSVD}
		\end{algorithm}

		Different from the low-rank matrix projection, which can be efficiently and exactly computed via truncated SVD, performing low-rank tensor projection exactly, even for $\br=1$, can be NP-hard in general. We thus introduce the following quasi-projection property and the approximation constant $\delta(d+m)$. 
		
		\begin{Definition}[Quasi-projection of $\cH_{\br}$ and Approximation Constant $\delta(d+m)$]\label{def: quasi-projection Hr}
			Let $P_{\bbM_{\br}}(\cdot)$ be the projection map from $\bbR^{p_1 \times \cdots \times p_{d+m}}$ to the tensor space of Tucker rank at most $\br:= (r_1,\ldots,r_{d+m})$, i.e., for any $\bcZ\in \mathbb{R}^{p_1\times\cdots\times p_{d+m}}$ and $\widehat{\bcZ}$ of Tucker rank at most $\br$, one always has $\|\bcZ - \widehat{\bcZ}\|_{\F}\geq\|\bcZ - P_{\bbM_{\br}}(\bcZ)\|_{\F}$. 
			
			We say $\cH_{\br}$ satisfies the quasi-projection property with approximation constant $\delta(d+m)$ if $\|\bcZ - \cH_{\br}(\bcZ)\|_{\F} \leq \delta(d+m)\|\bcZ - P_{\bbM_{\br}}(\bcZ) \|_{\F}$ for any $\bcZ \in \bbR^{p_1 \times \cdots \times p_{d+m}}$. 
		\end{Definition}
		It is known that T-HOSVD and ST-HOSVD satisfy the \emph{quasi-projection property} (Chapter 10 in \cite{hackbusch2012tensor}). 
		\begin{Proposition}[Quasi-projection property of T-HOSVD and ST-HOSVD] \label{prop: quasi-projection}
			T-HOSVD and ST-HOSVD satisfy the quasi-projection property with the approximation constant $\delta(d+m)=\sqrt{d+m}$. That is for any $\bcZ \in \bbR^{p_1 \times \cdots \times p_{d+m}}$,
			$$\|\bcZ - \cH_{\br}(\bcZ)\|_{\F} \leq \sqrt{d+m}\|\bcZ - P_{\bbM_{\br}}(\bcZ) \|_{\F}.$$
			Here $\cH_{\br}$ is either T-HOSVD or ST-HOSVD.
		\end{Proposition}
		We note in the matrix setting by taking the rank $r$ truncated SVD as the retraction operator, the approximation constant in Definition \ref{def: quasi-projection Hr} is $1$, and this fact is used in Corollary \ref{coro: local contraction general setting-RGD-RGN-matrix}.

		\section{Additional Notation and Preliminaries} \label{sec: additional-notation}
		Let $\bbN = \{0,1,2,\ldots \}$ be the set of natural numbers. For any $\U\in \mathbb{O}_{p, r}$, we use $\U_\perp\in \mathbb{O}_{p, p-r}$ to represent the orthonormal complement of $\U$. The matricization operation $\mathcal{M}_k(\cdot)$ unfolds an order-$d$ tensor along mode $k$ to a matrix, say $\bcA\in\mathbb{R}^{p_1\times \cdots \times p_d}$ to $\mathcal{M}_k(\bcA)\in \mathbb{R}^{p_k\times p_{-k}}$, where $p_{-k} = \prod_{j\neq k}p_j$. Specifically, $\mathcal{M}_k(\bcA) \in \mathbb{R}^{p_k\times p_{-k}},$
		\begin{equation}\label{eq: tensor-matricization}
			\left(\mathcal{M}_k(\bcA)\right)_{\left[i_k, j\right]} =\bcA_{[i_1,\ldots, i_d]}, \quad j = 1 + \sum_{\substack{l=1\\l\neq k}}^d\left\{(i_l-1)\prod_{\substack{z=1\\z\neq k}}^{l-1}p_z\right\}
		\end{equation}
		for any $1\leq i_l \leq p_l, l=1,\ldots, d$. We also denote $\cT_k(\cdot)$ as the mode-$k$ tensorization, i.e., the reverse operator of $\cM_k(\cdot)$: $\cT_k(\cM_k(\bcA)) = \bcA$ for any $k =1,\ldots,d$ The mode-$k$ product of $\bcA \in \mathbb{R}^{p_1 \times \cdots \times p_d}$ with a matrix $\B\in \mathbb{R}^{r_k\times p_k}$, denoted by $\bcA \times_k \B$, is a $p_1 \times \cdots \times p_{k-1}\times r_k \times p_{k+1}\times \cdots \times p_d$-dimensional tensor, defined as
		\begin{equation}\label{eq: tensor-matrix product}
			(\bcA \times_k \B)_{[i_1, \ldots, i_{k-1}, j, i_{k+1}, \ldots, i_d]} = \sum_{i_k=1}^{p_k} \bcA_{[i_1, i_2, \ldots, i_d]} \B_{[j, i_k]}.
		\end{equation} The inner product $\langle \cdot, \cdot \rangle$ of any two tensors is defined as $\langle\bcA, \bcB\rangle = \sum_{i_1,\ldots, i_d} \bcA_{[i_1,\ldots, i_d]}\bcB_{[i_1,\ldots, i_d]}$. The following equality connects the tensor-matrix product and matricizations \cite[Section 4]{kolda2001orthogonal}: 
		\begin{equation}\label{eq: matricization relationship}
			\begin{split}
				\mathcal{M}_k \left( \bcS\times_1 \U_1 \times \cdots \times_d \U_d \right)
				= \U_k \mathcal{M}_k(\bcS) (\U_d^\top \otimes \cdots \otimes \U_{k+1}^\top \otimes \U_{k-1}^\top \otimes \cdots \otimes \U_{1}^\top),
			\end{split}
		\end{equation}
		where ``$\otimes$'' is the matrix Kronecker product. Recall the contracted tensor inner product is defined as follows
		\begin{equation}\label{def: contracted-tensor-inner-product}
			\left(\langle \bcA_i, \bcX^* \rangle_* \right)_{[j_1,\ldots, j_m]} = \sum_{\substack{k_l=1,\\l=1,\ldots,d}}^{p_l}  \bcA_{i[k_1,\ldots,k_d]}  \bcX^*_{[k_1,\ldots,k_d,j_1,\ldots,j_m]}.
		\end{equation}
		
		Next, we provide the explicit formulas for the tangent space of $\bbM_{\br}$ at $\bcX$ \citep{koch2010dynamical,luo2021low}:
		\begin{equation}\label{eq: tangent-space-representation}
			\begin{split}
				& T_{\bcX} \bbM_{\br} = \left\{  \bcB \times_{k=1}^{d+m} \U_k + \sum_{k=1}^{d+m} \cT_k(\U_{k\perp}\D_k \W_k^{\top}):\begin{array}{l}
					\bcB \in \bbR^{r_1 \times \cdots \times r_{d+m}}, \D_k \in \bbR^{(p_k - r_k)\times r_k},\\
					k=1,\ldots,(d+m)
				\end{array}\right\}\\
				= & \left\{\bcB \times_{k=1}^{d+m} \U_k + \sum_{k=1}^{d+m} \bcS \times_k  \U_{k\perp} \D_k \times_{j \neq k} \U_j : \begin{array}{l}
					\bcB \in \bbR^{r_1 \times \cdots \times r_{d+m}}, \D_k \in \bbR^{(p_k-r_k)\times r_k}, \\ k=1,\ldots,d+m
				\end{array}\right\},
			\end{split}
		\end{equation} 
		where $\cT_k(\cdot)$ is the mode-$k$ tensorization operator and $\W_k$ is given in \eqref{def: Wk}.

\section{Rank Selection} \label{sec: rank-selection}
To ensure a small estimation error of RGD and RGN, it is important to select a parsimonious rank $\br$, while also ensuring that $\br \geq \br^*$. In this section, we further introduce a data-driven method for choosing $\br$. Suppose $n$ is even for simplicity of presentation. One straightforward idea is to use the first half of the samples $\{\bcY_i, \bcA_i \}_{i=1}^{n/2}$ to get an estimator of $\bcX^*$ under different input rank and then compute the prediction error of the estimator on the second half of the data. Then we can choose the rank that minimizes the out-of-sample prediction error. 

In the rest of this section, we explore the empirical performance of this method in the scalar-on-tensor regression setting considered in Section \ref{sec: numerical RGD/RGN property} with $p = 30, r^* = 5, \sigma \in \{0.01,0.1,1\}$ and $n = 5000$. We consider different input ranks in $\{2,\ldots,9\}$ and let $\widehat{\bcX}_r$ denote the estimator returned by applying RGD or RGN on the first half of the data with input rank $r$. In the following Figure \ref{fig: rank-selection}, we plot the our-of-sample prediction error of the estimator $\widehat{\bcX}_r$ on the holdout data, i.e., $\sqrt{\sum_{i=n/2+1}^{n}(\bcY_i - \langle \bcA_i, \widehat{\bcX}_r \rangle_* )^2}$, versus different input rank $r$. We can see that the out-of-sample prediction error consistently achieves its minimum at the true rank value.  
		\begin{figure}[ht]
			\centering
			\subfigure[RGD]{\includegraphics[width=0.46\textwidth]{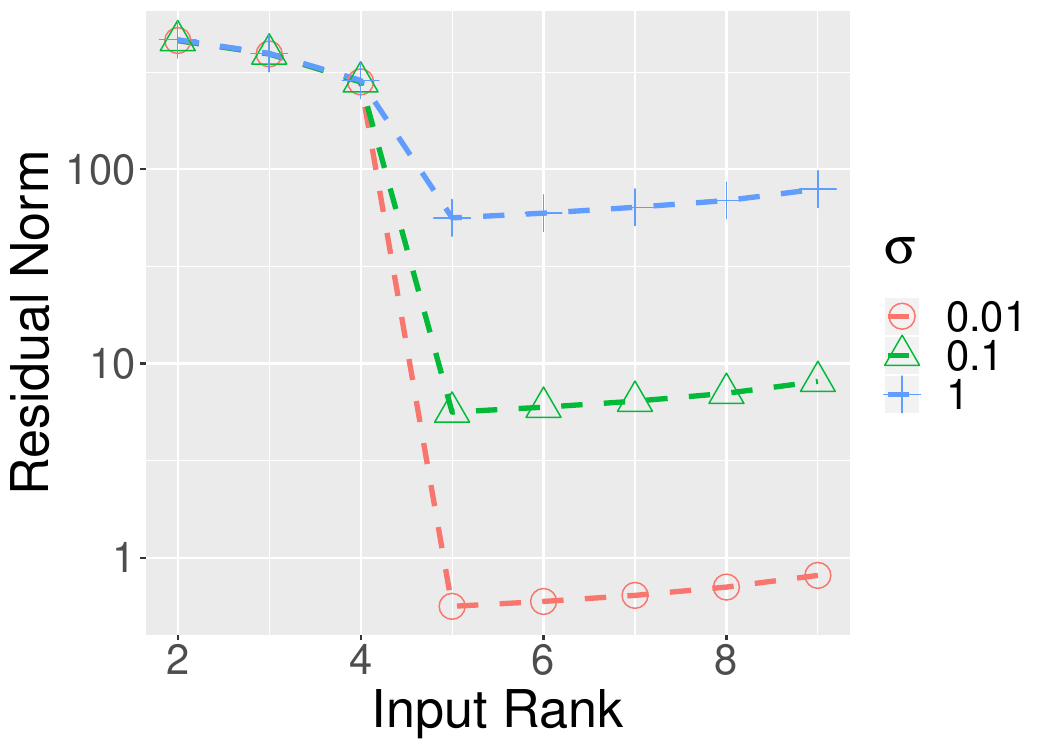}
			}\qquad 
			\subfigure[RGN]{\includegraphics[width=0.46\textwidth]{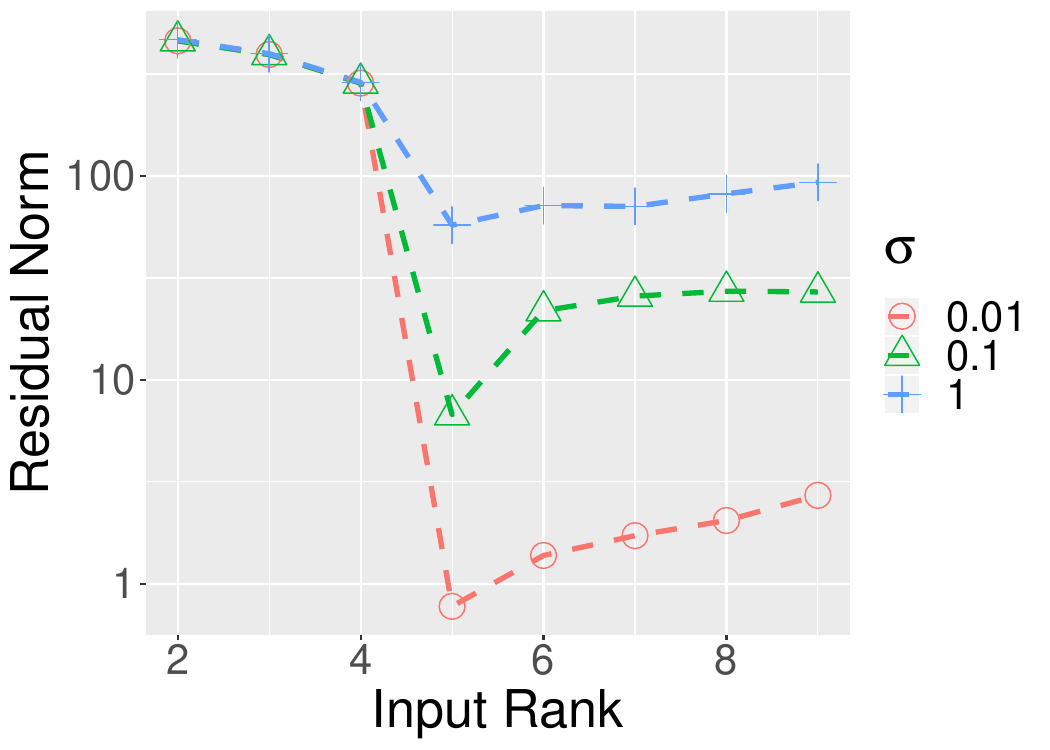}
			}
			\caption{Out-of-sample prediction error with different input rank in training. Here $p = 30, r^* = 5,n = 5000,\sigma \in \{0.01,0.1,1\}, r \in \{2,3,4,5,6,7,8,9\}$.}\label{fig: rank-selection}
		\end{figure}

		\section{Proofs in Section \ref{sec: Algorithm} }\label{sec:proof-algorithm}
		\subsection{Proof of Lemma \ref{lm: gradient}}\label{sec:proof-lm-gradient}
		Since $\bbM_{\br}$ is an embedded submanifold of $\bbR^{p_1\times \cdots \times p_{d+m}}$ and the Euclidean gradient of $f$ in the ambient space is $ \nabla f(\bcX) = \scA^*( \scA(\bcX) - \bcY )  $, the result follows from \cite[(3.37)]{absil2009optimization}. \quad $\blacksquare$
		\subsection{Proof of Lemma \ref{lm: gauss-newton}}\label{sec:proof-lm-gauss-newton}
		In our objective \eqref{eq:minimization}, the RGN update $\eta^{\RGN} \in T_{\bcX^t} \bbM_{\br}$ should solve the following RGN equation \citep[Chapter 8.4]{absil2009optimization},
		\begin{equation} \label{eq: Riemannian Gauss-newton equation}
			-\grad f(\bcX^t)= P_{T_{\bcX^t}} \left( \scA^* (\scA (\eta^{\RGN})) \right).
		\end{equation}
		In view of the Riemannian Gauss-Newton equation in \eqref{eq: Riemannian Gauss-newton equation} and the Riemannian gradient in Lemma \ref{lm: gradient}, to prove the claim, we only need to show
		\begin{equation} \label{eq: Thm Rie G-N need to proof}
			P_{T_{\bcX^t}}(\scA^*(\scA( \bcX^t + \eta^{\RGN}) -\bcY)) = 0.
		\end{equation} 
		
		From the optimality condition of the least squares problem $\min_{\eta \in T_{\bcX^t} \bbM_{\br} } \frac{1}{2}\|\bcY - \scA P_{T_{\bcX^t}}(\bcX^t + \eta)\|_2^2$, we have
		\begin{equation*}
			P_{T_{\bcX^t}}\scA^*\left(\scA P_{T_{\bcX^t}} (\bcX^t + \eta^{\RGN})
			- \bcY\right) = P_{T_{\bcX^t}}\scA^*\left(\scA (\bcX^t + \eta^{\RGN})
			- \bcY\right)=0.
		\end{equation*} This finishes the proof. \quad $\blacksquare$

		\section{Proofs in Section \ref{sec:theory} }\label{sec:proof-theory}
		We begin by introducing a few preliminary results and then give the proof in subsections. The following lemma shows that any tensor in $T_{\bcX}\bbM_{\br}$ is at most Tucker rank $2\br$
		\begin{Lemma}\label{lm: tangent vector rank 2r property}[\cite[Lemma 1]{luo2021low}] For any $\bcX \in \bbM_\br$, any tensor $\bcZ \in T_{\bcX}\bbM_{\br}$ is at most Tucker rank $2\br$.
		\end{Lemma}
		
		The following lemma provides a user-friendly perturbation bound for singular subspaces under perturbation and is critical in our analysis.
		\begin{Lemma} \label{lm: user-friendly-subspace-perturb}
			Let $\A$ be a rank $r$ matrix in $\bbR^{m\times n}$ with economic singular value decomposition (SVD) $\U \bSigma \V^\top $. Suppose $\B = \A + \Z \in \bbR^{m \times n}$ for some perturbation matrix $\Z$, and the top $r$ truncated SVD of $\B$ is given as $\widehat{\U} \widehat{\bSigma} \widehat{\V}^\top$. Then
			\begin{equation*}
				\max\left\{\|\widehat{\U} \widehat{\U}^\top - \U\U^\top\|, \|\widehat{\V} \widehat{\V}^\top - \V\V^\top)\|\right\} \leq \frac{2\| \Z \| }{\sigma_r(\A)}.
			\end{equation*}
		\end{Lemma}
		{\noindent \bf Proof}. First, notice $\|\widehat{\U} \widehat{\U}^\top - \U\U^\top\| = \|\widehat{\U}_\perp^\top \U\| = \| \sin \Theta (\widehat{\U}, \U) \| $ by \cite[Theorem 2.5.1]{golub2013matrix} and \cite[Lemma 1]{cai2018rate}. This lemma is a special case of \cite[Theorem 5]{luo2021schatten}.
		\quad $\blacksquare$
		
		Also for $\bcX \in \bbM_{\br}$ and the projector $P_{T_\bcX}(\cdot)$ in \eqref{eq: tangent space projector}, we let $P_{(T_{\bcX})_\perp}(\bcZ) := \bcZ - P_{T_{\bcX}}(\bcZ)$ be the orthogonal complement of the projector $P_{T_{\bcX}}$.
		The next lemma provides another useful upper bound which will be used frequently in our proof.
		\begin{Lemma}\label{lm: orthogonal-projection-bound}$\bcX^t$ is a Tucker rank $\br$ tensor and $\bcX^*$ is a Tucker rank $\br^*$ tensor. 
			Then $\| P_{T_{\bcX^t}} \scA^* \scA( P_{(T_{\bcX^t})_\perp} \bcX^* ) \|_{\F} \leq R_{2\br + \br^*} \| P_{(T_{\bcX^t})_\perp} \bcX^*\|_{\F}$. 
		\end{Lemma}
		{\noindent \bf Proof}.
		\begin{equation*}
			\begin{split}
				\| P_{T_{\bcX^t}} \scA^* \scA( P_{(T_{\bcX^t})_\perp} \bcX^* ) \|_{\F} =& \sup_{\bcZ: \|\bcZ\|_{\F}\leq 1} \langle P_{T_{\bcX^t}} \scA^* \scA( P_{(T_{\bcX^t})_\perp} \bcX^* ), \bcZ \rangle\\
				=& \sup_{\bcZ: \|\bcZ\|_{\F}\leq 1} \langle \scA( P_{(T_{\bcX^t})_\perp} \bcX^* ), \scA P_{T_{\bcX^t}} (\bcZ) \rangle\\
				\overset{(a)}\leq & \sup_{\bcZ: \|\bcZ\|_{\F}\leq 1}  R_{2\br + \br^*}\|P_{(T_{\bcX^t})_\perp} \bcX^*\|_{\F} \|P_{T_{\bcX^t}} (\bcZ) \|_{\F}\\
				\overset{(b)}\leq  & R_{2\br + \br^*}\|P_{(T_{\bcX^t})_\perp} \bcX^*\|_{\F}.
			\end{split}
		\end{equation*} Here (a) is due to Lemma \ref{lm:retricted orthogonal property}, $\langle P_{(T_{\bcX^t})_\perp} \bcX^*, P_{T_{\bcX^t}} (\bcZ) \rangle = 0$, $P_{(T_{\bcX^t})_\perp} \bcX^*$ and $P_{T_{\bcX^t}} (\bcZ)$ are of Tucker rank at most $\br^*$ and $2\br$, respectively; (b) is because $\|P_{T_{\bcX^t}}(\bcZ) \|_{\F} \leq \|\bcZ\|_{\F} \leq 1$. 
		\quad $\blacksquare$
		
		Next, we begin the proof for the results in Section \ref{sec:theory} one by one.
		
		\subsection{Proof of Proposition \ref{prop: TRIP-under-Gaussian-design}}
		To prove the statement, we need to show when $n \geq C(\sum_{i=1}^d (p_i-r_i)r_i + \prod_{i=1}^d r_i )\log(d)/R^2_{\br} $, for any tensor $\bcZ$ of Tucker rank at most $\br$, we have $(1-R_{\br}) \|\bcZ\|^2_{\F} \leq \|\scA(\bcZ)\|_\F^2 \leq (1+R_{\br}) \|\bcZ\|_{\F}^2$ holds with probability at least $1 - \exp(-c(\sum_{i=1}^d p_i))$.
		
		Denote $\bcZ^{(j_1,\ldots,j_m)} := \bcZ_{[:,\ldots,:, j_1,\ldots,j_m]} \in \bbR^{p_{d+1} \times \cdots \times p_{d+m}}$ for any $j_1 \in [p_{d+1}], \cdots, j_m \in [p_{d+m}]$. Since $\bcZ$ is of Tucker rank at most $\br$ and $\bcZ^{(j_1,\ldots,j_m)}$ can be rewritten as $\bcZ \times_{d+1} \e_{j_1}^\top \cdots \times_{d+m} \e^\top_{j_m} $ where $\e_i$ denotes the standard $i$-th base vector, we have $\bcZ^{(j_1,\ldots,j_m)}$ is of Tucker rank at most $\br' = (r_1,\ldots, r_d)$. Let us define another linear map $\scA': \bbR^{p_{d+1} \times \cdots \times p_{d+m}} \to \bbR^n$ such that $\scA'(\bcZ^{(j_1,\ldots,j_m)})_{[i]} = \langle \bcA_i, \bcZ^{(j_1,\ldots,j_m)} \rangle $. It is easy to check $\|\scA(\bcZ)\|_{\F}^2 = \sum_{j_1 \in [p_{d+1}]} \cdots \sum_{j_m \in [p_{d+m}]}  \| \scA'(\bcZ^{(j_1,\ldots,j_m)}) \|_2^2$. 
		
		On the other hand, following the same proof of \cite[Theorem 2]{rauhut2017low} by replacing their Lemma 2 with a tighter covering number for the low Tucker rank tensor space given in Lemma \ref{lm: covering-number}, we have when $n \geq C(\sum_{i=1}^d (p_i-r_i)r_i + \prod_{i=1}^d r_i )\log(d)/R^2_{\br} $, then with probability at least $1 - \exp(-c (\sum_{i=1}^d p_i) )$, $(1-R_{\br}) \|\bcZ'\|^2_{\F} \leq \|\scA'(\bcZ')\|_2^2 \leq (1+R_{\br}) \|\bcZ'\|_{\F}^2$ holds for any tensor $\bcZ'$ of Tucker rank at most $\br'$. So 
		\begin{equation*}
			\begin{split}
				(1 - R_\br ) \|\bcZ\|^2_\F & =(1 - R_\br ) \sum_{j_1,\ldots,j_m} \|\bcZ^{(j_1,\ldots,j_m)}\|_\F^2 \leq \sum_{j_1,\ldots,j_m} \| \scA'(\bcZ^{(j_1,\ldots,j_m)}) \|_2^2\\
				&= \|\scA(\bcZ)\|_{\F}^2 \leq (1 + R_\br ) \sum_{j_1,\ldots,j_m} \|\bcZ^{(j_1,\ldots,j_m)}\|_\F^2 \\
				& = (1 + R_\br ) \|\bcZ\|^2_\F,
			\end{split}
		\end{equation*} holds with probability at least $1 - \exp(-c(\sum_{i=1}^d p_i))$. This finishes the proof of this proposition. \quad $\blacksquare$
		
		\subsection{Proof of Theorem \ref{th: local contraction general setting-RGD}}
		First, notice the convergence result in the noiseless setting follows easily from the noisy setting by setting $\bcE = 0$. We prove the convergence result in the noisy case. Recall $\bcX^{t+0.5} = \bcX^t - \alpha_t P_{T_{\bcX^t}} \scA^*( \scA(\bcX^t) - \bcY ) $ where  $\alpha_t = \|P_{T_{\bcX^t}} (\bcA^*( \bcA(\bcX^t) - \bcY ))\|_\F^2/\|\scA P_{T_{\bcX^t}} (\bcA^*( \bcA(\bcX^t) - \bcY ))\|_{\F}^2$ and $ \bcX^{t+1} = \cH_\br(\bcX^{t+0.5})$. Then
		\begin{equation} \label{ineq: RGD-Xt+1 -Xstar}
			\begin{split}
				&\|\bcX^{t+1} - \bcX^*\|_{\F}\\ 
				=& \| \cH_{\br}(\bcX^{t+0.5}) - \bcX^* \|_{\F}\\
				\leq & \| \cH_{\br}(\bcX^{t+0.5}) - \bcX^{t+0.5}\|_{\F} + \|\bcX^{t+0.5} - \bcX^*\|_{\F} \\
				\overset{(a)}\leq & \sqrt{d+m} \|P_{\bbM_{\br}}(\bcX^{t+0.5}) - \bcX^{t+0.5}\|_{\F} + \|\bcX^{t+0.5} - \bcX^*\|_{\F}\\
				\overset{(b)}\leq & (\sqrt{d+m} + 1)\|\bcX^{t+0.5} - \bcX^*\|_{\F}\\
				=  & (\sqrt{d+m} + 1)\| \bcX^t - \alpha_t P_{T_{\bcX^t}} \scA^*( \scA(\bcX^t) - \bcY ) - \bcX^*\|_{\F}\\
				\overset{(c)}=  & (\sqrt{d+m} + 1)\| \bcX^t - \alpha_t P_{T_{\bcX^t}} \scA^*( \scA(\bcX^t - \bcX^*) - \bcE ) - \bcX^*\|_{\F}\\
				\overset{(d)}\leq & (\sqrt{d+m} + 1) \Big ( \underbrace{\| P_{(T_{\bcX^t})_\perp} (\bcX^t - \bcX^*)\|_{\F}}_{(A1)} + \underbrace{\| (P_{T_{\bcX^t}}-  \alpha_t P_{T_{\bcX^t}} \scA^* \scA P_{T_{\bcX^t}}) (\bcX^t - \bcX^*)\|_{\F}}_{(A2)}\\
				& +  \alpha_t \underbrace{\| P_{T_{\bcX^t}} \scA^* \scA P_{(T_{\bcX^t})_\perp} (\bcX^t - \bcX^*)\|_{\F}}_{(A3)} + \alpha_t \underbrace{\|  P_{T_{\bcX^t}} \scA^*(\bcE )  \|_\F}_{(A4)} \Big).
			\end{split}
		\end{equation} here (a) is by the quasi-projection property of T-HOSVD and ST-HOSVD in Proposition \ref{prop: quasi-projection}; (b) is by the projection property of $P_{\bbM_{\br}}(\cdot)$; (c) is because $\bcY = \scA(\bcX^*) + \bcE$; (d) is by triangle inequality.
		
		Next, we bound $(A1), (A2), (A3)$ and $(A4)$ separately.
		\begin{itemize}
			\item $(A1) = \| P_{(T_{\bcX^t})_\perp} (\bcX^t - \bcX^*)\|_{\F} = \|P_{(T_{\bcX^t})_\perp} (\bcX^*)\|_{\F} \overset{ \text{Lemma } \ref{lm: orthogonal projection} } \leq \frac{2(d+m) \| \bcX^t - \bcX^* \|^2_{\F} }{\underline{\lambda}} $.
			\item \begin{equation*} 
				\begin{split}
					(A2) \leq & \| (P_{T_{\bcX^t}}-  \alpha_t P_{T_{\bcX^t}} \scA^* \scA P_{T_{\bcX^t}})\| \|\bcX^t - \bcX^*\|_{\F}\\
					\leq & \left( \| P_{T_{\bcX^t}}- P_{T_{\bcX^t}} \scA^* \scA P_{T_{\bcX^t}}\| + |1-\alpha_t| \|P_{T_{\bcX^t}} \scA^* \scA P_{T_{\bcX^t}}\|   \right) \|\bcX^t - \bcX^*\|_{\F} \\
					\overset{(a)}\leq & (R_{2\br} + \frac{R_{2\br}}{1-R_{2\br}} (1+ R_{2\br}) ) \|\bcX^t - \bcX^*\|_{\F}\\
					= & \frac{2R_{2\br}}{1-R_{2\br}} \|\bcX^t - \bcX^*\|_{\F},
				\end{split}
			\end{equation*} where $(a)$ is because $\|P_{T_{\bcX^t}} \scA^* \scA P_{T_{\bcX^t}}\| \leq 1+R_{2\br}$ by Lemma \ref{lm: spectral norm bound of Atop A}, $ (1+R_{2\br})^{-1} \leq \alpha_t \leq (1-R_{2\br})^{-1} \Longrightarrow |1-\alpha_t| \leq R_{2\br}/(1-R_{2\br}) $ by the TRIP assumption on $\scA$, and 
			\begin{equation}\label{ineq: Pxt-minus-PAAP}
				\begin{split}
					\| P_{T_{\bcX^t}}- P_{T_{\bcX^t}} \scA^* \scA P_{T_{\bcX^t}}\| &= \sup_{\bcZ: \|\bcZ\|_\F = 1} |\langle ( P_{T_{\bcX^t}}- P_{T_{\bcX^t}} \scA^* \scA P_{T_{\bcX^t}})(\bcZ), \bcZ \rangle|\\
					& = \sup_{\bcZ: \|\bcZ\|_\F = 1} \left| \|P_{T_{\bcX^t}}(\bcZ)\|_\F^2 - \|\scA P_{T_{\bcX^t}}(\bcZ)\|_{\F}^2  \right|\\
					& \overset{(a')}\leq R_{2\br}  \sup_{\bcZ: \|\bcZ\|_\F = 1} \|P_{T_{\bcX^t}}(\bcZ)\|_\F^2 \leq R_{2\br},
				\end{split}
			\end{equation} where $(a')$ is by the TRIP assumption on $\scA$.
			\item $(A3) = \| P_{T_{\bcX^t}} \scA^* \scA P_{(T_{\bcX^t})_\perp} (\bcX^*)\|_{\F} \overset{ \text{Lemma } \ref{lm: orthogonal-projection-bound} }\leq R_{2\br + \br^*} \|P_{(T_{\bcX^t})_\perp} (\bcX^*)\|_\F \overset{ \text{Lemma } \ref{lm: orthogonal projection} }  \leq \frac{2(d+m) R_{2\br + \br^*}  \| \bcX^t - \bcX^* \|^2_{\F} }{\underline{\lambda}} $.
			\item $(A4) \leq \|(\scA^*(\bcE ))_{\max(2\br)}\|_{\F}$ since $\|  P_{T_{\bcX^t}} \scA^*(\bcE )  \|_\F$ lies in $T_{\bcX^t} \bbM_\br$ and is of rank at most $2\br$ by Lemma \ref{lm: tangent vector rank 2r property}.
		\end{itemize}
		
		By plugging upper bounds of $(A1) - (A4)$ into \eqref{ineq: RGD-Xt+1 -Xstar}, we have
		\begin{equation} \label{ineq: RGD-iterative-bound}
			\begin{split}
				\|\bcX^{t+1} - \bcX^*\|_{\F} &\leq \frac{2(\sqrt{d+m} + 1)}{1-R_{2\br}}\left( R_{2\br} + \frac{(1 +R_{2\br + \br^*} - R_{2\br})(d+m)}{\underline{\lambda}} \|\bcX^{t} - \bcX^*\|_{\F} \right)\|\bcX^{t} - \bcX^*\|_{\F} \\
				&\quad +   \frac{\sqrt{d+m} + 1}{1-R_{2\br}}  \|(\scA^*(\bcE ))_{\max(2\br)}\|_{\F}.
			\end{split}
		\end{equation} Next, based on \eqref{ineq: RGD-iterative-bound}, we show inductively that for all $t \geq 0$, \eqref{ineq: RGD-error-bound} and $	\|\bcX^t - \bcX^*\|_{\F} \leq \frac{ R_{2\br} }{(d+m)(1+R_{2\br + \br^*}-R_{2\br})}\underline{\lambda}$ hold. First, it is clear the statements are true when $t = 0$. Suppose now \eqref{ineq: RGD-error-bound} and  $	\|\bcX^{t_0} - \bcX^*\|_{\F} \leq \frac{ R_{2\br} }{(d+m)(1+R_{2\br + \br^*}-R_{2\br})}\underline{\lambda}$ hold when $t = t_0$. Then
		\begin{equation} \label{ineq: RGD-induction}
			\begin{split}
				\|\bcX^{t_0+1} - \bcX^*\|_{\F} & \overset{ (a) } \leq   \frac{4 R_{2\br}(\sqrt{d+m} + 1)}{1-R_{2\br}}\|\bcX^{t_0} - \bcX^*\|_{\F} +   \frac{\sqrt{d+m} + 1}{1-R_{2\br}}  \|(\scA^*(\bcE ))_{\max(2\br)}\|_{\F} \\
				&  \overset{(b) } \leq \frac{1}{2} \|\bcX^{t_0} - \bcX^*\|_{\F} +   \frac{\sqrt{d+m} + 1}{1-R_{2\br}}  \|(\scA^*(\bcE ))_{\max(2\br)}\|_{\F}\\
				& \overset{(c)} \leq \frac{1}{2^{t_0 + 1}} \|\bcX^0 - \bcX^*\|_{\F} + \frac{2(\sqrt{d+m} + 1)}{1 - R_{2\br}} \|(\scA^*(\bcE ))_{\max(2\br)}\|_{\F}.
			\end{split}
		\end{equation} Here (a) is based on \eqref{ineq: RGD-iterative-bound} and the inductive assumption; (b) is because $ \frac{4 R_{2\br}(\sqrt{d+m} + 1)}{1-R_{2\br}} \leq 1/2$ based on the assumption $R_{2\br} \leq \frac{1}{8(\sqrt{d+m} + 1) + 1}$ and (c) is based on the inductive assumption. 
		
		Finally since  $ \underline{\lambda} \geq \frac{2(1 + R_{2\br + \br^*} - R_{2\br} ) (\sqrt{d+m} + 1)(d+m)}{R_{2\br} (1 - R_{2\br})} \|(\scA^*(\bcE ))_{\max(2\br)}\|_{\F} $, and $\|\bcX^{t_0} - \bcX^*\|_{\F} \leq  \frac{ R_{2\br} }{(d+m)(1+R_{2\br + \br^*}-R_{2\br})}\underline{\lambda}$, we have
		\begin{equation*}
			\|\bcX^{t_0+1} - \bcX^*\|_{\F} \overset{ \eqref{ineq: RGD-induction} }\leq \frac{1}{2} \|\bcX^{t_0} - \bcX^*\|_{\F} +   \frac{\sqrt{d+m} + 1}{1-R_{2\br}}  \|(\scA^*(\bcE ))_{\max(2\br)}\|_{\F} \leq \frac{ R_{2\br} }{(d+m)(1+R_{2\br + \br^*}-R_{2\br})}\underline{\lambda}.
		\end{equation*}
		This ends the induction and also finishes the proof of this theorem.
		\quad $\blacksquare$
		
		\subsection{Proof of Theorem \ref{th: local contraction general setting-RGN}}
		Similar to the proof of Theorem \ref{th: local contraction general setting-RGD}, we just need to prove the convergence result in the noisy case.
		
		First, the least squares in \eqref{eq: alg least square} can be viewed as an unconstrained least squares in the vector space $T_{\bcX^t} \bbM_\br$ and by Lemma \ref{lm: spectral norm bound of Atop A}, we know $\bcX^{t+0.5}$ can be compactly written as $\bcX^{t+0.5} = (P_{T_{\bcX^t}} \scA^* \scA P_{T_{\bcX^t}})^{-1} P_{T_{\bcX^t}} \scA^* (\bcY)$ and it is unique. So
		\begin{equation} \label{ineq: Xt+1 -Xstar}
			\begin{split}
				&\|\bcX^{t+1} - \bcX^*\|_{\F} \\
				=& \| \cH_{\br}(\bcX^{t+0.5}) - \bcX^* \|_{\F}\\
				\leq & \| \cH_{\br}(\bcX^{t+0.5}) - \bcX^{t+0.5}\|_{\F} + \|\bcX^{t+0.5} - \bcX^*\|_{\F} \\
				\overset{(a)}\leq & \sqrt{d+m} \|P_{\bbM_{\br}}(\bcX^{t+0.5}) - \bcX^{t+0.5}\|_{\F} + \|\bcX^{t+0.5} - \bcX^*\|_{\F}\\
				\overset{(b)}\leq & (\sqrt{d+m} + 1)\|\bcX^{t+0.5} - \bcX^*\|_{\F}\\
				=  & (\sqrt{d+m} + 1)\| (P_{T_{\bcX^t}} \scA^* \scA P_{T_{\bcX^t}})^{-1} P_{T_{\bcX^t}} \scA^* (\scA(\bcX^*) + \bcE) - \bcX^* \|_{\F}\\
				= & (\sqrt{d+m} + 1)\| (P_{T_{\bcX^t}} \scA^* \scA P_{T_{\bcX^t}})^{-1} P_{T_{\bcX^t}} \scA^* (\scA(P_{(T_{\bcX^t})\perp}\bcX^*) + \bcE)  - P_{(T_{\bcX^t})\perp}\bcX^* \|_{\F} \\
				\leq & (\sqrt{d+m} + 1) \left( \| (P_{T_{\bcX^t}} \scA^* \scA P_{T_{\bcX^t}})^{-1} P_{T_{\bcX^t}} \scA^* (\scA(P_{(T_{\bcX^t})\perp}\bcX^*) + \bcE) \|_\F + \| P_{(T_{\bcX^t})\perp}\bcX^* \|_{\F} \right) \\
				\overset{ \text{Lemma } \ref{lm: spectral norm bound of Atop A} }\leq & (\sqrt{d+m} + 1) \left( \frac{1}{1- R_{2\br}} \left( \| P_{T_{\bcX^t}} \scA^* \scA(P_{(T_{\bcX^t})\perp}\bcX^*) \|_\F +  \|(\scA^*(\bcE ))_{\max(2\br)}\|_{\F}\right) +  \| P_{(T_{\bcX^t})\perp}\bcX^* \|_{\F}  \right) \\
				\overset{ \text{Lemma } \ref{lm: orthogonal-projection-bound}, \ref{lm: orthogonal projection} } \leq  & (\sqrt{d+m} + 1) \left( \underbrace{\frac{2(1 + R_{2\br + \br^*} - R_{2\br})(d+m)}{(1-R_{2\br})\underline{\lambda}} \| \bcX^t - \bcX^* \|_\F^2}_{(A1)} +   \underbrace{\|(\scA^*(\bcE ))_{\max(2\br)}\|_{\F}/(1 - R_{2\br})}_{(A2)} \right),
			\end{split}
		\end{equation} here (a) is by the quasi-projection property of T-HOSVD and ST-HOSVD in Proposition \ref{prop: quasi-projection}; (b) is by the projection property of $P_{\bbM_{\br}}(\cdot)$.
		
		Denote $\Delta := \sqrt{\frac{ \underline{\lambda} \|(\scA^*(\bcE ))_{\max(2\br)}\|_{\F} }{  2(d+m)(1+R_{2\br + \br^*}-R_{2\br}) }}$.
		Notice, if $\|\bcX^t - \bcX^*\|_{\F} \geq \Delta$, (A1) dominates (A2) and if $\|\bcX^t - \bcX^*\|_{\F} \leq \Delta$, (A2) dominates (A1). 
		
		By \eqref{ineq: Xt+1 -Xstar}, when $\|\bcX^t - \bcX^*\|_{\F} \geq \Delta$ we have the error shrinks in each iteration as
		\begin{equation} \label{ineq: quadratic-shrinkage}
			\begin{split}
				\|\bcX^{t+1} - \bcX^*\|_{\F} \leq (\sqrt{d+m} + 1)\frac{4(1 + R_{2\br + \br^*} - R_{2\br})(d+m)}{(1-R_{2\br})\underline{\lambda}} \| \bcX^t - \bcX^* \|_\F^2.
			\end{split}
		\end{equation} We show by induction that when $\|\bcX^t - \bcX^*\|_{\F} \geq \Delta$, $\|\bcX^{t} - \bcX^* \|_{\F} \leq 2^{-2^{t}} \|\bcX^{0} - \bcX^*\|_{\F}$. This is true when $t = 0$. Suppose it holds when $t = t_0$, then
		\begin{equation*}
			\begin{split}
				\|\bcX^{t_0+1} - \bcX^*\|_{\F} \overset{ \eqref{ineq: quadratic-shrinkage}} \leq & (\sqrt{d+m} + 1)\frac{4(1 + R_{2\br + \br^*} - R_{2\br})(d+m)}{(1-R_{2\br})\underline{\lambda}} \| \bcX^{t_0} - \bcX^* \|_\F^2 \\
				\overset{(a)}\leq &  (\sqrt{d+m} + 1)\frac{4(1 + R_{2\br + \br^*} - R_{2\br})(d+m)}{(1-R_{2\br})\underline{\lambda}} \cdot 2^{-2^{t_0+1}} \|\bcX^{0} - \bcX^*\|^2_{\F}\\
				\overset{(b)} \leq & 2^{-2^{t_0+1}} \|\bcX^{0} - \bcX^*\|_{\F},
			\end{split}
		\end{equation*} where (a) is by the inductive assumption and (b) is due to the initialization condition. This finishes the induction.

		When $\|\bcX^t - \bcX^*\|_{\F} \leq \Delta$, the iteration error satisfies $ \|\bcX^{t+1} - \bcX^* \|_{\F} \leq \frac{2(\sqrt{d+m} + 1)}{1 - R_{2\br}} \|(\scA^*(\bcE ))_{\max(2\br)}\|_{\F}$. Combining two phases, we have
		\begin{equation*}
			\|\bcX^{t} - \bcX^* \|_{\F} \leq 2^{-2^{t}} \|\bcX^{0} - \bcX^*\|_{\F} + \frac{2(\sqrt{d+m} + 1)}{1 - R_{2\br}} \|(\scA^*(\bcE ))_{\max(2\br)}\|_{\F},\quad \forall t \geq 0.
		\end{equation*}
		\quad $\blacksquare$

		\subsection{Proof of Theorem \ref{th: MLE-guarantee} }
		First, simple computation yields
		\begin{equation} \label{ineq: likelihood-compare}
			\begin{split}
				0 & \geq \| \bcY - \scA(\widehat{\bcX}) \|_\F^2 - \| \bcY - \scA(\bcX^*) \|_\F^2\\
				& = 2 \langle \scA^*(\bcY - \scA(\bcX^*) ), \bcX^* - \widehat{\bcX}  \rangle + \langle \scA^* \scA(\widehat{\bcX} - \bcX^*), \widehat{\bcX} - \bcX^* \rangle \\
				&=2 \langle \scA^*(\bcY - \scA(\bcX^*) ), \bcX^* - \widehat{\bcX}  \rangle + \|\scA(\widehat{\bcX} - \bcX^*)\|_\F^2.
			\end{split}
		\end{equation} 
		
		Since $\scA$ satisfies $2\br$-TRIP, we have
		\begin{equation} \label{ineq: mle-guarantee-proof-ineq2}
			\begin{split}
				(1 - R_{2\br}) \| \widehat{\bcX} - \bcX^* \|_\F^2 \leq \|\scA(\widehat{\bcX} - \bcX^*)\|_\F^2 &\overset{\eqref{ineq: likelihood-compare} }\leq 2 \langle \scA^*(\bcY - \scA(\bcX^*) ), \widehat{\bcX} - \bcX^* \rangle \\
				&= 2 \langle \scA^*(\bcE), \widehat{\bcX} - \bcX^* \rangle\\
				& = 2 \|\widehat{\bcX} - \bcX^* \|_\F  \langle \scA^*(\bcE), \frac{\widehat{\bcX} - \bcX^*}{\|\widehat{\bcX} - \bcX^* \|_\F} \rangle \\
				& \overset{(a)} \leq 2 \|\widehat{\bcX} - \bcX^* \|_\F \|(\scA^*(\bcE ))_{\max(2\br)}\|_{\F},
			\end{split}
		\end{equation} here (a) is because $\frac{\widehat{\bcX} - \bcX^*}{\|\widehat{\bcX} - \bcX^* \|_\F}$ is of Tucker rank at most $2\br$ and Frobenius norm at most $1$ and by definition $\|(\scA^*(\bcE ))_{\max(2\br)}\|_\F$ is equal to $\sup_{\bcZ: \tuckerrank(\bcZ) \leq 2\br, \|\bcZ\|_\F \leq 1 } \langle \bcZ, \scA^*(\bcE ) \rangle $. The result follows from \eqref{ineq: mle-guarantee-proof-ineq2}. \quad $\blacksquare$

		\subsection{Proof of Theorem \ref{th: lower-bound} }
		We prove the two statements in Part 1 and Part 2 separately.
		
		{\bf Part 1}. The proof is based on a $\epsilon$-net argument. Let us first bound $\|(\scA^*(\bcE ))_{\max(\br)}\|_{\F}$. Recall
		\begin{equation} \label{eq: equi-xi-representation}
			\begin{split}
				\|(\scA^*(\bcE ))_{\max(\br)}\|_{\F} &:=  \sup_{ \substack{\U_k \in \bbO_{p_k,r_k} } } \|\scA^*(\bcE ) \times_{k=1}^{d +m} \U_k \U_k^\top \|_\F\\
				&=\sup_{ \substack{\bcS \in \bbR^{r_1 \times \cdots \times r_{d+m}}, \|\bcS\|_\F \leq 1 \\ \U_k \in \bbO_{p_k,r_k} } } \left \langle \scA^*(\bcE ), \bcS \times_{k=1}^{d+m} \U_k \right \rangle\\
				& = \sup_{ \substack{\bcS \in \bbR^{r_1 \times \cdots \times r_{d+m}}, \|\bcS\|_\F \leq 1 \\ \Q_k \in \bbR^{p_k\times r_k}, \|\Q_k\| \leq 1 } } \left \langle \scA^*(\bcE ), \bcS \times_{k=1}^{d+m} \Q_k \right \rangle
			\end{split}
		\end{equation}
		
		Since $\bcE_i$ has i.i.d. $N(0, \frac{\sigma^2}{n} )$ entries, for any fixed $\bcS$ with $\|\bcS\|_\F \leq 1$, and $\{\U_k \in \bbO_{p_k,r_k} \}_{k=1}^{d+m}$, condition on $\scA$, we have 
		\begin{equation*}
			\left \langle \scA^*(\bcE ), \bcS \times_{k=1}^{d+m} \U_k \right \rangle \Big| \scA = \left \langle \bcE, \scA(   \bcS \times_{k=1}^{d+m} \U_k) \right \rangle \Big| \scA \sim N(0, \frac{\tau^2 \sigma^2}{n} ), 
		\end{equation*} where $\tau^2 = \| \scA(   \bcS \times_{k=1}^{d+m} \U_k)\|_\F^2$. Thus, by the tail bound for the Gaussian random variable (see \cite[Chapter 2.1.2]{wainwright2019high}), we have
		\begin{equation}
			\begin{split}
				\bbP\left(  \left \langle \scA^*(\bcE ), \bcS \times_{k=1}^{d+m} \U_k \right \rangle  \geq t\Big| \scA \right) \leq \exp( - \frac{n t^2}{2 \tau^2 \sigma^2 } ).
			\end{split}
		\end{equation} When $n \geq C(\sum_{i=1}^d (p_i-r_i)r_i + \prod_{i=1}^d r_i )\log(d)$, by Proposition \ref{prop: TRIP-under-Gaussian-design}, $\scA$ satisfies TRIP with probability at least $1 - \exp(-c \underline{p})$ and denote the event that $\scA$ satisfies the TRIP property as $A$. So under $A$,  $\tau^2 \leq c \| \bcS \times_{k=1}^{d+m} \U_k \|_\F^2 \leq c$ for some $c > 0$ and
		\begin{equation}\label{ineq: single-point-Gaussian-upper}
			\bbP\left(  \left \langle \scA^*(\bcE ), \bcS \times_{k=1}^{d+m} \U_k \right \rangle  \geq t\Big| A \right) \leq \exp( - \frac{n t^2}{2 c \sigma^2 } ).
		\end{equation}

		By \cite[Lemma 7]{zhang2018tensor}, we can construct a $\epsilon$-net ($0<\epsilon < 1$) $\{ \bcS^{(1)},\cdots, \bcS^{(N_{\bcS})} \}$ for $\{ \bcS' \in \bbR^{r_1 \times \cdots \times r_{d+m}}: \|\bcS'\|_\F \leq 1 \}$ such that $$\sup_{\bcS': \|\bcS'\|_\F \leq 1 } \min_{i \leq N_{\bcS}} \|\bcS' - \bcS^{(i)}\|_\F \leq \epsilon $$
		with $N_{\bcS} \leq (\frac{3}{\epsilon})^{ \prod_{i=1}^{d+m} r_i } $.
		
		At the same time, by \cite[Proposition 8]{szarek1982nets}, for each $k = 1,\ldots, d+m$, we can construct a $\epsilon$-net $\{ \U_k^{(1)},\ldots, \U_k^{(N_k)} \}$ on the Grassmann manifold of $r_k$-dimensional subspaces in $\bbR^{p_k}$ with the metric $d(\U_1, \U_2) = \| \U_1 \U_1^\top - \U_2 \U_2^\top \|$ such that 
		\begin{equation*}
			\sup_{\U_k \in \bbO_{p_k, r_k}} \min_{i \leq N_k} d( \U_k, \U_k^{(i)} ) \leq \epsilon
		\end{equation*} with $N_k \leq (\frac{c_0}{\epsilon})^{r_k(p_k - r_k)} $ for some absolute constant $c_0 > 0$.
		
		Suppose $$(\bcS^*, \{ \U_k^* \}_{k=1}^{d+m} ) = \argmax_{ \substack{\bcS \in \bbR^{r_1 \times \cdots \times r_{d+m}}, \|\bcS\|_\F \leq 1 \\ \U_k \in \bbO_{p_k,r_k} } } \left \langle \scA^*(\bcE ), \bcS \times_{k=1}^{d+m} \U_k \right \rangle$$
		and denote $T := \left \langle \scA^*(\bcE ), \bcS^* \times_{k=1}^{d+m} \U^*_k \right \rangle = \|(\scA^*(\bcE ))_{\max(\br)}\|_{\F}$.
		For $k \in [d+m]$, we can find $\U_k^{(i_k)} $ in the corresponding $\epsilon$-net such that $d(\U_k^{(i_k)}, \U_k^*) \leq \epsilon$. Let $\O_k = \argmin_{\O \in \bbO_{r_k}} \| \U_k^* \O - \U_k^{(i_k)} \|$. By \cite[Lemma 1]{cai2018rate}, we have $\| \U_k^* \O_k - \U_k^{(i_k)} \| \leq \sqrt{2} d(\U_k^{(i_k)}, \U_k^*) \leq \sqrt{2} \epsilon$. Denote $\widebar{\bcS} = \bcS^* \times_{k=1}^{d+m} \O_k^\top$ and let $\bcS^{(i_0)}$ be the one in the core tensor $\epsilon$-net such that $\|\bcS^{(i_0)} - \widebar{\bcS}\|_\F \leq \epsilon$. Thus
		\begin{equation} \label{ineq: telescoping}
			\begin{split}
				&\left \langle \scA^*(\bcE ), \widebar{\bcS} \times_{k=1}^{d+m} \U^*_k \O_k \right \rangle - \left \langle \scA^*(\bcE ), \bcS^{(i_0)} \times_{k=1}^{d+m} \U^{(i_k)}_k \right \rangle \\
				=&  \left \langle \scA^*(\bcE ), (\widebar{\bcS} - \bcS^{(i_0)})  \times_{k=1}^{d+m} \U^*_k \O_k \right \rangle \\
				& + \sum_{k=1}^{d+m} \left \langle \scA^*(\bcE ), \bcS^{(i_0)} \times_{j < k} \U_j^{(i_j)} \times_k (\U^*_k \O_k - \U^{(i_k)}_k ) \times_{j > k}  \U^*_j \O_j \right \rangle\\
				\leq & (\sqrt{2}(d+m) + 1)\epsilon T.
			\end{split}
		\end{equation}
		So for $\epsilon < 1/(2 (\sqrt{2}(d+m) + 1) ) $ we have
		\begin{equation*}
			\begin{split}
				T &= \left \langle \scA^*(\bcE ), \bcS^{(i_0)} \times_{k=1}^{d+m} \U^{(i_k)}_k \right \rangle + T - \left \langle \scA^*(\bcE ), \bcS^{(i_0)} \times_{k=1}^{d+m} \U^{(i_k)}_k \right \rangle\\
				& = \left \langle \scA^*(\bcE ), \bcS^{(i_0)} \times_{k=1}^{d+m} \U^{(i_k)}_k \right \rangle  + \left \langle \scA^*(\bcE ), \widebar{\bcS} \times_{k=1}^{d+m} \U^*_k \O_k \right \rangle - \left \langle \scA^*(\bcE ), \bcS^{(i_0)} \times_{k=1}^{d+m} \U^{(i_k)}_k \right \rangle\\
				& \overset{ \eqref{ineq: telescoping} }\leq  \left \langle \scA^*(\bcE ), \bcS^{(i_0)} \times_{k=1}^{d+m} \U^{(i_k)}_k \right \rangle + T/2.
			\end{split}
		\end{equation*} This implies $T \leq 2 \left \langle \scA^*(\bcE ), \bcS^{(i_0)} \times_{k=1}^{d+m} \U^{(i_k)}_k \right \rangle$.
		
		Then by union bound, we have
		\begin{equation*}
			\begin{split}
				\bbP\left( T \geq t | A  \right) &\leq \bbP\left( \max_{i'_0,\ldots,i'_{d+m}}  \left \langle \scA^*(\bcE ), \bcS^{(i'_0)} \times_{k=1}^{d+m} \U^{(i'_k)}_k \right \rangle \geq t/2 \Big| A \right)\\
				& \overset{ \eqref{ineq: single-point-Gaussian-upper} } \leq  \left(N_{\bcS} \prod_{i=1}^{d+m} N_i \right) \exp( - \frac{n t^2}{2c \sigma^2} ) = \exp( C \cdot df - nt^2/(2c\sigma^2)  ),
			\end{split}
		\end{equation*} where $df = \sum_{i=1}^{d+m} r_i(p_i - r_i) + \prod_{i=1}^{d+m} r_i$. Notice here $C$ will depend on $d,m$ and for simplicity we omit them here. So we have $\bbP(T \geq c \sqrt{\frac{df}{n}} \sigma  | A ) \leq \exp( - c' \cdot df )$. Overall,
		\begin{equation*}
			\bbP(T \leq  c  \sqrt{\frac{df}{n}} \sigma ) \geq \bbP(A) \bbP(T \leq t | A  ) \geq (1- \exp( - c' \cdot df ))(1 - \exp(-c_1 \underline{p} )) \geq 1 - \exp(-c \underline{p} ).
		\end{equation*} Now to bound $\|(\scA^*(\bcE ))_{\max(2\br)}\|_{\F}$, we just need to replace $\br$ by $2\br$ in the bound for $T$ and up to a constant (depending on $d$ and $m$ only), the same upper bound for $T$ holds for $\|(\scA^*(\bcE ))_{\max(2\br)}\|_{\F}$. 
		
		Combining the upper bound of $\|(\scA^*(\bcE ))_{\max(2\br)}\|_{\F}$ with Theorem \ref{th: MLE-guarantee}, we have given event $A$, $\|\widehat{\bcX} - \bcX^*\|_\F \geq c t$ happens with probability at most $\exp(-nt^2/(2c\sigma^2))$ for $t \geq c' \sigma \sqrt{\frac{df}{n}} $. Thus
		\begin{equation*}
			\begin{split}
				& \bbE(\|\widehat{\bcX} - \bcX^*\|_\F | A ) \leq \sqrt{ \bbE (\|\widehat{\bcX} - \bcX^*\|^2_\F|A) } \\
				= & \sqrt{ \int_{0}^{\infty} \bbP(\|\widehat{\bcX} - \bcX^*\|^2_\F \geq t |A ) }dt\\
				= & \sqrt{ \int_{0}^{a} \bbP(\|\widehat{\bcX} - \bcX^*\|^2_\F \geq t |A )dt + \int_{a}^{\infty} \bbP(\|\widehat{\bcX} - \bcX^*\|^2_\F \geq t |A )dt }\\
				\leq & \sqrt{ a +  \int_{a}^{\infty} \bbP(\|\widehat{\bcX} - \bcX^*\|_\F \geq \sqrt{t} |A )dt }\\
				\overset{\text{let } a = c^{'2} \sigma^2 \frac{df}{n} }\leq & \sqrt{ a + \int^{\infty}_{a} \exp(-nt/(2c\sigma^2))dt } \\
				\leq & C \sigma \sqrt{\frac{df}{n}}.
			\end{split}
		\end{equation*} So $\bbE(\|\widehat{\bcX} - \bcX^*\|_\F ) \leq C \sigma \sqrt{\frac{df}{n}}$. This finishes the proof for part 1.
		
		{\bf Part 2.} Suppose we can find a set of $\{\bcX^i \}_{i=1}^N \in \mathcal{F}_{\bp,\br}$ such that $\min_{i\neq j} \| \bcX^i - \bcX^j \|_\F \geq s$, by the standard argument of reducing the problem of providing a minimax risk lower bound to lower bounding the probability of error in a multiple hypothesis testing problem \citep[Chapter 2]{tsybakov2009introduction}, we have
		\begin{equation} \label{ineq: general-minimax-lower-bd-argu}
			\inf_{\widehat{\bcX}} \sup_{ \bcX \in \mathcal{F}_{\bp,\br} }\bbE \| \widehat{\bcX} - \bcX \|_\F \geq \frac{s}{2} \inf_{\widehat{\bcX} } \frac{1}{N} \sum_{i=1}^N \bbP_{\bcX^i} ( \widehat{\bcX} \neq \bcX^i ).   
		\end{equation} 
		
		We will use Fano's Lemma to lower bound the right-hand side of \eqref{ineq: general-minimax-lower-bd-argu}. Before that, let us first compute the Kullback-Leibler (KL) divergence between two different distributions in our setting. For $i = 1,\ldots, N$, let $\bbQ^i$ denotes the conditional distribution of $\bcY$ given $\scA$ and $\bcX^* = \bcX^i$. Suppose $\bcY \sim \bbQ^i$, for $j_1 \in [p_{d+1}],\ldots,j_m \in [p_{d+m}]$, let $\bcY^{(j_1,\ldots,j_m)} := ( \bcY_{[1,j_1,\ldots,j_m]}, \ldots, \bcY_{[n,j_1,\ldots,j_m]} )^\top \in \bbR^n$ and $\bcE^{(j_1,\ldots,j_m)}$ be defined in the same way.  In addition, let $\bcX^{i(j_1,\ldots,j_m)} := \bcX^i_{[:,\ldots,:,j_1,\ldots,j_m]} \in \bbR^{p_1 \times \cdots \times p_d}$ and define $\scA': \bbR^{p_1 \times \cdots \times p_d} \to \bbR^n$ such that $(\scA'(\bcX^{i(j_1,\ldots,j_m)}))_{[k]} = \langle \bcX^{i(j_1,\ldots,j_m)}, \bcA_k \rangle $ where $\{ \bcA_k \}_{k=1}^n$ are the same set of tensor covariates given in $\scA$. Thus, we have $\bcY^{(j_1,\ldots,j_m)} = \scA'(\bcX^{i(j_1,\ldots,j_m)}) + \bcE^{(j_1,\ldots,j_m)}$ and denote $\bbQ^{i(j_1,\ldots,j_m)}$ as the conditional distribution of $\bcY^{(j_1,\ldots,j_m)}$ condition on $\scA$ and $\bcX^* = \bcX^i$. Condition on $\scA$, $\bcY^{(j_1,\ldots,j_m)}$ are independent for different $(j_1,\ldots,j_m)$s. Due to the Gaussian ensemble design we assume, i.e., $\bcE_i$ has i.i.d. $N(0, \sigma^2/n )$ entries, and the linear regression model we have on $\bcY^{(j_1,\ldots,j_m)}$, the KL divergence between $\bbQ^{i_1}$ and $\bbQ^{i_2}$, denoted by $\KL(\bbQ^{i_1}|| \bbQ^{i_2} )$, can be computed as follows 
		\begin{equation} \label{eq: KL-divergence}
			\begin{split}
				\KL(\bbQ^{i_1}|| \bbQ^{i_2} ) &= \sum_{j_1 \in [p_{d+1}]} \cdots \sum_{j_m \in [p_{d+m}] } \KL(\bbQ^{i_1(j_1,\ldots,j_m)} || \bbQ^{i_2(j_1,\ldots,j_m)} ) \\
				&= \sum_{j_1,\ldots,j_m} \sum_{i=1}^n \frac{n}{2\sigma^2} \left \langle \bcX^{i_1(j_1,\ldots,j_m)} - \bcX^{i_2 (j_1,\ldots,j_m)}, \bcA_i  \right \rangle^2.
			\end{split}
		\end{equation}
		Then by the Fano's Lemma, e.g., see \citep[Chapter 15.3.2]{wainwright2019high}, we have
		\begin{equation} \label{ineq: general-lower-bd-continue}
			\begin{split}
				& ~~ \frac{s}{2} \inf_{\widehat{\bcX} } \frac{1}{N} \sum_{i=1}^N \bbP_{\bcX^i} ( \widehat{\bcX} \neq \bcX^i ) = \frac{s}{2} \inf_{\widehat{\bcX} } \frac{1}{N} \sum_{i=1}^N \bbE_{\scA} (\bbP_{\bcX^i} ( \widehat{\bcX} \neq \bcX^i | \scA )) \\
				& \geq \frac{s}{2}\bbE_{\scA} \inf_{\widehat{\bcX} } \frac{1}{N} \sum_{i=1}^N (\bbP_{\bcX^i} ( \widehat{\bcX} \neq \bcX^i | \scA )) \\
				& \geq \frac{s}{2} \bbE_\scA \left( 1 -  \frac{ \max_{i_1 \neq i_2} \KL( \bbQ^{i_1} || \bbQ^{i_2}) + \log 2   }{\log N}   \right)\\
				& \overset{ (a) } = \frac{s}{2} \left(  1 - \frac{ \frac{1}{2\sigma^2} \max_{i_1 \neq i_2} \sum_{i=1}^n  \sum_{j_1,\ldots,j_m}\left \| \bcX^{i_1(j_1,\ldots,j_m)} - \bcX^{i_2 (j_1,\ldots,j_m)}\right \|_\F^2 + \log 2  }{\log N}  \right)\\
				& =  \frac{s}{2} \left(  1 - \frac{ \frac{n}{2\sigma^2} \max_{i_1 \neq i_2}  \left \| \bcX^{i_1} - \bcX^{i_2 }\right \|_\F^2 + \log 2  }{\log N}  \right),
			\end{split}
		\end{equation} here (a) is because of \eqref{eq: KL-divergence} and $\bcA_i$ has i.i.d. $N(0,1/n)$ entries.
		
		Next, we consider two constructions for the sets of $\{ \bcX^i \}_{i=1}^N$ so that we can have a proper lower bound for $\min_{i_1\neq i_2} \| \bcX^{i_1} - \bcX^{i_2} \|_\F $ and a proper upper bound for $\max_{i_1 \neq i_2}  \left \| \bcX^{i_1} - \bcX^{i_2 }\right \|_\F$.
		
		{\noindent \bf Construction 1}. For $k = 1,\ldots, d+m$, pick $\U_k \in \bbO_{p_k, r_k}$. Given any $\delta > 0$ and $\min_i r_i \geq C'$, by a slight modified version of \cite[Lemma 5]{agarwal2012noisy}, we can construct a set of $r_1 \times \cdots \times r_{d+m}$ full Tucker rank tensors $\{ \bcS^{(1)}, \ldots, \bcS^{(N_0)} \}$ with cardinality $N_0 \geq \frac{1}{4} \exp( \prod_{i=1}^{d+m} r_i/128 )  $ such that: (1) $\|\bcS^{(i)}\|_\F = \delta $ holds for all $i = 1,\ldots,N_0$, (2) $\| \bcS^{(i_1)} - \bcS^{(i_2)} \|_\F \geq \delta$ for all $i_1, i_2 \in [N_0], i_1 \neq i_2$. Let $\bcX^i := \bcS^{(i)} \times_{k=1}^{d+m} \U_k $ for $i \in [N_0]$ and clearly $\bcX^i \in \mathcal{F}_{\bp,\br}$.
		
		For this set, we have 
		\begin{equation} \label{ineq: lower-bound-construction-1}
			\begin{split}
				\max_{i_1 \neq i_2}  \left \| \bcX^{i_1} - \bcX^{i_2 }\right \|_\F &\leq \max_{i_1 \neq i_2} ( \left \| \bcX^{i_1}\|_\F + \| \bcX^{i_2 }\right \|_\F ) = 2 \delta,\\
				\min_{i_1\neq i_2} \| \bcX^{i_1} - \bcX^{i_2} \|_\F &=  \min_{i_1\neq i_2} \|(\bcS^{(i_1)} - \bcS^{(i_2)} ) \times_{k=1}^{d+m} \U_k \|_\F \geq \delta.
			\end{split}
		\end{equation}
		
		Plug in the results in \eqref{ineq: lower-bound-construction-1} into \eqref{ineq: general-minimax-lower-bd-argu} and \eqref{ineq: general-lower-bd-continue}, we have 
		\begin{equation} \label{ineq: minimax-lower-bound-part-1}
			\begin{split}
				\inf_{\widehat{\bcX}} \sup_{ \bcX \in \mathcal{F}_{\bp,\br} }\bbE \| \widehat{\bcX} - \bcX \|_\F \geq \frac{\delta}{2} \left(  1 - \frac{ 2n \delta^2/\sigma^2 + \log 2  }{c_1 \prod_{i=1}^{d+m} r_i }   \right)\overset{(a)}\geq c \sqrt{\frac{\prod_{i=1}^{d+m} r_i}{n}} \sigma, 
			\end{split}
		\end{equation} where in (a) we pick $\delta^2 $ to be $c'\frac{\prod_{i=1}^{d+m} r_i}{n} \sigma^2$ for some $c' > 0$.
		
		{\noindent \bf Construction 2}. Let $\bcS \in \bbR^{r_1 \times \cdots \times r_{d+m}}$ be a fixed core tensor of full Tucker rank such that 
		\begin{equation*}
			c \lambda \leq \min_{i} \sigma_{r_i}( \cM_i(\bcS) ) \leq \max_i \sigma_1( \cM_i(\bcS) ) \leq C \lambda, 
		\end{equation*} for some $c, C, \lambda >0$. For $k = 1,\ldots, d+m$, pick $\U_k \in \bbO_{p_k, r_k}$. 
		
		Let us first focus on mode-$1$. Consider the metric space $(\mathcal{G}_{p_1,r_1}, d(\cdot, \cdot) )$, where $ \mathcal{G}_{p_1,r_1} $ denotes the Grassmann manifold of $r_1$ dimensional subspaces in $\bbR^{p_1}$ and $d(\U_1, \U_2) = \|\U_1\U_1^\top - \U_2 \U_2^\top \|_\F$. Consider the following ball of radius $\epsilon > 0$ and centered at $\U_1$:
		\begin{equation*}
			B(\U_1, \epsilon ) = \{ \U' \in \bbO_{p_1,r_1}: d(\U',\U_1) \leq \epsilon  \}.
		\end{equation*} By \cite[Lemma 1]{cai2013sparse}, for $0 < \alpha < 1$ and $0 < \epsilon \leq 1$, there exists $\U_1^{'(1)}, \ldots, \U_{1}^{'(N_1)} \in B(\U_1, \epsilon )$ such that $\min_{1 \leq i < j \leq N_1} d(\U_1^{'(i)}, \U_1^{'(j)} ) \geq \alpha \epsilon $ with $N_1 \geq (c_0/ \alpha )^{r_1(p_1 - r_1)}  $ for some absolute constant $c_0$. 
		
		Let $\O_i = \argmin_{\O \in \bbO_{r_1} } \| \U_1^{'(i)} \O - \U_1 \|_\F$ and let $\U_1^{(i)} = \U_1^{'(i)} \O_i$ for $i = 1,\ldots, N_1$. By \cite[Lemma 1]{cai2018rate}, we have $\| \U_1^{(i)} - \U_1 \|_\F \leq d(\U_1^{'(i)}, \U_1) \leq \epsilon$ for all $i \in [N_1]$. We then construct $\bcX^i = \bcS \times_1 \U_1^{(i)} \times_{k = 2}^{d+m} \U_k $ for $i \in [N_1]$. For this set, we have 
		\begin{equation} \label{ineq: lower-bound-construction-2}
			\begin{split}
				\max_{i_1 \neq i_2}  \left \| \bcX^{i_1} - \bcX^{i_2 }\right \|_\F &\leq \max_{i_1 \neq i_2} \| \bcS \times_1 (\U_1^{(i_1)} - \U_1^{(i_2)} ) \times_{k = 2}^{d+m} \U_k \|_\F \\
				&  \leq \max_{i_1 \neq i_2} (\| \bcS \times_1 (\U_1^{(i_2)} - \U_1)  \times_{k = 2}^{d+m} \U_k \|_\F + \| \bcS \times_1 (\U_1^{(i_1)} - \U_1)  \times_{k = 2}^{d+m} \U_k \|_\F ) \\
				&\leq  2C\epsilon \lambda  ,\\
				\min_{i_1\neq i_2} \| \bcX^{i_1} - \bcX^{i_2} \|_\F &=  \min_{i_1\neq i_2} \| \bcS \times_1 (\U_1^{(i_1)} - \U_1^{(i_2)} ) \times_{k = 2}^{d+m} \U_k \|_\F\\
				&  \overset{(a)}\geq \min_{i_1\neq i_2} \min_{\O \in \bbO_{r_1}} \| \bcS \times_1 (\U_1^{(i_1)} - \U_1^{(i_2)} \O ) \times_{k = 2}^{d+m} \U_k \|_\F \geq c \alpha \epsilon \lambda/\sqrt{2}.
			\end{split}
		\end{equation} Here (a) is because $\min_{\O \in \bbO_{r_1}} \|\U_1^{(i_1)} - \U_1^{(i_2)} \O \|_\F \geq d(\U_1^{(i_1)}, \U_2^{(i_2)})/\sqrt{2}$ by \cite[Lemma 1]{cai2018rate} and $ d(\U_1^{(i_1)}, \U_2^{(i_2)}) =  d(\U_1^{'(i_1)}, \U_2^{'(i_2)}) \geq \alpha \epsilon$.

		Plug \eqref{ineq: lower-bound-construction-2} into \eqref{ineq: general-minimax-lower-bd-argu} and \eqref{ineq: general-lower-bd-continue}, we have 
		\begin{equation} \label{ineq: minimax-lower-bound-part-2}
			\begin{split}
				\inf_{\widehat{\bcX}} \sup_{ \bcX \in \mathcal{F}_{\bp,\br} }\bbE \| \widehat{\bcX} - \bcX \|_\F \geq \frac{c \alpha \epsilon \lambda }{2\sqrt{2}} \left(  1 - \frac{ Cn \epsilon^2 \lambda^2 /\sigma^2 + \log 2  }{\log(c_0/\alpha) r_1(p_1 - r_1) }   \right)\overset{(a)}\geq c_1 \sqrt{\frac{r_1(p_1 - r_1)}{n}} \sigma, 
			\end{split}
		\end{equation} where in (a) we pick $\epsilon^2 $ to be $c'\frac{r_1(p_1 - r_1)}{n \lambda^2 } \sigma^2$, $\alpha = (c_0 \wedge 1 )/8 $. 
		
		We can apply similar procedures on modes $2$ to $(d+m)$ and get 
		\begin{equation} \label{ineq: minimax-lower-bound-part-3}
			\begin{split}
				\inf_{\widehat{\bcX}} \sup_{ \bcX \in \mathcal{F}_{\bp,\br} }\bbE \| \widehat{\bcX} - \bcX \|_\F \geq c \sqrt{\frac{\sum_{k=1}^{d+m}r_k(p_k - r_k)}{n}} \sigma. 
			\end{split}
		\end{equation} Combining it with \eqref{ineq: minimax-lower-bound-part-1}, we obtain the full minimax risk lower bound. \quad $\blacksquare$

		\section{Proofs in Section \ref{sec: applications-initializations}}
		
		\subsection{Proof of Theorem \ref{th: initialization-scalar-on-tensor} }
		Denote $\widetilde{\sigma}^2 =  \|\bcA\|_{\F}^2 + \sigma^2$. We first check the condition required for applying Theorem \ref{th: perturbation-HOOI}. Recall $\widetilde{\U}_k^0 = \SVD_{r_k}( \cM_k( \scA^*(\y) ) )$ and let $\widetilde{\U}_k^0 = [\widebar{\U}_k^0 \quad \widecheck{\U}_k^0]$ where $\widebar{\U}_k^0$ is composed of the first $r_k^*$ columns of $\widetilde{\U}_k^0$ and $\widecheck{\U}_k^0$ contains the remaining $(r_k - r_k^*)$ columns of $\widetilde{\U}_k^0$. By the proof of \cite[Theorem 4, step 1]{zhang2020islet}, we have with probability at least $1 - \underline{p}^{-C}$ for some $C > 0$ such that when  $n \geq c(d) (\|\bcX^*\|_{\F}^2 + \sigma^2)  \frac{ (\prod_{i=1}^d p_i  )^{1/2} + \bar{p} }{\underline{\lambda}^2}$, we have 
		\begin{equation} \label{ineq: tensor regression two fact}
			\begin{split}
				&\left\|\sin\Theta(\widebar{\U}_k^0,\U_k)\right\|= \|\widebar{\U}_{k\perp}^{0\top} \U_k\| \leq \frac{\sqrt{p_k/n}\tilde{\sigma} \underline{\lambda} + (\prod_{k=1}^d p_k)^{1/2} \widetilde{\sigma}^2 /n}{\underline{\lambda}^2} \leq \frac{\sqrt{2}}{2}. 
			\end{split}
		\end{equation}
		Moreover, 
		\begin{equation} \label{eq: overparameterization-deal-trick}
			\begin{split}
				\|\widetilde{\U}_{k\perp}^{0\top} \U_k\| = \|(\I - P_{\widetilde{\U}_{k}^{0}} )\U_k\| = \|(\I - P_{\widecheck{\U}_{k}^{0}})(\I - P_{\widebar{\U}_{k}^{0}}) \U_k\| \leq \|(\I - P_{\widebar{\U}_{k}^{0}}) \U_k\| =  \|\widebar{\U}_{k\perp}^{0\top} \U_k\| \leq \frac{\sqrt{2}}{2}.
			\end{split}
		\end{equation}

		Consider applying Theorem \ref{th: perturbation-HOOI} in our setting, we have $\bcT = \bcX^*$, $\widetilde{\bcT} = \scA^*(\y)$ and $\bcZ = \scA^*(\y) - \bcX^* = \scA^*\scA(\bcX^*) - \bcX^* + \scA^*(\bvarepsilon )$. With probability at least $1 - \exp(-c\underline{p} )$, we have 
		\begin{equation*}
			\begin{split}
				\|\bcZ_{\max(\br)}\|_{\F} &= \sup_{\U'_k \in \bbO_{p_k, r_k}, k=1,\ldots,d} \|\bcZ \times_{k=1}^d P_{\U'_k} \|_{\F} \\
				&= \sup_{\U'_k \in \bbO_{p_k, r_k}, \|\bcW\|_\F \leq 1} \langle \bcZ \times_{k=1}^d P_{\U'_k}, \bcW \rangle \\
				&= \sup_{\U'_k \in \bbO_{p_k, r_k}, \|\bcW\|_\F \leq 1} \langle \bcZ , \bcW \times_{k=1}^d P_{\U'_k} \rangle\\
				&\leq \sup_{\U'_k \in \bbO_{p_k, r_k}, \|\bcW\|_\F \leq 1} \langle \scA^*\scA (\bcX^*) - \bcX^*  , \bcW \times_{k=1}^d P_{\U'_k} \rangle \\
				& \quad + \sup_{\U'_k \in \bbO_{p_k, r_k}, \|\bcW\|_\F \leq 1} \langle \scA^*(\bvarepsilon )  , \bcW \times_{k=1}^d P_{\U'_k} \rangle \\
				& \overset{ \text{Lemma } \ref{lm:retricted orthogonal property}, \text{Theorem } \ref{th: lower-bound} } \leq  R_{2\br} \|\bcX^*\|_{\F} + C \sigma\sqrt{\frac{\sum_{i=1}^d r_i (p_i-r_i) + \prod_{i=1}^d r_i }{n}}\\
				& \overset{ \text{Proposition }\ref{prop: TRIP-under-Gaussian-design} } \leq C( \|\bcX^*\|_{\F} + \sigma ) \sqrt{\frac{\sum_{i=1}^d r_i (p_i-r_i) + \prod_{i=1}^d r_i }{n}}.
			\end{split}
		\end{equation*} 
		
		So when $n$ satisfies the condition indicated in the proposition, we have $\|\bcZ_{\max(\br)}\|_{\F} \leq c \underline{\lambda}$. Thus by Theorem \ref{th: perturbation-HOOI}, we have the output of the algorithm satisfies
		\begin{equation*}
			\|\bcX^0 - \bcX^*\|_{\F} \leq c'(d) \|\bcZ_{\max(\br)}\|_{\F} \leq c \underline{\lambda} 
		\end{equation*} as long as $n \geq c(d)  \left( \frac{(\|\bcX^*\|_{\F}^2 + \sigma^2)}{\underline{\lambda}^2} \left((\prod_{i=1}^d p_i  )^{1/2} +  \sum_{i=1}^d (p_i-r_i) r_i + \prod_{i=1}^d r_i \right) \right)$ for sufficiently large $c(d)$.

		Since $ \|(\scA^*(\bvarepsilon ))_{\max(2\br)}\|_{\F} \leq c \sigma \sqrt{df/n} $ by Theorem \ref{th: lower-bound}, when $ \underline{\lambda} \leq c_1 \sigma \sqrt{df/n}   $, the estimator is already optimal after initialization, no further refinement is needed. If $ \underline{\lambda} \geq c_1 \sigma \sqrt{df/n}$, TRIP holds for $n$ indicated in the proposition, then by Theorem \ref{th: local contraction general setting-RGD}, after $t_{\max} \geq \log\left( \frac{\underline{\lambda} \sqrt{n/df} }{ c_1(d)\sigma }   \right)$, we have $ \| \bcX^{t_{\max} } - \bcX^* \|_\F \leq c_3(d) \sigma \sqrt{\frac{df}{n}}$. A similar argument applies to the study of RGN.
		\quad $\blacksquare$

		\subsection{Proof of Theorem \ref{th: initialization-tensor-on-vector}}
		Since $\A$ is a $n$-by-$p_1$ matrix with i.i.d. $N(0,1/n)$ entries and $n \geq Cp_1$, we have the $\R_\A$ factor of the QR decomposition of $\A$ has the spectrum bound
		\begin{equation} \label{ineq: spectral-bound-A}
			1-c \leq  \sigma_{\min}(\R_\A) \leq \sigma_{\max}(\R_\A) \leq 1+c  
		\end{equation}
		for some $1 > c > 0$ with probability at least $1 - \exp(-c_1 p_1)$ by the standard random matrix theory, see \cite[Corollary 5.35]{vershynin2010introduction}.
		
		Let $\widebar{\bcY}:=\bcY \times_1 \Q_\A^\top = \bcX^* \times_1 \R_\A + \widebar{\bcE}$ where $\widebar{\bcE} = \bcE \times_1 \Q_\A^\top$ has i.i.d. $N(0,\sigma^2/n)$ entries. Denote $\widebar{\bcX}^* := \bcX^* \times_1 \R_\A$, and its Tucker decomposition as $\widebar{\bcS} \times_{i=1}^{m+1} \widebar{\U}_i $, we have $ (1+c) \underline{\lambda} \geq \underline{\lambda}':= \min_{i} \sigma_{r_i^*}( \cM_i( \widebar{\bcX}^* ) ) \geq (1-c) \underline{\lambda} $. Let $\widetilde{\U}_k^0 = [\widebar{\U}_k^0 \quad \widecheck{\U}_k^0]$ where $\widebar{\U}_k^0$ is composed of the first $r_k^*$ columns of $\widetilde{\U}_k^0$ and $\widecheck{\U}_k^0$ contains the remaining $(r_k - r_k^*)$ columns of $\widetilde{\U}_k^0$. By the proof of Theorem 1 in \cite{zhang2018tensor}, we have for the sample size indicated in the proposition, with probability at least $1-\exp(-c\underline{p})$, the following inequalities hold for $k =1,\ldots,m+1$:
		\begin{equation} \label{ineq: tensor SVD two fact}
			\begin{split}
				&\left\|\sin\Theta(\widebar{\U}_k^0,\U_k)\right\| = \|\widebar{\U}_{k\perp}^{0\top} \U_k\| \leq \frac{\sigma\sqrt{p_k}}{\underline{\lambda}' \sqrt{n} } + \frac{(\prod_{k=1}^{m+1} p_k)^{1/2} \sigma^2 }{\underline{\lambda}^{'2}n } \leq \frac{\sqrt{2}}{2}.
			\end{split}
		\end{equation} 
		Moreover, 
		\begin{equation*}
			\begin{split}
				\|\widetilde{\U}_{k\perp}^{0\top} \U_k\| = \|(\I - P_{\widetilde{\U}_{k}^{0}} )\U_k\| = \|(\I - P_{\widecheck{\U}_{k}^{0}})(\I - P_{\widebar{\U}_{k}^{0}}) \U_k\| \leq \|(\I - P_{\widebar{\U}_{k}^{0}}) \U_k\| =  \|\widebar{\U}_{k\perp}^{0\top} \U_k\| \leq \frac{\sqrt{2}}{2}.
			\end{split}
		\end{equation*}

		Consider applying Theorem \ref{th: perturbation-HOOI} in our setting, we have $\bcT = \widebar{\bcX}^*$, $\widetilde{\bcT} = \widebar{\bcY}$ and $\bcZ = \widebar{\bcY} - \widebar{\bcX}^* = \widebar{\bcE}$. Moreover, by a similar proof of Theorem \ref{th: lower-bound} part I, we have $\|\bcZ_{\max(\br)}\|_{\F} = \| \widebar{\bcE}_{\max(\br)} \|_\F \leq C \sigma \sqrt{ \frac{df}{n}}$ for some $C > 0$ with probability at least $1 - \exp(-c \underline{p} ).$ So when $n$ satisfies the condition indicated in the proposition, we have $\|\bcZ_{\max(\br)}\|_{\F} \leq c \underline{\lambda}'$. Thus by Theorem \ref{th: perturbation-HOOI}, we have for the output of OHOOI, $\widebar{\bcX}^0$, satisfies
		\begin{equation*}
			\|\widebar{\bcX}^0 - \widebar{\bcX}^*\|_{\F} \leq c'(d) \|\bcZ_{\max(\br)}\|_{\F} \leq c' \underline{\lambda}' \leq C' \underline{\lambda}. 
		\end{equation*}
		So $\| \bcX^0 - \bcX^* \|_\F = \|(\widebar{\bcX}^0 - \widebar{\bcX}^*) \times_1 \R_\A^{-1} \|_\F \leq c \|\widebar{\bcX}^0 - \widebar{\bcX}^*\|_{\F} \leq C' \underline{\lambda}$.
		
		Since $ \|(\scA^*(\bcE ))_{\max(2\br)}\|_\F = \|(\bcE \times_1 \A^\top)_{\max(2\br)} \|_\F =  \|(\widebar{\bcE} \times_1 \R_\A^\top)_{\max(2\br)}\|_{\F} \leq c \sigma \sqrt{df/n} $, when $ \underline{\lambda} \leq c_1 \sigma \sqrt{df/n}$, the initialization is already optimal, no further refinement is needed. If $ \underline{\lambda} \geq c_1 \sigma \sqrt{df/n}$, we consider the local refinement. By the spectral bounds in \eqref{ineq: spectral-bound-A} for $\A$, following a similar proof of Proposition \ref{prop: TRIP-under-Gaussian-design}, we have TRIP holds for $\scA$ holds with probability at least $1 - \exp(-cp_1)$ when $n \geq C p_1$. Then by Theorem \ref{th: local contraction general setting-RGD}, after $t_{\max} \geq \log\left( \frac{\underline{\lambda} \sqrt{n/df} }{c_1(m) \sigma }   \right)$, we have $ \| \bcX^{t_{\max} } - \bcX^* \|_\F \leq c_3(m) \sigma \sqrt{\frac{df}{n}}$. A similar argument applies to the study of RGN. \quad $\blacksquare$

		\subsection{Proof of Theorem \ref{th: matrix-initialization}}
		Let $\Q_0 \in \bbR^{p_1 \times (r+r^*)}$ be the orthogonal matrix spans the column spaces of $\X^0$ and $\X^*$. Since
		\begin{equation*}
			\|\X^0 - \scA^*(\y)\|_\F^2 = \|\X^0 - P_{\Q_0} (\scA^*(\y)) \|_\F^2 + \| P_{\Q_{0\perp}} (\scA^*(\y)) \|_\F^2
		\end{equation*} and 
		\begin{equation*}
			\|\X^* - \scA^*(\y)\|_\F^2 = \|\X^* - P_{\Q_0} (\scA^*(\y)) \|_\F^2 + \| P_{\Q_{0\perp}} (\scA^*(\y)) \|_\F^2,
		\end{equation*} the SVD property $\|\X^0 - \scA^*(\y)\|_\F^2 \leq \|\X^* - \scA^*(\y)\|_\F^2$ implies that
		\begin{equation*}
			\|\X^0 - P_{\Q_0} (\scA^*(\y)) \|_\F^2 \leq \|\X^* - P_{\Q_0} (\scA^*(\y)) \|_\F^2.
		\end{equation*} 
		Note that following the same proof as \eqref{ineq: Pxt-minus-PAAP}, we have 
		\begin{equation} \label{ineq: spec bnd in init}
			\begin{split}
				\|P_{\Q_0} - P_{\Q_0} \scA^*\scA P_{\Q_0}\| \leq R_{2r}.
			\end{split}
		\end{equation} Hence, 
		\begin{equation} \label{ineq: bound for initialization}
			\begin{split}
				\|\X^0 - \X^*\|_\F &\leq \|\X^0 - P_{\Q_0} (\scA^*(\y)) \|_\F + \|\X^* - P_{\Q_0} (\scA^*(\y)) \|_\F\\
				& \leq 2\|\X^* - P_{\Q_0} (\scA^*(\y)) \|_\F\\
				& \overset{(a)}= 2 \| \X^* - P_{\Q_0} (\scA^*(\scA(\X^*) + \bvarepsilon)) \|_\F\\
				& = 2 \| P_{\Q_0}\X^* - P_{\Q_0} \scA^*\scA(P_{\Q_0}\X^*) -  P_{\Q_0} (\scA^*(\bvarepsilon)) \|_\F\\
				& \leq 2 \left( \| (P_{\Q_0} - P_{\Q_0} \scA^*\scA P_{\Q_0})\X^*\|_\F + \|  P_{\Q_0} (\scA^*(\bvarepsilon)) \|_\F \right)\\
				& \overset{(b)}\leq 2 R_{2r} \|\X^*\|_\F + 2\sqrt{2} \|(\scA^*(\bvarepsilon))_{\max(r)}\|_\F,
			\end{split}
		\end{equation} where (a) is due to the model of $\y$ and (b) is due to that $P_{\Q_0} (\scA^*(\bvarepsilon))$ is a at most rank $2r$ matrix and the spectral norm bound for the operator $(P_{\Q_0} - P_{\Q_0} \scA^*\scA P_{\Q_0})$ in \eqref{ineq: spec bnd in init}.
		
		By Theorem \ref{th: lower-bound} part I, we have $\| (\scA^*(\bvarepsilon))_{\max(r)} \|_\F \leq c'\sqrt{\frac{(p_1 + p_2-r)r}{n}} \sigma$ holds with probability at least $1- \exp(-c\underline{p} )$. At the same time, by Proposition \ref{prop: TRIP-under-Gaussian-design}, there exists $C > 0$ such that when $n \geq C (p_1 + p_2-r)r \|\X^*\|^2_\F/ \sigma^2_{r^*}(\X^*)$, we have $R_{2r} \leq c_1 \frac{\sigma_{r^*}(\X^*)}{\|\X^*\|_\F}$. So for $n$ indicated in the proposition, we have $\|\X^0 - \X^*\|_\F \leq  c \sigma_{r^*} (\X^*) $. Since $\| (\scA^*(\bvarepsilon))_{\max(r)} \|_\F \leq c'\sqrt{\frac{df}{n}} \sigma$, if $ \sigma_{r^*}(\X^*) \leq c_1 \sigma \sqrt{df/n}   $, the estimator is already optimal after initialization, no further refinement is needed. If $ \sigma_{r^*}(\X^*) \geq c_1 \sigma \sqrt{df/n}$, then by Corollary \ref{coro: local contraction general setting-RGD-RGN-matrix}, after $t_{\max} \geq\log\left( \frac{\sigma_{r^*}(\X^*) \sqrt{n/df} }{ c_1\sigma }   \right)$, we have $ \| \bcX^{t_{\max} } - \bcX^* \|_\F \leq c_3 \sigma \sqrt{\frac{df}{n}}$. A similar argument applies to the study of RGN.
		\quad $\blacksquare$

		\subsection{Proof of Theorem \ref{th: initialization-rank-1-tensor-tensor}}
First, notice that in this special setting, $\bcY^1$ and $\bcY^2$ can be compactly represented as $\bcY^j = \lambda \|\w_j\|_2 \cdot \frac{\w_j}{\|\w_j\|_2} \circ \u_{d+1} \circ \cdots \circ \u_{d+m} + \bcE^j$ for $j = 1$ and $2$, where $\bcE^1$ and $\bcE^2$ collect the noise in the first and second half of the data and
\begin{equation} \label{eq: w-formula}
\begin{split}
	&\w_1 = ( \langle \bcA_1, \u_1 \circ \cdots \circ \u_d \rangle, \cdots, \langle \bcA_n, \u_1 \circ \cdots \circ \u_d \rangle )^\top, \\
	&\w_2 = ( \langle \bcA_{n+1}, \u_1 \circ \cdots \circ \u_d \rangle, \cdots, \langle \bcA_{2n}, \u_1 \circ \cdots \circ \u_d \rangle )^\top.
\end{split}
\end{equation} Notice that $\w_1$ and $\w_2$ consist i.i.d. $N(0,1/n)$ random variables. So we have $n \|\w\|^2_2/2 \sim \chi_{n/2}^2$. Then for any $k =d+1, \ldots, d+m$, we have
\begin{equation*}
	\begin{split}
		& \bbP\left( \|\widetilde{\u}_k^{0 \top} \u_{k\perp} \| \geq \sqrt{2}/2  \right) \\
		\leq & \bbP \left( \|\widetilde{\u}_k^{0 \top} \u_{k\perp} \| \geq \sqrt{2}/2 , n \|\w_1\|^2_2 - n \geq -n/16 \right)  + \bbP( n \|\w_1\|^2_2 - n \leq -n/16 )\\
		\overset{(a)}\leq & \bbP \left( \|\widetilde{\u}_k^{0 \top} \u_{k\perp} \| \geq \sqrt{2}/2 \Big| n \|\w_1\|^2_2 - n \geq -n/16 \right) \bbP( n \|\w_1\|^2_2 - n \geq -n/16 )  + \exp(-c n) \\
		\overset{(b)}\leq & \exp(-c \underline{p} ) \bbP( n \|\w_1\|^2_2 - n \geq -n/16 ) + \exp(-c n) \\
		\leq & C \exp(-c \underline{p}),
	\end{split}
\end{equation*}  where (a) is by the standard concentration for $\chi^2$-distribution; (b) is because given $n \|\w_1\|^2_2 - n \geq -n/16$, i.e., $\|\w_1\|_2^2 \geq 15/16$, we have $\lambda \|\w_1\|_2/(\sigma/\sqrt{n} ) \geq C \left[(n \prod_{k=d+1}^{d+m} p_i )^{1/4} + \sqrt{\bar{p}} + \sqrt{n} \right]$ under the condition of $n \geq c(d,m)  \left( \bar{p}  + \frac{\sigma^4}{\lambda^4}\left(\prod^{d+m}_{i=d+1} p_i + \bar{p} \right) \right)$ stated in the theorem and by the Step 1 in the proof of \cite[Theorem]{zhang2018tensor}. By union bound, we have $\bbP\left( \max_{k=d+1, \ldots, d+m} \|\widetilde{\u}_k^{0 \top} \u_{k\perp} \| \geq \sqrt{2}/2  \right) \leq Cm \exp(-c \underline{p})$.

Now, we move onto the guarantee in the second step. Notice that condition on $\{\widetilde{\u}_k^{0} \}_{k=d+1}^{d+m}$, we have
\begin{equation*}
	\y_i' \overset{i.i.d.}\sim \lambda \prod_{k=d+1}^{d+m} ( \widetilde{\u}_k^{0\top } \u_k ) \langle \bcA_i, \u_1 \circ \cdots \circ \u_d \rangle + \epsilon_i', \quad \text{ for } i = n/2+1, \ldots, n,
\end{equation*} and $\epsilon_i' \overset{i.i.d.}\sim N(\sigma^2/n)$. Moreover, with probability  $1-  Cm \exp(-c \underline{p})$,  $\prod_{k=d+1}^{d+m} | \widetilde{\u}_k^{0\top } \u_k | \asymp 1$.  By the proof of \cite[Theorem 4, step 1]{zhang2020islet}, we have with probability at least $1 - \underline{p}^{-C}$ for some $C > 0$ such that when  $n \geq c(d) (\lambda^2+ \sigma^2)  \frac{ (\prod_{i=1}^d p_i  )^{1/2} + \bar{p} }{\lambda^2}$, we have 
		\begin{equation*} 
			\begin{split}
				\|\widetilde{\u}_k^{0 \top} \u_{k\perp}\| \leq \frac{\sqrt{p_k/n} \sqrt{\lambda^2 + \sigma^2} \lambda + (\prod_{k=1}^d p_k)^{1/2} (\lambda^2 + \sigma^2) /n}{\lambda^2} \leq \frac{\sqrt{2}}{2}.
			\end{split}
		\end{equation*} By the union bound, we have condition on  $\{\widetilde{\u}_k^{0} \}_{k=d+1}^{d+m}$ and $\prod_{k=d+1}^{d+m} | \widetilde{\u}_k^{0\top } \u_k | \asymp 1$, $ \max_{k=1,\ldots,d} \|\widetilde{\u}_k^{0 \top} \u_{k\perp} \| \geq \sqrt{2}/2 $ holds with probability at most $Cd p^{-C'}$.  In summary, with probability at least $ 1- \underline{p}^{-C}$, we have $\max_{k=1,\ldots, d+m} \|\widetilde{\u}_k^{0 \top} \u_{k\perp} \| \leq \sqrt{2}/2$.

	Consider applying Theorem \ref{th: perturbation-HOOI} in our setting, we have $\bcT = \bcX^*$, $\widetilde{\bcT} = \scA^*(\bcY)$ and $\bcZ = \scA^*(\bcY) - \bcX^* = \scA^*\scA(\bcX^*) - \bcX^* + \scA^*(\bcE )$. With probability at least $1 - \exp(-c\underline{p} )$, we have 
		\begin{equation*}
			\begin{split}
				\|\bcZ_{\max(\br)}\|_{\F} &= \sup_{\u'_k \in \bbO_{p_k, 1}, k=1,\ldots,d+m} \|\bcZ \times_{k=1}^{d+m} P_{\u'_k} \|_{\F} \\
				&= \sup_{\U'_k \in \bbO_{p_k, 1}, \|\bcW\|_\F \leq 1} \langle \bcZ \times_{k=1}^{d+m} P_{\u'_k}, \bcW \rangle \\
				&= \sup_{\U'_k \in \bbO_{p_k, 1}, \|\bcW\|_\F \leq 1} \langle \bcZ , \bcW \times_{k=1}^{d+m} P_{\u'_k} \rangle\\
				&\leq \sup_{\U'_k \in \bbO_{p_k, 1}, \|\bcW\|_\F \leq 1} \langle \scA^*\scA (\bcX^*) - \bcX^*  , \bcW \times_{k=1}^{d+m} P_{\u'_k} \rangle \\
				& \quad + \sup_{\U'_k \in \bbO_{p_k, 1}, \|\bcW\|_\F \leq 1} \langle \scA^*(\bcE )  , \bcW \times_{k=1}^{d+m} P_{\u'_k} \rangle \\
				& \overset{ \text{Lemma } \ref{lm:retricted orthogonal property}, \text{Theorem } \ref{th: lower-bound} } \leq  R_{2\br} \|\bcX^*\|_{\F} + C \sigma\sqrt{\frac{\sum_{i=1}^{d+m} p_i  }{n}}\\
				& \overset{ \text{Proposition }\ref{prop: TRIP-under-Gaussian-design} } \leq C( \lambda + \sigma ) \sqrt{\frac{\sum_{i=1}^{d+m} p_i }{n}}.
			\end{split}
		\end{equation*} 
		
		So when $n$ satisfies the condition indicated in the proposition, we have $\|\bcZ_{\max(\br)}\|_{\F} \leq c \underline{\lambda}$. Thus by Theorem \ref{th: perturbation-HOOI}, we have the output of the algorithm satisfies
		\begin{equation*}
			\|\bcX^0 - \bcX^*\|_{\F} \leq c'(d) \|\bcZ_{\max(\br)}\|_{\F} \leq c \lambda
		\end{equation*} with probability at least $1 - \underline{p}^{-C}$.

		Since $ \|(\scA^*(\bcE ))_{\max(2\br)}\|_{\F} \leq c \sigma \sqrt{df/n} $ by Theorem \ref{th: lower-bound}, when $ \lambda\leq c_1 \sigma \sqrt{df/n}   $, the estimator is already optimal after initialization, no further refinement is needed. If $ \lambda\geq c_1 \sigma \sqrt{df/n}$, TRIP holds for $n$ indicated in the proposition, then by Theorem \ref{th: local contraction general setting-RGD}, after $t_{\max} \geq \log\left( \frac{\lambda \sqrt{n/df} }{ c_1(d,m)\sigma }   \right)$, we have $ \| \bcX^{t_{\max} } - \bcX^* \|_\F \leq c_3(d,m) \sigma \sqrt{\frac{df}{n}}$. A similar argument applies to the study of RGN.

		\section{Proofs in Section \ref{sec: computational-limit}} \label{sec: proof-comp-limit}
		
		\subsection{An Equivalence Formulation For Scalar-on-tensor Regression} \label{sec: equiva-scalar-on-tensor-reg}
		In this section, we show without loss of generality, we can assume $\bvarepsilon_i \overset{i.i.d.}\sim N(0,\sigma^2)$ with $0 \leq \sigma^2 < 1$, $\bcA_i \overset{i.i.d.}\sim N(0,1)$ and $\|\bcX^*\|_\F + \sigma^2 = 1$ in establishing the computational lower bound for scalar-on-tensor regression under Gaussian design. Consider \eqref{eq: scalar-on-tensor-model} and suppose we are in a simpler setting that $\|\bcX^*\|_\F$ and $\sigma$ are known. Then we can rescale the problem by multiplying $\sqrt{n/(\|\bcX^*\|_\F + \sigma^2)}$ on both sides of \eqref{eq: scalar-on-tensor-model} and get 
		\begin{equation} \label{eq: scalar-on-tensor-model2}
			\y_i' = \langle \bcA', \bcX^{*'} \rangle + \bvarepsilon_i', \quad i=1,\ldots, n,
		\end{equation} 
		where $ \y_i' = \sqrt{n}\y_i/\sqrt{ \|\bcX^*\|_\F + \sigma^2 } \sim N(0,1) $, $\bcA' = \sqrt{n} \bcA $ has i.i.d. $N(0,1)$ entries and $\bcX^{*'} = \bcX^*/\sqrt{ \|\bcX^*\|_\F + \sigma^2 }$, $\bvarepsilon_i'$ follows i.i.d. $N(0, \sigma^2/\sqrt{ \|\bcX^*\|_\F + \sigma^2 }  )$ satisfying $\|\bcX^{*'}\|_\F^2 + \var(\bvarepsilon_i' ) = 1$. It is not hard to see that any lower bound established for estimating $\bcX^{*'}$ based on $\{ \y_i', \bcA' \}_{i=1}^n$, say $c_{\textrm{lower}}$, also implies the lower bound for estimating $\bcX^*$. This is because for any estimator $\widehat{\bcX}$ of $\bcX^*$, we have
		\begin{equation*}
			\| \widehat{\bcX} - \bcX^* \|_\F = \sqrt{ \|\bcX^*\|_\F + \sigma^2 } \left\| \frac{\widehat{\bcX}}{\sqrt{ \|\bcX^*\|_\F + \sigma^2 }} - \bcX^{*'} \right\|_\F \geq c_{\textrm{lower}}  \sqrt{ \|\bcX^*\|_\F + \sigma^2 }.
		\end{equation*} Thus, without loss of generality, we can assume $\bvarepsilon_i \overset{i.i.d.}\sim N(0,\sigma^2)$ with $0 \leq \sigma^2 < 1$, $\bcA_i \overset{i.i.d.}\sim N(0,1)$ and $\|\bcX^*\|_\F + \sigma^2 = 1$.

		\subsection{More Backgroud and Preliminaries for Low-degree Polynomials Method} \label{sec: preliminary-low-degree-hermite}
		
		Next, we provide a few preliminary facts for low-degree polynomials methods. Given data $X$, consider the simple hypothesis testing problem: $H_0$ v.s. $H_1$. We have the following result for the low-degree likelihood ratio.
		\begin{Proposition}[Page 35 of \cite{hopkins2018statistical} or \cite{kunisky2019notes} Proposition 1.15]\label{prop: trun likelihood ration}
			Let likelihood ratio be $\LR(x) = \frac{p_{H_1}(x)}{p_{H_0} (x)}: \Omega^n \to \bbR$. For every $D \in \bbN$, we have 
			\begin{equation*}
				\frac{ \LR^{\leq D} - 1}{\|\LR^{\leq D} - 1\|} = \argmax_{\substack{f : f \text{ has degree at most }D\\ \bbE_{H_0} f^2(X) = 1,\bbE_{H_0} f(X) = 0}} \bbE_{ H_1} f(X)
			\end{equation*}
			and 
			\begin{equation*}
				\|\LR^{\leq D} - 1\| = \max_{\substack{f:f \text{ has degree at most }D\\ \bbE_{H_0} f^2(X) = 1,\\ \bbE_{H_0} f(X) = 0}} \bbE_{H_1} f(X),
			\end{equation*}
			where $\|f\| = \sqrt{\bbE_{H_0}f^2(X)}$ and $f^{\leq D}$ is the projection of a function $f$ to the linear subspace of degree-$D$ polynomials, where the projection is orthonormal with respect to the inner product induced under $H_0$. 
		\end{Proposition}
		Here the key quantity we are interested to bound is $\|\LR^{\leq D} - 1\|$. Suppose $D \geq 1$ is fixed, $f_0, f_1, \ldots, f_q: \Omega^n \to R$ are orthonormal basis for degree $D$ functions (with respect to $\langle \cdot, \cdot \rangle_{H_0}$), and that $f_0(x) = 1$ is a constant function. Then by the property of basis functions, we have
		\begin{equation} \label{eq: evalutation-truncated-likelihood-ratio}
			\begin{split}
				\|\LR^{\leq D} - 1\|^2 & = \sum_{i=1}^q \left(  \bbE_{H_0}(f_i(X) (\LR^{\leq D} - 1) )  \right)^2\\
				& \overset{(a)}= \sum_{i=1}^D \left(  \bbE_{H_0}(f_i(X) \LR(X) )  \right)^2 = \sum_{i=1}^D (\bbE_{H_1} f_i(X))^2,
			\end{split}
		\end{equation} here (a) is because $\LR - (\LR^{\leq D} - 1)$ is orthogonal to $f_i$ for $ i \in [1,D]$ by assumption. So the typical main task in the low-degree polynomials method boils down to bound $\sum_{i=1}^D (\bbE_{H_1} f_i(X))^2$. In the hypothesis testing formulation for scalar-on-tensor regression \eqref{eq: hypo-test-scalar-on-tensor}, the data are i.i.d. Gaussian under the null hypothesis. So a natural choice for the basis functions are Hermite polynomials, which are orthogonal polynomials with respect to the Gaussian measure \citep{szeg1939orthogonal}.
		
		Recall $\bbN = \{0,1,2,\ldots \}$, let $\{h_k \}_{k \in \bbN}$ be the normalized univariate Hermite polynomials $h_k = \frac{1}{\sqrt{k!}} H_k$ where $\{ H_k \}_{k \in \bbN}$ are univariate Hermite polynomials which are defined by the following recurrence:
		\begin{equation*}
			H_0(x) = 1, \quad H_1(x) = x, \quad H_{k+1}(x) = x H_k(x) - k H_{k-1}(x) \quad \text{  for } k \geq 1.
		\end{equation*}
		The normalized univariate Hermite polynomials satisfy the following key property:
		\begin{equation*}
			\begin{split}
				\bbE_{Z \sim N(0,1)} [h_k(Z)] = 0 \quad \forall k \geq 0, \quad \text{ and } \quad \bbE_{Z \sim N(0,1)} [ h_{k_1}(Z) h_{k_2} (Z) ] = 1(k_1 = k_2 ),
			\end{split}
		\end{equation*} where $1(\cdot)$ is the indicator function. In addition, for $\balpha = (\alpha_1,\ldots,\alpha_n) \in \bbN^n $ and $\Z = (Z_1,\ldots,Z_n) \in \bbR^n $, let $h_{\balpha}(\Z) = \prod_{i \in [n]} h_{\alpha_i}(Z_i) $ be the normalized $n$-variate Hermite polynomials. They form an orthogonal basis with respect to $N(0,1)^{\otimes n}$, i.e., if $\Z$ has i.i.d. $N(0,1)$ entries, then $\bbE[ h_{\balpha}(\Z) h_{\bbeta}(\Z) ] = 1(\balpha = \bbeta)$. Next, we introduce two additional properties regards Hermite polynomials.
		\begin{Lemma}(Expansion for Shifted Hermite Polynomials (\cite[Proposition 3.1]{schramm2020computational})) \label{lm: hermite-poly-expansion} For any $k \in \bbN$ and $z, \mu \in \bbR$, then
			\begin{equation*}
				H_k(z + \mu) = \sum_{l=0}^k {k \choose l} \mu^{k-l} H_{l}(z), 
			\end{equation*} and 
			\begin{equation*}
				h_k(z + \mu) = \sum_{l=0}^k \sqrt{ \frac{l!}{k!} }   {k \choose l} \mu^{k-l} h_{l}(z).
			\end{equation*}
		\end{Lemma} 
		
		\begin{Lemma}(Gaussian Integration by Parts \cite[Proposition 2.10]{kunisky2019notes}) \label{lm: hermite-poly-integration-by-part} If $f: \bbR \to \bbR$ is a $k$ ($k \geq 0$) times continuously differentiable and $f(y)$ and its first $k$ derivatives are bounded by $O(\exp(|y|^\alpha))$ for some $\alpha \in (0,2)$, then
			\begin{equation*}
				\bbE_{Y \sim N(0,1)}[ h_k(Y) f(Y) ] = \frac{1}{\sqrt{k!}} \bbE_{Y \sim N(0,1)} \left[ \frac{d^k f}{dY^k}(Y) \right].
			\end{equation*}
		\end{Lemma}

		\subsection{Proof of Theorem \ref{th: lower-degree-polynomial-tensor-regression} } \label{sec: low-degree-lower-bound-proof}
		First, it is easy to check under $H_1$, we have
		\begin{equation} \label{eq: joint-distribution-H1}
			(\y_i, \rmvec(\bcA_i) )| \bcX^* \overset{i.i.d.}\sim N\left(0, \left( \begin{array}{cc}
				1 & (\rmvec(\bcX^*))^\top\\
				\rmvec(\bcX^*) & \I_{p^d}
			\end{array}  \right) \right).
		\end{equation} Moreover, let us denote the distribution of $\bcX^*$ under $H_1$ as $\cP_{\bcX^*}$. Since the data is i.i.d. standard Gaussian under $H_0$, Hermite polynomials are a natural choice for the orthogonal polynomial basis under $H_0$. Let the degrees of $\{\y_i, \bcA_i \}_{i=1}^n$ to be $\{\alpha_i, \bbeta_i \}_{i=1}^n$ defined in the following way:
		\begin{equation*}
			\balpha \in \bbN^n = \{\alpha_i \}_{i=1}^n, \quad \bbeta \in \bbN^{np^d} = \{ \beta_{i,\j  } \}_{i \in [n], \j = (j_1,\ldots,j_d)  \in [p]^{\otimes d}}, \quad \bbeta_i = \{ \beta_{i,\j} \}_{\j = (j_1,\ldots,j_d)  \in [p]^{\otimes d} }.
		\end{equation*} Finally, throughout the proof, given any vector $\a = (a_1,\ldots,a_q) \in \bbR^q, \b = (b_1,\ldots,b_q) \in \bbR^q $, let $|\a| = \sum_{j=1}^q a_j$ and $ \a! = \prod_{j=1}^q a_j! $, $\b^\a = \prod_{j=1}^q b_j^{a_j} $.
		
		By Proposition \ref{prop: trun likelihood ration} and the properties of Hermite polynomials introduced in Appendix \ref{sec: preliminary-low-degree-hermite}, we have
		\begin{equation} \label{eq: low-degree-main-bound}
			\begin{split}
				\sup_{\text{polynomial } f: ~ \substack{deg(f) \leq D  \\
						\mathbb{E}_{H_0} f(\{y_i, \bcA_i\}_{i=1}^n) = 0,\\ \text{Var}_{H_0} f(\{y_i, \bcA_i\}_{i=1}^n) = 1}}\mathbb{E}_{H_1}f(\{y_i, \bcA_i\}_{i=1}^n) =& \| \text{LR}^{\leq D} -1 \| \\
				\overset{ \eqref{eq: evalutation-truncated-likelihood-ratio} } =& \sqrt{ \sum_{(\balpha, \bbeta): 1\leq |\balpha| + |\bbeta| \leq D}  \left(\bbE_{H_1} h_{\balpha,\bbeta}( \{\y_i,\bcA_i \}_{i=1}^n ) \right)^2 }. 
			\end{split}
		\end{equation} 
		Then
		\begin{equation} \label{eq: low-degree-poly-exp-com}
			\begin{split}
				\bbE_{H_1} \left( h_{\balpha,\bbeta}( \{\y_i,\bcA_i \}_{i=1}^n )\right)
				\overset{(a)}=&\bbE_{\bcX^*\sim \cP_{\bcX^*}}  \left(\prod_{i=1}^n  \bbE \left( h_{\alpha_i,\bbeta_i}( \y_i,\bcA_i )| \bcX^*\right) \right) \\
				\overset{ \eqref{eq: joint-distribution-H1}, \text{Lemma }\ref{lm:correlated-hermitian} }	= &\bbE_{\bcX^*\sim \cP_{\bcX^*}}  \left( \prod_{i=1}^n  \sqrt{ \frac{\alpha_i!}{\bbeta_i!} } \bcX^{* \bbeta_i} 1(\alpha_i = |\bbeta_i|)   \right),
			\end{split}
		\end{equation} where (a) is because condition on $\bcX^*$, $(\y_i, \bcA_i)$ are independent.
		
		We divide the rest of the proof into two steps.

		{\bfseries Step 1} In this step, we bound the right-hand side of \eqref{eq: low-degree-main-bound}.
		\begin{equation} \label{eq: low-degree-second-argument}
			\begin{split}
				&\sum_{(\balpha, \bbeta): |\balpha| + |\bbeta| \leq D}  \left(\bbE_{H_1} h_{\balpha,\bbeta}( \{\y_i,\bcA_i \}_{i=1}^n ) \right)^2 \\
				\overset{ \eqref{eq: low-degree-poly-exp-com} }	=& \sum_{(\balpha, \bbeta): \alpha_i = |\bbeta_i|, \forall i \in [n], |\balpha|\leq D/2}  \left( \bbE_{\bcX^*\sim \cP_{\bcX^*}}  \left( \prod_{i=1}^n  \sqrt{ \frac{\alpha_i!}{\bbeta_i!} } \bcX^{* \bbeta_i}  \right) \right)^2 \\
				\overset{(a)}=& \sum_{(\balpha, \bbeta): \alpha_i = |\bbeta_i|, \forall i \in [n], |\balpha|\leq D/2} \left(\bbE_{\bcX^*\sim \cP_{\bcX^*}}  \left( \prod_{i=1}^n  \frac{\alpha_i!}{\bbeta_i!} \bcX_1^{* \bbeta_i} \bcX_2^{*\bbeta_i} \right) \right) \\
				= &  \sum_{\balpha: |\balpha| \leq D/2}\bbE_{\bcX_1^*, \bcX_2^* \sim \cP_{\bcX^*} } \left( \sum_{\bbeta_1: |\bbeta_1| = \alpha_1} \cdots \sum_{\bbeta_n: |\bbeta_n| = \alpha_n}  \prod_{i=1}^n  \frac{\alpha_i!}{\bbeta_i!}\bcX_1^{* \bbeta_i} \bcX_2^{*\bbeta_i}  \right) \\
				= & \sum_{\balpha: |\balpha| \leq D/2}\bbE_{\bcX_1^*, \bcX_2^* \sim \cP_{\bcX^*} } \left( \prod_{i=1}^n \left( \sum_{\bbeta_i: |\bbeta_i| = \alpha_i}  \frac{\alpha_i!}{\bbeta_i!}\bcX_1^{* \bbeta_i} \bcX_2^{*\bbeta_i} \right) \right) \\
				\overset{(b)} =& \sum_{\balpha: |\balpha| \leq D/2}\bbE_{\bcX_1^*, \bcX_2^* \sim \cP_{\bcX^*} } \left( \prod_{i=1}^n  \langle \bcX^*_1, \bcX^*_2 \rangle^{\alpha_i}  \right) \\
				\overset{(c)}=& \sum_{\balpha: |\balpha| \leq D/2}  \bbE_{\bcX_1^*, \bcX_2^* \sim \cP_{\bcX^*} } \left( \left((1-\sigma^2) \langle \x^*_1, \x^*_2 \rangle^d \right)^{|\balpha|}  \right).
			\end{split}
		\end{equation} Here (a) is by performing the ``replica'' manipulation where $\bcX_1^*$, $\bcX_2^*$ are drawn independently from $\cP_{\bcX^*}$; (b) is by the sum of multinomials; (c) is by the generating process of $\bcX_1^*$ and $\bcX_2^*$.
		
		Given $(\balpha, \bbeta)$, let $I_{\balpha} = \{i \in [n]: \alpha_i \neq 0  \}$ and $|I_{\balpha} |$ be the cardinality of the set $I_{\balpha}$. Also notice that since $\x^*_1$ and $\x^*_2$ has i.i.d. entries generated from Uniform$( 1/\sqrt{p}, -1/\sqrt{p} )$, we have $p \langle \x^*_1, \x^*_2 \rangle$ has i.i.d. Uniform$( 1, -1 )$ entries. By the property of sub-gaussian random variable, for any integer $k \geq 1$, we have
		\begin{equation} \label{ineq: sub-gaussian-abs-bound}
			\bbE(|p \langle \x^*_1, \x^*_2 \rangle|^{k}) \leq p^{k/2} k^{k/2}
		\end{equation} by \cite[Eq. (5.11)]{vershynin2010introduction} and the fact $p \langle \x^*_1, \x^*_2 \rangle$ has sub-gaussian norm $\sqrt{p}$. Then we have
		\begin{equation} \label{ineq: low-degree-third-argument}
			\begin{split}
				&\sum_{(\balpha, \bbeta): |\balpha| + |\bbeta| \leq D}  \left(\bbE_{H_1} h_{\balpha,\bbeta}( \{\y_i,\bcA_i \}_{i=1}^n ) \right)^2 \\
				&\overset{\eqref{eq: low-degree-second-argument}} \leq \sum_{\balpha: |\balpha| \leq D/2} \frac{(1-\sigma^2)^{|\balpha|}}{p^{d|\balpha|}}  \bbE_{\bcX_1^*, \bcX_2^* \sim \cP_{\bcX^*} } \left( \left(  |p\langle \x^*_1, \x^*_2 \rangle| \right)^{d|\balpha|}  \right) \\
				& \overset{ \eqref{ineq: sub-gaussian-abs-bound} } \leq \sum_{\balpha: |\balpha| \leq D/2}  \frac{(1-\sigma^2)^{|\balpha|}}{p^{d|\balpha|}} p^{d|\balpha|/2} (d|\balpha|)^{d|\balpha|/2}\\
				& \leq \sum_{\widebar{D} = 1}^{D/2} \sum_{a=1}^{ \widebar{D} } \sum_{\balpha: |\balpha| = \widebar{D}, |I_{\balpha}| = a } \frac{(1-\sigma^2)^{\widebar{D}}}{p^{d\widebar{D}}} p^{d\widebar{D}/2} (d\widebar{D})^{d\widebar{D}/2} \\
				& \overset{(a)}\leq \sum_{\widebar{D} = 1}^{D/2} \sum_{a=1}^{ \widebar{D} } {n \choose a} {\widebar{D} - 1 \choose a-1 }   \frac{(1-\sigma^2)^{ \widebar{D} }}{p^{d \widebar{D} /2}}  (d\widebar{D})^{d \widebar{D} /2}
			\end{split}
		\end{equation} here (a) is because the set $\{ \balpha: |I_{\balpha}| = a, |\balpha| = \widebar{D} \}$ has cardinality at most ${n \choose a} {\widebar{D} - 1 \choose a-1 }$. 
		
		{\bfseries Step 2}. In this step, we bound $\sum_{a=1}^{ \widebar{D} } {n \choose a} {\widebar{D} - 1 \choose a-1 } $ and complete the proof.
		\begin{equation} \label{ineq: bound-a-sum-part}
			\begin{split}
				&\sum_{a=1}^{ \widebar{D} } {n \choose a} {\widebar{D} - 1 \choose a-1 } \\
				\leq & \sum_{a=1}^{ \widebar{D} } n^a {\widebar{D} - 1 \choose a-1 } \leq n^{\widebar{D}}  \sum_{a=1}^{ \widebar{D} }  {\widebar{D} - 1 \choose a-1 } \leq (2n)^{\widebar{D}}.
			\end{split}
		\end{equation}
		
		By plugging \eqref{ineq: bound-a-sum-part} into \eqref{ineq: low-degree-third-argument}, we have
		\begin{equation*}
			\begin{split}
				\sum_{(\balpha, \bbeta): |\balpha| + |\bbeta| \leq D}  \left(\bbE_{H_1} h_{\balpha,\bbeta}( \{\y_i,\bcA_i \}_{i=1}^n ) \right)^2 &\leq  \sum_{\widebar{D} = 1}^{D/2} \left( 2n (1-\sigma^2) (d \widebar{D}/p )^{d/2} \right)^{\widebar{D}} \\
				&\overset{(a)}\leq  \sum_{\widebar{D} = 1}^{D/2} \delta^{\widebar{D}} \leq \frac{\delta}{1-\delta},
			\end{split}
		\end{equation*} here (a) is due to $n \leq \frac{\delta}{2(1-\sigma^2) } (p/dD)^{d/2}$ by our assumption. The proof is finished by observing \eqref{eq: low-degree-main-bound}. \quad $\blacksquare$

		\subsection{Hardness of Hypothesis Testing Implies Hardness of Estimation} \label{sec: supp-test-to-estimation}
		\begin{Proposition}[Hardness of Hypothesis Testing Implies Hardness of Estimation] \label{prop: testing-imply-estimation} If there does not exist a polynomial-time tester distinguishing between $H_0$ and $H_1$ in scalar-on-tensor regression  with Type I + II error tending to zero as $n \to \infty$, then there is no polynomial-time estimator $\widehat{\bcX}$ such that $\| \widehat{\bcX} - \bcX^* \|_\F \leq \frac{1}{4} \|\bcX^*\|_\F $ as $n \to \infty$.
		\end{Proposition}
		We prove this by using the sample splitting and contradiction argument. Suppose $\bcX^*$ and $\{\y_i, \bcA_i \}_{i=1}^{2n}$ are generated under $H_1$ described in \eqref{eq: hypo-test-scalar-on-tensor}. 
		
		If there is polynomial-time estimator $\widehat{\bcX}$ based on $\{\y_i,\bcA_i \}_{i=1}^n$ such that $\|\widehat{\bcX} - \bcX^*\|_\F \leq \frac{1}{4}\|\bcX^*\|_\F $, then 
		\begin{equation*}
			\frac{3}{4}\|\bcX^*\|_\F \leq \|\bcX^*\|_\F - \|\bcX^* - \widehat{\bcX}\|_\F \leq \| \widehat{\bcX} \|_\F \leq \|\bcX^*\|_\F+ \|\bcX^* - \widehat{\bcX}\|_\F \leq \frac{5}{4} \|\bcX^*\|_\F.
		\end{equation*}
		
		Given the first half of the data, we consider performing the test based on the statistic
		\begin{equation*}
			T = \sum_{i=n+1}^{2n} (\y_i^2 - 1) \langle \bcA_i, \widehat{\bcX}\rangle^2.
		\end{equation*}
		
		Under $H_0$, $\bbE_{H_0}(T|\widehat{\bcX}) = 0$ and
		\begin{equation*}
			\begin{split}
				\var_{H_0}(T|\widehat{\bcX}) = n \var_{H_0}( (\y_i^2-1)\langle \bcA_i, \widehat{\bcX}\rangle^2 )& = n \bbE_{H_0}((\y_i^2-1)^2) \bbE_{H_0}(\langle \bcA_i, \widehat{\bcX}\rangle^2) \\
				&= 2n \|\widehat{\bcX}\|_\F^2  \leq \frac{25}{8}n \|\bcX^*\|^2_\F.
			\end{split}
		\end{equation*}
		
		Under $H_1$, we have
		\begin{equation*}
			\begin{split}
				&\bbE_{H_1}(T| \widehat{\bcX} )\\
				=& n \bbE_{H_1} \left( \left( \langle \bcA_i , \bcX^* \rangle^2 + \bvarepsilon_i^2 + 2 \bvarepsilon_i \langle \bcA_i , \bcX^* \rangle - 1 \right) \langle \bcA_i, \widehat{\bcX}\rangle^2  \right)\\
				=& n \bbE_{H_1} \left( \langle \bcA_i , \bcX^* \rangle^2 \langle \bcA_i, \widehat{\bcX}\rangle^2  \right)\\
				=& n \bbE_{H_1} \Bigg\{ \left( \sum_{ \j \in [p]^{\otimes d} } \bcA^2_{i,\j} \bcX^{*2}_{\j} + \sum_{ \substack{\j \in [p]^{\otimes d}; \z  \in [p]^{\otimes d} \\ \j \neq \z } } \bcA_{i,\j} \bcA_{i,\z} \bcX^*_{\j} \bcX^*_{\z}  \right)\\
				& \cdot \left(  \sum_{ \j \in [p]^{\otimes d} } \bcA^2_{i,\j} \widehat{\bcX}^{2}_{\j} + \sum_{ \substack{\j \in [p]^{\otimes d}; \z  \in [p]^{\otimes d} \\ \j \neq \z } } \bcA_{i,\j} \bcA_{i,\z} \widehat{\bcX}_{\j} \widehat{\bcX}_{\z} \right)   \Bigg\}\\
				=& n \left(  \sum_{ \substack{\j \in [p]^{\otimes d}; \z  \in [p]^{\otimes d} \\ \j \neq \z } } \bcX^{*2}_{\j} \widehat{\bcX}^2_{\z} + 3 \sum_{\j \in [p]^{\otimes d}}  \bcX^{*2}_{\j}\widehat{\bcX}^{2}_{\j} +  2\sum_{ \substack{\j \in [p]^{\otimes d}; \z  \in [p]^{\otimes d} \\ \j \neq \z } } \bcX^*_{\j} \bcX^*_{\z} \widehat{\bcX}_{\j} \widehat{\bcX}_{\z} \right)\\
				= & n ( \|\bcX^*\|_\F^2\|\widehat{\bcX}\|_\F^2 + 2\langle \bcX^*, \widehat{\bcX} \rangle^2 )\\
				\overset{(a)}\geq & \frac{27}{16}n \|\bcX^*\|_\F^4,
			\end{split}
		\end{equation*} here (a) is because $\|\widehat{\bcX}\|_\F^2 \geq \frac{9}{16} \|\bcX^*\|_\F^2 $ and 
		\begin{equation*}
			\begin{split}
				\langle \bcX^*, \widehat{\bcX} \rangle =& \frac{1}{4}( \| \bcX^* + \widehat{\bcX} \|_\F^2 - \| \bcX^* - \widehat{\bcX} \|_\F^2 ) \\
				\geq & \frac{1}{4}\left( (\|2\bcX^*\|_\F - \| \widehat{\bcX} - \bcX^* \|_\F )^2 - \frac{1}{16}\|\bcX^*\|_\F^2  \right) \\
				\geq & \frac{1}{4}\left( (7/4)^2 \|\bcX^*\|_\F^2 - 1/16 \|\bcX^*\|_\F^2 \right) = \frac{3}{4} \|\bcX^*\|_\F^2.
			\end{split}
		\end{equation*}
		Moreover, $\var_{H_1}(T|\widehat{\bcX}) = n \var_{H_1}( (\y_i^2-1)\langle \bcA_i, \widehat{\bcX}\rangle^2 ) = O(n)$.
		
		So condition on $\{\y_i, \bcA_i \}_{i=1}^n$, we consider the test of rejecting if $T \geq \frac{27}{32} n\|\bcX^*\|_\F^4 $, not rejecting otherwise. Then
		\begin{equation*}
			\begin{split}
				\text{Type I error} &= \bbP_{H_0}( T \geq \frac{27}{32} n\|\bcX^*\|_\F^4|\widehat{\bcX}) \overset{(a)}\leq \frac{\var_{H_0}(T|\widehat{\bcX})}{(\frac{27}{32} n\|\bcX^*\|_\F^4)^2} \overset{n \to \infty}\to 0, \\
				\text{Type II error} &= \bbP_{H_1} (T \leq \frac{27}{32} n\|\bcX^*\|_\F^4|\widehat{\bcX}) = \bbP_{H_1} (T - \frac{27}{16} n\|\bcX^*\|_\F^4 \leq -\frac{27}{32} n\|\bcX^*\|_\F^4|\widehat{\bcX} ) \\
				& \leq \frac{\var_{H_1}(T|\widehat{\bcX})}{(\frac{27}{32} n\|\bcX^*\|_\F^4)^2} \overset{n \to \infty}\to 0,
			\end{split}
		\end{equation*} where (a) is by the Chebysev's inequality.
		
		So this polynomial-time test based on $T$ achieves distinguishing between $H_0$ and $H_1$ with Type I + II errors goes to $0$ as $n \to \infty$. This contradicts our assumption. So there does not exist a polynomial-time estimator $\widehat{\bcX}$ such that $\|\widehat{\bcX} - \bcX^*\|_\F \leq \frac{1}{4}\|\bcX^*\|_\F $ and this finishes the proof. \quad $\blacksquare$

		\section{Proofs in Section \ref{sec:technical contribution} } \label{sec: proof-techni-contribution}
		
		\subsection{Proof of Lemma \ref{lm: orthogonal projection}}
		Suppose $\bcX^t$ has Tucker rank $\br$ decomposition $\llbracket \bcS^t; \U_1^t,\ldots,\U_{d+m}^t \rrbracket $ and $\bcX^*$ has Tucker rank $\br^*$ decomposition $\llbracket \bcS; \U_1,\ldots,\U_{d+m} \rrbracket $, respectively. Recall 
		\begin{equation*} 
			\begin{split}
				\W_k  := (\U_{d+m} \otimes \cdots \otimes \U_{k+1} \otimes \U_{k-1} \otimes \ldots \otimes \U_1) \V_k \in \bbO_{p_{-k}, r^*_{k}};\\
				\W^t_k  := (\U^t_{d+m} \otimes \cdots \otimes \U^t_{k+1} \otimes \U^t_{k-1} \otimes \ldots \otimes \U^t_1) \V^t_k \in \bbO_{p_{-k}, r_{k}},
			\end{split}
		\end{equation*} where $p_{-k} = \prod_{i\neq k} p_i$, $\V_k = \QR(\cM_k(\bcS)^\top), \V_k^t = \QR(\cM_k(\bcS^t)^\top)$. For $\bcX^*$, it can be decomposed in the following way
		\begin{equation} \label{eq: decom of Xstar}
			\begin{split}
				\bcX^* = & \bcX^* \times_1 P_{\U_{1\perp}^t} + \bcX^* \times_1 P_{\U_{1}^t} \times_2 P_{\U_{2\perp}^t} + \cdots +  \bcX^* \times_{l=1}^{k-1} P_{\U_{l}^t} \times_k P_{\U_{k\perp}^t} + \cdots + \bcX^* \times_{l=1}^{d+m} P_{\U_{l}^t}\\
				=& \sum_{k=1}^{d+m} \bcX^* \times_{l=1}^{k-1} P_{\U_{l}^t} \times_k P_{\U_{k\perp}^t} + \bcX^* \times_{l=1}^{d+m} P_{\U_{l}^t}.
			\end{split}
		\end{equation}
		For $k = 1,\ldots, d$, let us denote $\U_k^t = [\widebar{\U}_k^t \quad \widecheck{\U}_k^t  ]$ where $\widebar{\U}_k^t$ is composed of the first $r_k^*$ columns of $\U_k^t$ and $\widecheck{\U}_k^t$ is composed of the remaining $(r_k - r_k^*) $ columns of $\U_k^t$. 
		
		Then
		\begin{equation} \label{eq: P-Tperp Xstar}
			\begin{split}
				P_{(T_{\bcX^t})_\perp} \bcX^* =& \bcX^* - P_{T_{\bcX^t}} \bcX^*\\
				\overset{\eqref{eq: tangent space projector}}=	& \bcX^* - (\bcX^* \times_{k=1}^{d+m} P_{\U_k^t} + \sum_{k=1}^{d+m} \cT_k( P_{\U_{k\perp}^t} \cM_k(\bcX^*) P_{\W_k^t} )  )\\
				\overset{\eqref{eq: decom of Xstar}}= & \sum_{k=1}^{d+m} \left(\bcX^* \times_{l=1}^{k-1} P_{\U_{l}^t} \times_k P_{\U_{k\perp}^t} - \cT_k( P_{\U_{k\perp}^t} \cM_k(\bcX^*) P_{\W_k^t} ) \right)\\
				\overset{ \eqref{eq: matricization relationship} }= &  \sum_{k=1}^{d+m} \left( \cT_k \left( P_{\U_{k\perp}^t} \cM_k(\bcX^*) ( \otimes_{l=d+m}^{k+1}\I_{p_l}\otimes_{l=k-1}^{1} P_{\U_l^t} - P_{\W_k^t}  )   \right)  \right) \\
				= & \sum_{k=1}^{d+m} \left( \cT_k \left( (\I_{p_k} - P_{\U_{k}^t}) \cM_k(\bcX^*) ( \otimes_{l=d+m}^{k+1}\I_{p_l}\otimes_{l=k-1}^{1} P_{\U_l^t} - P_{\W_k^t}  )   \right)  \right)\\
				\overset{ \eqref{eq: subspace-decomposition} } = & \sum_{k=1}^{d+m} \left( \cT_k \left(  (\I_{p_k} - P_{\widecheck{\U}_k^t} ) (\I_{p_k} - P_{\widebar{\U}_k^t} ) \cM_k(\bcX^*) ( \otimes_{l=d+m}^{k+1}\I_{p_l}\otimes_{l=k-1}^{1} P_{\U_l^t} - P_{\W_k^t}  )   \right)  \right)\\
				\overset{(a)} = & \sum_{k=1}^{d+m} \left( \cT_k \left(  (\I_{p_k} - P_{\widecheck{\U}_k^t} ) (P_{\U_k} - P_{\widebar{\U}_k^t} ) \cM_k(\bcX^* ) ( \otimes_{l=d+m}^{k+1}\I_{p_l}\otimes_{l=k-1}^{1} P_{\U_l^t} - P_{\W_k^t}  )   \right)  \right)\\
				\overset{(b)} = & \sum_{k=1}^{d+m} \left( \cT_k \left(  (\I_{p_k} - P_{\widecheck{\U}_k^t} ) (P_{\U_k} - P_{\widebar{\U}_k^t} ) \cM_k(\bcX^* - \bcX^t) ( \otimes_{l=d+m}^{k+1}\I_{p_l}\otimes_{l=k-1}^{1} P_{\U_l^t} - P_{\W_k^t}  )   \right)  \right),
			\end{split}
		\end{equation} here (a) is because the $\U_k$ spans the column space of $\cM_k(\bcX^*)$, (b) is because $\cM_k(\bcX^t) ( \otimes_{l=d+m}^{k+1}\I_{p_l}\otimes_{l=k-1}^{1} P_{\U_l^t} - P_{\W_k^t} ) = 0$.
		
		It is easy to check $\otimes_{l=d+m}^{k+1}\I_{p_l}\otimes_{l=k-1}^{1} P_{\U_l^t} - P_{\W_k^t}$ is a projection matrix. So from \eqref{eq: P-Tperp Xstar}, we have
		\begin{equation*}
			\begin{split}
				\|P_{(T_{\bcX^t})_\perp} \bcX^*\|_{\F}\leq &  \sum_{k=1}^{d+m} \| \cT_k \left( (P_{\U_k} - P_{\widebar{\U}_k^t} ) \cM_k(\bcX^* - \bcX^t) ( \otimes_{l=d+m}^{k+1}\I_{p_l}\otimes_{l=k-1}^{1} P_{\U_l^t} - P_{\W_k^t}  ) \right)\|_{\F}\\
				\leq &\sum_{k=1}^{d+m} \|(\bcX^* - \bcX^t) \times_k (P_{\U_k} - P_{\widebar{\U}_k^t} ) \|_{\F}\\
				\leq &  (d+m)\|\bcX^* - \bcX^t\|_{\F} \max_{k=1,\ldots,d+m} \|P_{\U_k} - P_{\widebar{\U}_k^t} \|\\
				\overset{\text{Lemma } \ref{lm: user-friendly-subspace-perturb} } \leq  & 2(d+m)\| \bcX^t - \bcX^* \|_{\F} \max_{k=1,\ldots,d+m} \frac{\|\cM_k(\bcX^*) - \cM_k(\bcX^t)\|}{\sigma_{r^*_k}(\cM_k(\bcX^*))}\\
				\leq & 2(d+m)\| \bcX^t - \bcX^* \|_{\F} \max_{k=1,\ldots,d+m} \frac{\|\cM_k(\bcX^*) - \cM_k(\bcX^t)\|_\F}{\sigma_{r^*_k}(\cM_k(\bcX^*))}\\
				\leq  & \frac{2(d+m) \| \bcX^t - \bcX^* \|^2_{\F}}{\underline{\lambda}}.
			\end{split} 
		\end{equation*}

		In the special matrix setting, i.e., $d + m = 2$, we can get a sharper bound. 
		\begin{equation} \label{eq: decomposition of X-bar complement}
			\begin{split}
				& ~~ P_{(T_{\X^t})_\perp} \X^* = \X^* - P_{T_{\X^t}} \X^*\\
				&  = P_{\U}\X^* + \X^* P_{\V} - P_{\U}\X^*P_{\V} - P_{\U^t} \X^* - \X^*P_{\V^t} + P_{\U^t} \X^* P_{\V^t}\\
				& = (P_{\U} - P_{\U^t})\X^* + \X^*(P_{\V} - P_{\V^t}) - P_{\U}\X^*P_{\V} + P_{\U}\X^*P_{\V^t}-P_{\U}\X^*P_{\V^t} + P_{\U^t} \X^* P_{\V^t}\\
				& = (P_{\U} - P_{\U^t})\X^* (\I - P_{\V^t}) + (\I - P_{\U}) \X^* (P_{\V} - P_{\V^t})\\
				& \overset{(a)}= (P_{\U} - P_{\U^t})\X^* (\I - P_{\V^t})\\
				& \overset{(b)} = (P_{\U} - P_{\U^t})(\X^*- \X^t) (\I - P_{\V^t}),
			\end{split}
		\end{equation} where $\U \in \bbO_{p_1,r^*}, \V \in \bbO_{p_2,r^*}$ are left and right singular vectors of $\X^*$, (a) is due to the fact that $(\I - P_{\U})\X^* (P_{\V} - P_{\V^t}) = 0$ and (b) is because $\X^t(\I - P_{\V^t} ) = 0$. Thus
		\begin{equation*}
			\begin{split}
				\|P_{(T_{\X^t})_\perp} \X^*\|_\F \overset{ \eqref{eq: decomposition of X-bar complement} } \leq \|P_{\U} - P_{\U^t}\| \|\X^*- \X^t\|_\F \overset{\text{Lemma } \ref{lm: user-friendly-subspace-perturb} } \leq \frac{2\|\X^*- \X^t\|_\F^2}{\sigma_{r^*}(\X^*)}. 
			\end{split}
		\end{equation*}

		This finishes the proof of this lemma.
		\quad $\blacksquare$

		\subsection{Proof of Theorem \ref{th: perturbation-HOOI}}
		Recall 
		\begin{equation*}
			\begin{split}
				\widetilde{\U}_i^{1} &= \SVD_{r_i} \left( \cM_{i}( \widetilde{\bcT} \times_{j < i} ( \widetilde{\U}_{j}^{0})^\top \times_{j > i} (\widetilde{\U}_{j}^{0})^\top )  \right)\\
				& = \SVD_{r_i} \left( \cM_i(\bcT \times_{j < i} ( \widetilde{\U}_{j}^{0})^\top \times_{j > i} (\widetilde{\U}_{j}^{0})^\top ) + \cM_i( \bcZ \times_{j < i} ( \widetilde{\U}_{j}^{0})^\top \times_{j > i} (\widetilde{\U}_{j}^{0})^\top )  \right).
			\end{split}
		\end{equation*}
		Suppose $\widetilde{\U}_i^{1} = [\widebar{\U}_i^{1} \quad \widecheck{\U}_i^{1}]$ where $\widebar{\U}_i^{1}$ is composed of the first $r_i^*$ columns of $\widetilde{\U}_i^{1}$ and $\widecheck{\U}_i^{1}$ contains the remaining $(r_i - r_i^*)$ columns of $\widetilde{\U}_i^{1}$. 
		First,
		\begin{equation} \label{ineq: orthogonal-projection-error}
			\begin{split}
				2 \|\bcZ_{\max(\br)}\|_\F &\geq 2\|\left( \cM_i( \bcZ \times_{j < i} ( \widetilde{\U}_{j}^{0})^\top \times_{j > i} (\widetilde{\U}_{j}^{0})^\top )  \right)_{\max(r_i)} \|_\F\\
				&\geq 2\|\left( \cM_i( \bcZ \times_{j < i} ( \widetilde{\U}_{j}^{0})^\top \times_{j > i} (\widetilde{\U}_{j}^{0})^\top )  \right)_{\max(r^*_i)} \|_\F\\
				& \overset{ \text{Lemma } \ref{lm:SVD-projection}  }\geq \| \widebar{\U}_{i\perp}^{1\top} \cM_i(\bcT)  \otimes_{j\neq i} \widetilde{\U}_j^{0}  \|_\F \\
				&  \overset{ \eqref{eq: subspace-decomposition} } \geq \| \widetilde{\U}_{i\perp}^{1\top} \cM_i(\bcT) \otimes_{j\neq i} \U_j (\otimes_{j\neq i} \U_j)^\top \otimes_{j\neq i} \widetilde{\U}_j^{0}  \|_\F \\
				& \geq \|\widetilde{\U}_{i\perp}^{1\top} \cM_i(\bcT)\|_\F \prod_{j\neq i} \sigma_{\min}(\U_j^\top \widetilde{\U}_j^0 )  \\
				& \overset{ (a) }\geq  \|\widetilde{\U}_{i\perp}^{1\top} \cM_i(\bcT)\|_\F (1/2)^{(d-1)/2}\\
				\Longrightarrow &  \|\widetilde{\U}_{i\perp}^{1\top} \cM_i(\bcT)\|_\F \leq 2^{(d+1)/2} \|\bcZ_{\max(\br)}\|_\F.
			\end{split}
		\end{equation} where (a) is because
		\begin{equation} \label{ineq: unequal-size-subspace-bound}
			\begin{split}
				\sigma_{\min}^2( \U_j^\top\widetilde{\U}_{j}^{0} ) = \min_{\a \in \bbR^{r^*_j} }  \frac{\|\a^\top \U_j^\top\widetilde{\U}_{j}^{0} \|_2^2}{\|\a\|_2^2} &= \min_{\a \in \bbR^{r^*_j} }  \frac{\|\a^\top \U_j^\top(\I - \widetilde{\U}_{j \perp}^{0}\widetilde{\U}_{j \perp}^{0\top}) \|_2^2}{\|\a\|_2^2}\\
				& = \min_{\a \in \bbR^{r^*_j} }  \frac{\|\a^\top \U_j^\top\|_2^2 - \|\a^\top \U_j^\top \widetilde{\U}_{j \perp}^{0}\widetilde{\U}_{j \perp}^{0\top} \|_2^2}{\|\a\|_2^2} \\
				&= \frac{\|\a^\top \U_j^\top  \|_2^2}{\|\a\|_2^2} - \max_{\a \in \bbR^{r^*_j} } \frac{\|\a^\top \U_j^\top \widetilde{\U}_{j \perp}^{0}\|_2^2}{\|\a\|_2^2}\\
				& = 1 - \|\U_j^\top \widetilde{\U}_{j \perp}^{0}\|^2 \geq 1/2,
			\end{split}
		\end{equation} 
		
		Finally, by Lemma \ref{lm: tensor estimation via projection}, we have
		\begin{equation*}
			\begin{split}
				\|\widehat{\bcT} - \bcT\|_\F &= \| \widetilde{\bcT} \times_{i=1}^d P_{ \widetilde{\U}^1_i } - \bcT \|_{\F} \\
				&\leq \| \bcZ  \times_{i=1}^d P_{ \widetilde{\U}^1_i } \|_\F + \sum_{i=1}^d \|\widetilde{\U}^{1\top}_{i \perp} \cM_i(\bcT)   \|_\F \overset{ \eqref{ineq: orthogonal-projection-error} }\leq \|\bcZ_{\max(\br)}\|_\F + d \cdot 2^{(d+1)/2} \|\bcZ_{\max(\br)}\|_\F.
			\end{split}
		\end{equation*} 
		This finishes the proof. \quad $\blacksquare$
		
		\subsection{Proof of Lemma \ref{lm:correlated-hermitian}. }
		We divide the proof into two steps: in step 1 we consider the setting $\alpha \leq \sum_{j=1}^w \beta_j $ and in step 2, we consider $\alpha > \sum_{j=1}^w \beta_w$.
		
		{\bfseries Step 1 ($\alpha \leq \sum_{j=1}^w \beta_j $)}. Given $(Y,\X) \sim \mathcal{N}\left(0, \begin{bmatrix}
			1 & \u^\top\\
			\u & \I_w
		\end{bmatrix}\right)$, let $Z = (Y - \sum_{j=1}^w u_w X_w )/\sqrt{1- \sum_{j=1}^w u_j^2 } $. First, by the construction and the property of multivariate Gaussian random vectors, we have $Z$ follows Gaussian distribution with $\bbE[Z] = 0$, $\var(Z) = 1$. Moreover $\text{cov}(Z,X_j) = 0$, so $Z$ is independent of $X_j$ for all $j=1,\ldots,w$ by the property of Gaussian random variables. Given $\bbeta \in \bbR^w$, let us define $\bbeta! = \prod_{i=1}^w \beta_i!$. Then
		\begin{equation*}
			\begin{split}
				\mathbb{E}\left( h_\alpha(Y)\prod_{j=1}^w h_{\beta_j}(X_j) \right) =& \mathbb{E}\left( h_\alpha \left( \sum_{j=1}^w u_w X_w + \sqrt{1- \sum_{j=1}^w u_j^2 } Z \right)\prod_{j=1}^w h_{\beta_j}(X_j) \right)\\
				=& \bbE \left(  \mathbb{E}\left( h_\alpha \left( \sum_{j=1}^w u_w X_w + \sqrt{1- \sum_{j=1}^w u_j^2 } Z \right)\prod_{j=1}^w h_{\beta_j}(X_j)  \Big| Z\right) \right) \\
				\overset{\text{Lemma } \ref{lm: hermite-poly-integration-by-part} } = & \bbE \left( \frac{1}{\sqrt{\bbeta!\alpha!}} \bbE \left( \frac{\partial \sqrt{\alpha!} h_\alpha(  \sum_{j=1}^w u_w X_w + \sqrt{1- \sum_{j=1}^w u_j^2 } Z  )  }{\partial X_1^{\beta_1} \cdots \partial X_w^{\beta_w} }   \Big| Z \right)   \right)\\
				\overset{(a) }= & \bbE \left( \frac{1}{\sqrt{\bbeta!\alpha!}} \bbE \left( \frac{\partial (  \sum_{j=1}^w u_w X_w   + \sqrt{1- \sum_{j=1}^w u_j^2 } Z )^\alpha  }{\partial X_1^{\beta_1} \cdots \partial X_w^{\beta_w} }   \Big| Z \right)   \right)\\
				\overset{(b) } = & \left\{\begin{array}{ll}
					\sqrt{\frac{\alpha!}{\prod_{j=1}^w \beta_j!}}\cdot\prod_{j=1}^w u_j^{\beta_j}, & \text{if } \alpha= \sum_{i=1}^w\beta_j\\
					0, & \text{if }\alpha < \sum_{j=1}^w \beta_j
				\end{array}\right..
			\end{split}
		\end{equation*} Here (a) is because for $\sqrt{\alpha!} h_\alpha(x)$ is an order-$\alpha$ polynomial and the coeffient with $x^\alpha$ is 1 and the fact that except the polynomial $ (  \sum_{j=1}^w u_w X_w  + \sqrt{1- \sum_{j=1}^w u_j^2 } Z )^\alpha$, other terms in $ \sqrt{\alpha!} h_\alpha(  \sum_{j=1}^w u_w X_w + \sqrt{1- \sum_{j=1}^w u_j^2 } Z  )$ will be zero after taking the derivative as $\alpha \leq \sum_{j=1}^w \beta_j $; (b) is because in the expansion of $ (  \sum_{j=1}^w u_w X_w  + \sqrt{1- \sum_{j=1}^w u_j^2 } Z )^\alpha$, the polynomial with respect to $(X_1,\ldots,X_w)$ has degree at most $\alpha$, the fact $\alpha \leq \sum_{j=1}^w \beta_j $ and the coefficient of polynomial $X_1^{\beta_1} \cdots X_q^{\beta_q}$ when $\alpha = \sum_{j=1}^w \beta_j $ is $ \frac{\alpha! \prod_{j=1}^w u_j^{\beta_j} }{\bbeta!} $.
		
		{\bfseries Step 2 ($\alpha > \sum_{j=1}^w \beta_w$)}. Given $(Y,\X) \sim \mathcal{N}\left(0, \begin{bmatrix}
			1 & \u^\top\\
			\u & \I_w
		\end{bmatrix}\right)$, let $Z_j = (X_j - u_j Y)/\sqrt{1-u_j^2} $ for $j = 1,\ldots, w$. Since $(Y,\X)$ is multivariate normal distributed, $Z_j$s are also normal distributed. Moreover, $\bbE(Z_j) = 0$ and $\var(Z_j) = 1$. Finally, $\{Z_j\}_{j=1}^w$ are independent of $Y$ as $\text{cov}(Z_j,Y) = 0$ for $j = 1,\ldots,w$. Then
		\begin{equation*}
			\begin{split}
				\mathbb{E}\left( h_\alpha(Y)\prod_{j=1}^w h_{\beta_j}(X_j) \right) &= \mathbb{E}\left( h_\alpha( Y )\prod_{j=1}^w h_{\beta_j}(u_j Y + \sqrt{1-u_j^2} Z_j ) \right)\\
				& = \bbE \left(  \mathbb{E}\left( h_\alpha( Y )\prod_{j=1}^w h_{\beta_j}(u_j Y + \sqrt{1-u_j^2} Z_j ) \Big| \{Z_j \}_{j=1}^w \right)   \right)\\
				& \overset{ \text{Lemma }\ref{lm: hermite-poly-integration-by-part} } = \bbE \left( \bbE \left(  \frac{\partial \prod_{j=1}^w h_{\beta_j}(u_j Y + \sqrt{1-u_j^2} Z_j ) }{\partial Y^\alpha}  \Big| \{Z_j \}_{j=1}^w \right)  \right)\\
				& \overset{(a)} = 0,
			\end{split}
		\end{equation*} where (a) is because $\prod_{j=1}^w h_{\beta_j}(u_j Y + \sqrt{1-u_j^2} Z_j )$ is a polynomial of degree at most $ \sum_{j=1}^w \beta_j $ in $Y$ and $\alpha > \sum_{j=1}^w \beta_j $. This finishes the proof of this lemma. \quad $\blacksquare$

		\section{Proofs in Section \ref{sec:implementation}} \label{sec: proof-implementation}
		
		\subsection{Properties of Contracted Tensor Inner Product}
		We first introduce the following lemma \ref{lm: contracted-tensor-inner-product}, which reveals a few useful properties of the contracted tensor inner product defined in \eqref{def: contracted-tensor-inner-product} and will be used frequently in deriving efficient implementation of RGN. 
		\begin{Lemma}[Properties of Contracted Tensor Inner Product] \label{lm: contracted-tensor-inner-product}
			Let $\bcX \in \bbR^{p_1 \times \cdots \times p_d \times p_{d+1} \times \cdots \times p_{d+m}}$, $\bcZ \in \bbR^{p_1 \times \cdots \times p_d}$, $\bcW \in \bbR^{p_1 \times \cdots \times p_{k-1} \times q_k \times p_{k+1} \times \cdots \times p_{d}}$ be tensors with $d \geq k \geq 1$, $m \geq 0$. For any $\A \in \bbR^{q_k \times p_k}$, we have
			\begin{equation} \label{eq: contracted-inner-product-inner-change}
				\langle \bcX \times_k \A, \bcW \rangle_* = \langle \bcX, \bcW \times_k \A^\top \rangle_*.
			\end{equation}
			For any $\B \in \bbR^{q_{d+j} \times p_{d+j} }$ with $1 \leq j \leq m$, we have
			\begin{equation} \label{eq: contracted-inner-product-in-out-change}
				\langle \bcX, \bcZ \rangle_* \times_j \B = \langle \bcX \times_{d+j} \B, \bcZ \rangle_*.
			\end{equation}
		\end{Lemma}
		{\noindent Proof of Lemma \ref{lm: contracted-tensor-inner-product}}. 
		We divide the proof into two steps. In step 1, we prove \eqref{eq: contracted-inner-product-inner-change} and in step 2, we prove \eqref{eq: contracted-inner-product-in-out-change}. Throughout the proof, let $\e_j$ be the standard $j$th base vector where $j$th entry is $1$ and others are zero.
		
		{\bfseries Step 1}. For any $i_1 \in [p_{d+1}], \ldots, i_m \in [p_{d+m}]$,
		\begin{equation*}
			\begin{split}
				\langle \bcX , \bcW \times_k \A^\top  \rangle_{*[i_1,\ldots,i_m]} \overset{ \eqref{def: contracted-tensor-inner-product} } =& \langle \bcX_{[:,\ldots,:,i_1,\ldots,i_m]} , \bcW \times_k \A^\top \rangle \\
				=& \langle \bcX \times_{l=d+1}^{d+m} \e_{i_{l-d}} , \bcW \times_k \A^\top \rangle\\
				\overset{ \eqref{eq: matricization relationship} } =& \langle \A \cM_k( \bcX \times_{l=d+1}^{d+m} \e_{i_{l-d}}  ) , \cM_k(\bcW) \rangle\\
				=& \langle \bcX \times_k \A \times_{l=d+1}^{d+m} \e_{i_{l-d}} , \bcW \rangle \\
				=& \langle (\bcX \times_k \A)_{[:,\ldots,:,i_1,\ldots,i_m]} , \bcW \rangle \\
				\overset{ \eqref{def: contracted-tensor-inner-product} }=& \langle \bcX \times_k \A , \bcW \rangle_{[i_1,\ldots,i_m]}. 
			\end{split}
		\end{equation*}
		
		{\bfseries Step 2}. For any $i_1 \in [p_{d+1}],\ldots,i_{j-1} \in [p_{d+j-1}], i_j \in [q_{d+j}], i_{j+1} \in [p_{d+j +1}], \ldots, i_m \in [p_{d+m}]$, 
		\begin{equation*}
			\begin{split}
				&\left( \langle \bcX, \bcZ \rangle_* \times_j \B  \right)_{[ i_1,\ldots,i_m ]} \\
				=&  \langle \bcX, \bcZ \rangle_* \times_{j} \e_{i_j}^\top \B \times_{l \neq j} \e_{i_l}^\top \\
				=& \sum_{z_1 \in [p_{d+1}], \ldots, z_m \in [p_{d+m}] } ( \e_{i_j}^\top \B)_{[z_j]} \left(\prod_{l \neq j} (\e_{i_l})_{[z_l]} \right)  \langle \bcX, \bcZ \rangle_{*[z_1,\ldots,z_m]}\\
				\overset{ \eqref{def: contracted-tensor-inner-product} }=&  \sum_{z_1 \in [p_{d+1}], \ldots, z_m \in [p_{d+m}] } \left(  ( \e_{i_j}^\top \B)_{[z_j]} \prod_{l \neq j} (\e_{i_l})_{[z_l]} \right) \left( \sum_{w_1 \in [p_1], \ldots,w_d \in [p_d]} \bcX_{[w_1,\ldots,w_d, z_1,\ldots,z_m]} \bcZ_{[w_1,\ldots,w_d]}  \right)\\
				=& \sum_{w_1 \in [p_1], \ldots,w_d \in [p_d]}  \bcZ_{[w_1,\ldots,w_d]} \left( \sum_{z_1 \in [p_{d+1}], \ldots, z_m \in [p_{d+m}] } \bcX_{[w_1,\ldots,w_d, z_1,\ldots,z_m]} ( \e_{i_j}^\top \B)_{[z_j]} \prod_{l \neq j} (\e_{i_l})_{[z_l]}   \right)\\
				=& \sum_{w_1 \in [p_1], \ldots,w_d \in [p_d]}  \bcZ_{[w_1,\ldots,w_d]} \cdot \left( \bcX_{[w_1,\ldots,w_d,:,\ldots,:]} \times_{d+j} \e_{i_j}^\top \B \times_{l = d+1, l \neq d+j}^{d+m} \e_{i_{l-d}}   \right)  \\
				=& \langle \bcZ, (\bcX \times_{d+j} \B) \times_{l=d+1 }^{d+m} \e_{i_{l-d}} \rangle= \langle \bcZ, (\bcX \times_{d+j} \B)_{[:,\ldots,:,i_1,\ldots,i_m]} \rangle \\
				\overset{ \eqref{def: contracted-tensor-inner-product} } = &  \langle \bcX \times_{d+j} \B, \bcZ \rangle_{*[:,\ldots,:,i_1,\ldots,i_m]}.
			\end{split} 
		\end{equation*} This finishes the proof of this lemma. \quad $\blacksquare$
		
		\subsection{Proof of Lemma \ref{lm: spectral norm bound of Atop A}}
		First, \eqref{eq: spectrum of LAAL inverse} follows from \eqref{eq: spectrum of LAAL} by the relationship between the spectrum of a linear operator and its inverse, so we just need to show \eqref{eq: spectrum of LAAL}.
		
		The claim \eqref{eq: spectrum of LAAL} is equivalent to say the spectrum of $P_{T_{\bcX^t}} \scA^* \scA P_{T_{\bcX^t}}$ is lower and upper bounded by $1 - R_{2\br}$ and $1+R_{2\br}$, respectively, for $\bcZ \in T_{\bcX^t} \bbM_\br$. Since $P_{T_{\bcX^t}} \scA^* \scA P_{T_{\bcX^t}}$ is a symmetric operator, its spectrum can be upper bounded by $\sup_{\bcZ \in T_{\bcX^t} \bbM_\br:\|\bcZ\|_{\F} =1 } \langle \bcZ, P_{T_{\bcX^t}} \scA^* \scA P_{T_{\bcX^t}}(\bcZ) \rangle $ and lower bounded by $\inf_{\bcZ \in T_{\bcX^t} \bbM_\br:\|\bcZ\|_{\F} =1 } \langle \bcZ, P_{T_{\bcX^t}} \scA^* \scA P_{T_{\bcX^t}} (\bcZ) \rangle $. Also
		\begin{equation*}
			\begin{split}
				&\sup_{\bcZ \in T_{\bcX^t} \bbM_\br:\|\bcZ\|_{\F} =1 } \langle \bcZ, P_{T_{\bcX^t}} \scA^* \scA P_{T_{\bcX^t}}(\bcZ) \rangle= \sup_{\bcZ \in T_{\bcX^t} \bbM_\br: \|\bcZ\|_{\F} =1 } \|\scA P_{T_{\bcX^t}} (\bcZ)\|^2_{\F} \overset{(a)} \leq 1+R_{2\br}\\
				& \inf_{\bcZ \in T_{\bcX^t} \bbM_\br:\|\bcZ\|_{\F} =1 } \langle \bcZ, P_{T_{\bcX^t}} \scA^* \scA P_{T_{\bcX^t}}(\bcZ) \rangle = \inf_{\bcZ \in T_{\bcX^t} \bbM_\br: \|\bcZ\|_{\F} = 1 } \|\scA P_{T_{\bcX^t}} (\bcZ)\|^2_{\F} \overset{(a)} \geq 1-R_{2\br}.
			\end{split}
		\end{equation*} Here (a) is by the TRIP condition for $\scA$, $P_{T_{\bcX^t}}(\bcZ)$ is at most Tucker rank $2\br$ by Lemma \ref{lm: tangent vector rank 2r property} and the fact $\|P_{T_{\bcX^t}} (\bcZ)\|_\F \leq \|\bcZ\|_{\F}$.
		\quad $\blacksquare$
		
		\subsection{Proof of Proposition \ref{prop: efficient-implementation-RGN} }
		First, for convenience of presentation, given $\U_k \in \bbO_{p_k,r_k}$, we define $\U_k^{(1)} = \U_k, \U_k^{(2)} = \U_{k\perp}$. By orthogonality, we know for any tensor $\bcT$ of dimension $p_1\times \cdots \times p_{d+m}$, 
		\begin{equation}\label{eq: tensor-norm-decomp}
			\begin{split}
				\|\bcT\|_\F^2 = & \|\bcT\times_{k=1}^{d+m} (P_{\U_k} + P_{\U_{k\perp}})\|_\F^2 = \left\|\sum_{s_1,\ldots, s_{d+m}=1}^{2} \bcT\times_{k=1}^{d+m} P_{\U_k^{(s_k)}}\right\|_\F^2\\ 
				= & \sum_{s_1,\ldots, s_{d+m}=1}^{2}\left\|\bcT\times_{k=1}^{d+m} P_{\U_k^{(s_k)}}\right\|_\F^2 = \sum_{s_1,\ldots, s_{d+m}=1}^{2}\left\|\bcT\times_{k=1}^{d+m} \U_k^{(s_k)\top }\right\|_\F^2.
			\end{split}
		\end{equation}
		
		By the parameterization form of the tangent space of $\bcX^t$ given in \eqref{eq: tangent-space-representation}, we know to solve the RGN update, i.e., $\bcX^{t+0.5} = \argmin_{\bcX \in T_{\bcX^t} \bbM_{\br} } \frac{1}{2}\|\bcY - \scA P_{T_{\bcX^t}}(\bcX)\|_\F^2$, it is equivalent to solve
		\begin{equation} \label{eq: reformulated-least-squares}
			\begin{split}
				&(\bcB^t,\{\D_k^{t}\}_{k=1}^{d+m})\\
				= &\argmin_{\substack{\bcB\in \bbR^{r_1 \times \cdots \times r_{d+m}}, \\ \D_k\in \bbR^{(p_k -r_k) \times r_k}, k=1,\ldots,d+m}} \sum_{i=1}^n \left\| \bcY_i - \langle \bcA_i, \bcB \times_{k=1}^{d+m} \U_k^t \rangle_* - \sum_{k=1}^{d+m} \langle \bcA_i, \bcS^t \times_k \U_{k\perp}^t \D_k \times_{j \neq k} \U_j^t \rangle_* \right\|_\F^2.
			\end{split}
		\end{equation}
		
		Next, let us decompose each term in the summand on the right-hand side of \eqref{eq: reformulated-least-squares}.
		\begin{equation} \label{eq: tensor-norm-decom-1}
			\begin{split}
				&\left\| \bcY_i - \langle \bcA_i, \bcB \times_{k=1}^{d+m} \U_k^t \rangle_* - \sum_{k=1}^{d+m} \langle \bcA_i, \bcS^t \times_k \U_{k\perp}^t \D_k \times_{j \neq k} \U_j^t \rangle_* \right\|_\F^2 \\
				\overset{ \eqref{eq: tensor-norm-decomp} }= & \sum_{s_{d+1},\ldots, s_{d+m}=1}^{2} \Big\| \bcY_i \times_{l=1}^m \U_{l+d}^{t(s_{l+d})\top} - \langle \bcA_i, \bcB \times_{k=1}^{d+m} \U_k^t \rangle_* \times_{l=1}^m \U_{l+d}^{t(s_{l+d})\top} \\
				& \quad - \sum_{k=1}^{d+m} \langle \bcA_i, \bcS^t \times_k \U_{k\perp}^t \D_k \times_{j \neq k} \U_j^t \rangle_* \times_{l=1}^m \U_{l+d}^{t(s_{l+d})\top}   \Big\|_\F^2 \\
				\overset{\text{Lemma } \ref{lm: contracted-tensor-inner-product} } =& \sum_{s_{d+1},\ldots, s_{d+m}=1}^{2} \Big\| \bcY_i \times_{l=1}^m \U_{l+d}^{t(s_{l+d})\top} - \langle \bcA_i, (\bcB \times_{k=1}^{d+m} \U_k^t) \times_{l=d+1}^{d+m} \U_{l}^{t(s_l)\top} \rangle_* \\
				& \quad - \sum_{k=1}^{d+m} \langle \bcA_i, (\bcS^t \times_k \U_{k\perp}^t \D_k \times_{j \neq k} \U_j^t) \times_{l=d+1}^{d+m} \U_{l}^{t(s_l)\top}  \rangle_*    \Big\|_\F^2 \\
				& =  \sum_{s_{d+1},\ldots, s_{d+m}=1}^{2} \Big\| \bcY_i \times_{l=1}^m \U_{l+d}^{t(s_{l+d})\top} - \underbrace{\langle \bcA_i, \bcB \times_{k=1}^{d} \U_k^t \times_{l=d+1}^{d+m} \U_{l}^{t(s_l)\top} \U_l^t \rangle_*}_{\text{(I)}} \\
				& \quad - \underbrace{\sum_{k=1}^{d} \langle \bcA_i, \bcS^t \times_k \U_{k\perp}^t \D_k \times_{j=1,j \neq k}^d \U_j^t \times_{l=d+1}^{d+m} \U_{l}^{t(s_l)\top} \U_l^t  \rangle_*}_{\text{(II)}}     \\
				& \quad - \underbrace{\sum_{k=d+1}^{d+m} \langle \bcA_i, \bcS^t \times_{j =1}^d \U_j^t \times_k \U_k^{t(s_k)\top}\U_{k\perp}^t \D_k \times_{l=d+1,l \neq k}^{d+m} \U_{l}^{t(s_l)\top} \U_l^t  \rangle_*}_{\text{(III)}}    \Big\|_\F^2.
			\end{split}
		\end{equation}
		Notice that (I) and (II) are non-zero only if $s_l = 1$ for all $l = d+1,\ldots,d+m$; (III) is non-zero only when there is one $z \in \{d+1,\ldots,d+m \}$ such that $s_z = 2$ and for all other $l \neq z$, $s_l = 1$; moreover, suppose $s_z = 2$ for some $z \in \{d+1,\ldots,d+m \}$, then
		\begin{equation*}
			\text{(III)} = \langle \bcA_i , \bcS^t \times_{j=1}^d \U_j^t \times_z \D_{z} \rangle_*. 
		\end{equation*} With the above observation, \eqref{eq: tensor-norm-decom-1} can be simplified as 
		\begin{equation*}
			\begin{split}
				&\left\| \bcY_i - \langle \bcA_i, \bcB \times_{k=1}^{d+m} \U_k^t \rangle_* - \sum_{k=1}^{d+m} \langle \bcA_i, \bcS^t \times_k \U_{k\perp}^t \D_k \times_{j \neq k} \U_j^t \rangle_* \right\|_\F^2 \\
				=& \left\| \bcY_i \times_{l=1}^m \U_{l+d}^{t\top} - \langle \bcA_i, \bcB \times_{k=1}^{d} \U_k^t \rangle_* -\sum_{k=1}^{d} \langle \bcA_i, \bcS^t \times_k \U_{k\perp}^t \D_k \times_{j=1,j \neq k}^d \U_j^t  \rangle_* \right\|_\F^2 \\
				& + \sum_{z=1}^m \left\|  \bcY_i \times_z \U_{z+d \perp}^{t\top} \times_{l\neq z} \U_{l+d}^{t\top} -  \langle \bcA_i , \bcS^t \times_{j=1}^d \U_j^t \times_{z+d} \D_{z+d} \rangle_* \right\|_\F^2 \\
				\overset{ \text{Lemma }\ref{lm: contracted-tensor-inner-product} }	= & \underbrace{\left\| \bcY_i \times_{l=1}^m \U_{l+d}^{t\top} - \langle \bcA_i \times_{k=1}^d \U_k^{t\top} , \bcB \rangle_* -\sum_{k=1}^{d} \langle \bcA_i \times_k \U_{k\perp}^{t\top} \times_{j=1,j \neq k}^d \U_j^{t\top} , \bcS^t \times_k  \D_k  \rangle_* \right\|_\F^2}_{(T_0)} \\
				& + \sum_{z=1}^m \underbrace{\left\|  \bcY_i \times_z \U_{z+d \perp}^{t\top} \times_{l\neq z} \U_{l+d}^{t\top} -  \langle \bcA_i  \times_{j=1}^d \U_j^{t\top}, \bcS^t \times_{z+d} \D_{z+d} \rangle_* \right\|_\F^2}_{(T_z)}
			\end{split}
		\end{equation*}
		We note that $(T_0)$ only involves $(\bcB,\{ \D_k\}_{k=1}^d)$ and $(T_z)$ only involves $\D_{z+d}$ for $z = 1,\ldots,m$. Moreover, we can view $(T_0)$ and $\{(T_z)\}_{z=1}^m$ as $(m+1)$ separate least squares for the reasons below. First,
		\begin{equation*}
			\begin{split}
				(T_0) &= \sum_{j_l \in [r_{d+l}], l = 1,\ldots,m } \Big(  \left(\bcY_i \times_{l=1}^m \U_{l+d}^{t\top} \right)_{[j_1,\ldots,j_m]}  -  \langle \bcA_i \times_{j=1}^d \U_j^{t\top} , \bcB \rangle_{*[j_1,\ldots,j_m]} - \\
				& \quad  - \sum_{k=1}^{d} \langle \bcA_i \times_k \U_{k\perp}^{t\top} \times_{j=1,j \neq k}^d \U_j^{t\top} , \bcS^t \times_k  \D_k  \rangle_{*[j_1,\ldots,j_m]}  \Big)^2\\
				& \overset{ \text{Lemma }\ref{lm: contracted-tensor-inner-product}, \eqref{eq: matricization relationship} } = \sum_{j_l \in [r_{d+l}], l = 1,\ldots,m } \Big(  \left(\bcY_i \times_{l=1}^m \U_{l+d}^{t\top} \right)_{[j_1,\ldots,j_m]}  -  \langle \bcA_i \times_{j=1}^d \U_j^{t\top} , \bcB_{[:,\ldots,:,j_1,\ldots,j_m]} \rangle - \\
				& \quad  - \sum_{k=1}^d \langle \U_{k\perp}^{t\top} \cM_k\left( \bcA_i \times_{j\neq k} \U_j^{t\top} \right) \left( \cM_k( \bcS^t_{[:,\ldots,:,j_1,\ldots,j_m]} ) \right)^\top  , \D_k \rangle  \Big)^2.
			\end{split}
		\end{equation*} Moreover,
		\begin{equation*}
			\begin{split}
				(T_z) \overset{ \text{Lemma }\ref{lm: contracted-tensor-inner-product} } =& \left\| \bcY_i \times_z \U_{z+d \perp}^{t\top} \times_{l\neq z} \U_{l+d}^{t\top} -  \langle \bcA_i  \times_{j=1}^d \U_j^{t\top}, \bcS^t \rangle_* \times_z  \D_{z+d} \right\|_\F^2\\
				\overset{ \eqref{eq: matricization relationship} } =& \left\| \U_{z+d \perp}^{t\top} \cM_z(\bcY_i \times_{l \neq z} \U_{l+d}^{t\top} ) - \D_{z+d} \cM_z(  \langle \bcA_i  \times_{j=1}^d \U_j^{t\top}, \bcS^t  \rangle_* ) \right\|_\F^2,
			\end{split}
		\end{equation*} notice that the above formula can be viewed as a multivariate least squares. Finally, the proposition follows by considering the sum of $(T_0)$, $\{ (T_z) \}_{z=1}^m$ over sample indices $ i$ in \eqref{eq: reformulated-least-squares}. \quad $\blacksquare$

		\subsection{Proof of Lemma \ref{lm: RGN-tensor-on-vector}}
		By computing the gradient of the objective in \eqref{eq: alg least square} and setting it to zero,  $\bcX^{t+0.5}$ can be uniquely solved by the following linear system:
		\begin{equation}\label{eq: RGN-linear-system-interpre}
			P_{T_{\bcX^t}} \scA^* \scA P_{T_{\bcX^t}} (\bcX^{t+0.5}) = P_{T_{\bcX^t}}( \scA^*(\bcY) ).
		\end{equation}
		In this setting, it reduces to
		\begin{equation} \label{eq: RGN-linear-sys-tensor-on-vector}
			P_{T_{\bcX^t}}\left( \bcX^{t+0.5} \times_1 \A^\top \A  \right) = P_{T_{\bcX^t}}( \bcY \times_1 \A^\top ).
		\end{equation} Since $\bcX^{t+0.5}\in T_{\bcX^t} \bbM_{\br}$, we known from \eqref{eq: tangent-space-representation} that $\bcX^{t+0.5}$ can be written as 
		\begin{equation} \label{eq: Xt0.5-tangent-representation}
			\bcX^{t+0.5} = \bcB^t \times_{k=1}^{1+m} \U_k^t + \sum_{k=1}^{1+m} \bcS^t \times_{k} \U_{k\perp}^t \D^t_k \times_{j \neq k} \U_j^t,
		\end{equation} for some to be determined $(\bcB^t, \{ \D_k^t \}_{k=1}^{1+m} )$.
		\begin{itemize}[leftmargin=*]
			\item Let us first compute $\cM_k(\bcS^t \times_{k} \U_{k\perp}^t \D^t_k \times_{j \neq k} \U_j^t)$ for $k = 2,\ldots,1+m$. By the projection operation onto the tangent space $T_{\bcX^t} \bbM_{\br}$ given in \eqref{eq: tangent space projector}, we have for $k = 2,\ldots,1+m$,
			\begin{equation} \label{eq: mode-2-m+1-projection}
				\begin{split}
					\U_{k\perp}^{t\top} \cM_k\left(  \bcX^{t+0.5} \times_1 \A^\top \A \right) \W_k^t = \U_{k\perp}^{t\top} \cM_{k}(\bcY \times_1 \A^\top) \W_k^t.
				\end{split}
			\end{equation} By plugging \eqref{eq: Xt0.5-tangent-representation} into \eqref{eq: mode-2-m+1-projection}, we get
			\begin{equation} \label{eq: mode-k-projection-eq}
				\begin{split}
					&\U_{k\perp}^{t\top} \cM_k\left(  \bcX^{t+0.5} \times_1 \A^\top \A \right) \W_k^t\\
					=&\U_{k\perp}^{t\top} \cM_k( \bcS^t \times_1 \A^\top \A \U_1^t \times_{k} \U_{k\perp}^t \D_k^t \times_{i \neq k, i \neq 1} \U_j^t ) \W_k^t \\
					=& \U_{k\perp}^{t\top}\U_{k\perp}^t \D_k^t \cM_k(\bcS^t) \left( \otimes_{i \neq 1, i \neq k} \U_j^{t\top} \otimes (\U_1^{t\top} \A^\top \A)  \right) \W_k^t  \\
					= & \D_k^t \cM_k(\bcS^t) \left( \otimes_{\substack{i \neq k,\\i\neq 1}} \I_{r_i} \otimes(\U_1^{t\top} \A^\top \A \U_1^t ) \right) \V_k^t \\
					\overset{(a)}= & \D_k^t \cM_k(\bcS^t) \V_k^t \V_k^{t\top} \left( \otimes_{\substack{i \neq k,\\i\neq 1}} \I_{r_i} \otimes(\U_1^{t\top} \A^\top \A \U_1^t ) \right) \V_k^t \\
					=& \U_{k\perp}^{t\top} \cM_{k}(\bcY \times_1 \A^\top) \W_k^t.
				\end{split}
			\end{equation} Here (a) is because $\V_k^t$ by definition spans the row space of $\cM_k(\bcS^t)$. The last equality in \eqref{eq: mode-k-projection-eq} implies
			\begin{equation} \label{eq: DktMkS}
				\D_k^t \cM_k(\bcS^t) \V_k^t  = \U_{k\perp}^{t\top} \cM_{k}(\bcY \times_1 \A^\top) \W_k^t \left( \V_k^{t\top} \left( \otimes_{\substack{i \neq k,\\i\neq 1}} \I_{r_i} \otimes(\U_1^{t\top} \A^\top \A \U_1^t ) \right) \V_k^t \right)^{-1}.
			\end{equation} So
			\begin{equation*}
				\begin{split}
					&\cM_k(\bcS^t \times_{k} \U_{k\perp}^t \D^t_k \times_{j \neq k} \U_j^t) \\
					=& \U_{k\perp}^t \D^t_k \cM_k(\bcS^t) \V_k^t \V_k^\top \otimes_{j \neq k} \U_j^{t\top} \\
					\overset{\eqref{eq: DktMkS} } =& \U_{k\perp}^t\U_{k\perp}^{t\top} \cM_{k}(\bcY \times_1 \A^\top) \W_k^t \left( \V_k^{t\top} \left( \otimes_{\substack{i \neq k,\\i\neq 1}} \I_{r_i} \otimes(\U_1^{t\top} \A^\top \A \U_1^t ) \right) \V_k^t \right)^{-1} \W_k^{t\top}.  
				\end{split}
			\end{equation*}

			\item Now we compute $\cM_1(\bcB^t)$ and $\cM_1(\bcS^t \times_{1} \U_{1\perp}^t \D^t_1 \times_{j \neq 1} \U_j^t)$. Similar to \eqref{eq: mode-2-m+1-projection}, here by projecting both sides of \eqref{eq: RGN-linear-sys-tensor-on-vector} onto the core and the first mode, we have
			\begin{equation} \label{eq: mode-1-projection}
				\begin{split}
					\U_{1\perp}^{t\top} \cM_1\left(  \bcX^{t+0.5} \times_1 \A^\top \A \right) \W_1^t &= \U_{1\perp}^{t\top} \cM_{1}(\bcY \times_1 \A^\top) \W_1^t, \\
					\left(  \bcX^{t+0.5} \times_1 \A^\top \A \right) \times_{k=1}^{1+m} \U_k^{t\top} & = (\bcY \times_1 \A^\top ) \times_{k=1}^{1+m} \U_k^{t\top}.
				\end{split}
			\end{equation} By plugging \eqref{eq: Xt0.5-tangent-representation} into \eqref{eq: mode-1-projection}, we get
			\begin{equation} \label{eq: mode-1-part}
				\U_{1\perp}^{t\top} \A^\top \A \U_1^t \cM_1(\bcB^t) \V_1^t + \U_{1\perp}^{t\top} \A^\top \A  \U_{1\perp}^{t} \D_1^t \cM_1(\bcS^t) \V_1^t = \U_{1\perp}^{t\top} \cM_{1}(\bcY \times_1 \A^\top) \W_1^t,
			\end{equation} and 
			\begin{equation} \label{eq: core-part}
				\bcB^t \times_1 \U_1^{t\top} \A^\top \A \U_1^t  + \bcS^t \times_1 \U_1^{t\top}\A^\top \A  \U_{1\perp}^{t} \D_1^t = (\bcY \times_1 \A^\top ) \times_{k=1}^{1+m} \U_k^{t\top}.
			\end{equation}
			
			Do matricization along mode-$1$ on both sides of \eqref{eq: core-part} yields
			\begin{equation} \label{eq: core-part2}
				\U_1^{t\top}  \A^\top \A \U_1^t \cM_1(\bcB^t) + \U_1^{t\top} \A^\top \A  \U_{1\perp}^{t} \D_1^t \cM_1(\bcS^t) = \U_{1}^{t\top} \cM_{1}(\bcY \times_1 \A^\top) \otimes_{j=(1+m)}^2 \U_j^t.
			\end{equation} By multiplying $\V_1^t$ on both sides of \eqref{eq: core-part2} and combining with \eqref{eq: mode-1-part} yields
			\begin{equation} \label{eq: mode-1-part2}
				\begin{split}
					& \A^\top \A \U_1^t \cM_1(\bcB^t) \V_1^t + \A^\top \A  \U_{1\perp}^{t} \D_1^t \cM_1(\bcS^t) \V_1^t = \cM_{1}(\bcY \times_1 \A^\top) \W_1^t \\
					\Longrightarrow & \U_1^t \cM_1(\bcB^t) \V_1^t + \U_{1\perp}^{t} \D_1^t \cM_1(\bcS^t) \V_1^t = (\A^\top \A)^{-1} \cM_{1}(\bcY \times_1 \A^\top) \W_1^t\\
					\overset{(a)}\Longrightarrow & \D_1^t \cM_1(\bcS^t) \V_1^t = \U_{1\perp}^{t\top} (\A^\top \A)^{-1} \cM_{1}(\bcY \times_1 \A^\top) \W_1^t,
				\end{split}
			\end{equation} here (a) is by multiplying $\U_{1\perp}^{t\top}$ on both sides of the equation. So
			\begin{equation*}
				\begin{split}
					&\quad \cM_1(\bcS^t \times_{1} \U_{1\perp}^t \D^t_1 \times_{j \neq 1} \U_j^t)\\ &= \U_{1\perp}^t \D^t_1 \cM_1(\bcS^t) \otimes_{j=(1+m)}^2 \U_j^{t\top}\\
					& = \U_{1\perp}^t \D^t_1 \cM_1(\bcS^t) \V_1^t \V_1^{t\top} \otimes_{j=(1+m)}^2 \U_j^{t\top}\\
					& \overset{ \eqref{eq: mode-1-part2} } =  \U_{1\perp}^t \U_{1\perp}^{t\top} (\A^\top \A)^{-1} \cM_{1}(\bcY \times_1 \A^\top) \W_1^t \W_1^{t\top}.
				\end{split}
			\end{equation*} Finally by \eqref{eq: core-part2}, we have
			\begin{equation*}
				\begin{split}
					& \cM_1(\bcB^t)  \\
					= & (\U_{1}^{t\top} \A^\top \A \U_1^t)^{-1} \left( \U_{1}^{t\top} \cM_{1}(\bcY \times_1 \A^\top) \otimes_{j=(1+m)}^2 \U_j^t- \U_1^{t\top} \A^\top \A  \U_{1\perp}^{t} \D_1^t \cM_1(\bcS^t) \right)\\
					\overset{ \eqref{eq: mode-1-part2} }= & (\U_{1}^{t\top} \A^\top \A \U_1^t)^{-1} \\
					& \cdot \left( \U_{1}^{t\top}\A^\top \cM_{1}(\bcY) \otimes_{j=(1+m)}^2 \U_j^t- \U_1^{t\top} \A^\top \A  \U_{1\perp}^{t} \U_{1\perp}^{t\top} (\A^\top \A)^{-1} \cM_{1}(\bcY \times_1 \A^\top) \W_1^t \V_1^{t\top} \right)\\
					= & ( \U_1^{t\top} \A^\top \A \U_1^t )^{-1} \U_1^{t\top} \A^\top\\
					& \cdot \left( \cM_1(\bcY) \otimes_{j=(1+m)}^2 \U_j^t -  \A \U_{1\perp}^t \U_{1\perp}^{t\top} (\A^\top \A)^{-1} \A^\top \cM_1(\bcY) \W^t_1 \V^{t\top}_1  \right).
				\end{split}
			\end{equation*} 
		\end{itemize} This finishes the proof. \quad $\blacksquare$

		\section{Additional Proofs and Lemmas} \label{sec: additional lemmas}
		
		\begin{Lemma}(Tensor Restricted Orthogonal Property)\label{lm:retricted orthogonal property}
			Let $\bcZ_1, \bcZ_2 \in \bbR^{p_1 \times \cdots \times p_d}$ be two low Tucker rank tensors with $\tuckerrank(\bcZ_1) = \br_1 := (r_1,\ldots,r_d)$, $ \tuckerrank(\bcZ_2) = \br_2 := (r_1',\ldots,r_d')$. Then, 
			\begin{equation} \label{ineq: restricted-orthogonal-property}
				|\langle \scA(\bcZ_1), \scA(\bcZ_2) \rangle - \langle \bcZ_1 , \bcZ_2 \rangle | \leq R_{\br_1 + \br_2} \|\bcZ_1\|_{\F} \|\bcZ_2\|_{\F}.
			\end{equation}
		\end{Lemma}
		{\noindent \bf Proof.}  Without loss of generality, assume $\|\bcZ_1\|_{\F} = 1$, $\|\bcZ_2\|_{\F} = 1$. Notice that $\bcZ_1 + \bcZ_2$ is of at most Tucker rank $\br_1 + \br_2$ as the matricization of $\bcZ_1 + \bcZ_2$ on each mode $k$ is of at most rank $r_k + r_k'$. Similarly, $\bcZ_1 - \bcZ_2$ is also at most Tucker rank $\br_1 + \br_2$. Then by the definition of TRIP constant, we have
		\begin{equation*}
			\begin{split}
				&(1-R_{\br_1 + \br_2}) \|\bcZ_1 \pm \bcZ_2 \|^2_{\F} \leq \| \scA(\bcZ_1 \pm \bcZ_2) \|_{\F}^2\\
				&(1+R_{\br_1 + \br_2}) \|\bcZ_1 \pm \bcZ_2 \|^2_{\F} \geq \| \scA(\bcZ_1 \pm \bcZ_2) \|_{\F}^2.
			\end{split}
		\end{equation*}Notice that the above inequalities hold in both settings when $R_{\br_1 + \br_2} > 1$ and $R_{\br_1 + \br_2} < 1$. Then
		\begin{equation*}
			\begin{split}
				\|\scA(\bcZ_1 + \bcZ_2)\|_\F^2 -  \|\scA(\bcZ_1 - \bcZ_2)\|_\F^2 &= 4 \langle \scA(\bcZ_1), \scA(\bcZ_2) \rangle \leq 4 \langle \bcZ_1, \bcZ_2 \rangle + 4 R_{\br_1 + \br_2}, \\
				\|\scA(\bcZ_1 + \bcZ_2)\|_\F^2 -  \|\scA(\bcZ_1 - \bcZ_2)\|_\F^2 &= 4 \langle \scA(\bcZ_1), \scA(\bcZ_2) \rangle \geq 4 \langle \bcZ_1, \bcZ_2 \rangle - 4 R_{\br_1 + \br_2}
			\end{split}
		\end{equation*} and this implies \eqref{ineq: restricted-orthogonal-property}. \quad $\blacksquare$

		\begin{Lemma}(Tensor Estimation from Projection \cite[Lemma 11]{luo2021low}) \label{lm: tensor estimation via projection}
			Given two order-$d$ tensors $\bcY, \bcX \in \bbR^{p_1 \times \cdots \times p_d}$. Suppose $\U_k^0 \in \bbO_{p_k,r_k}$, then
			\begin{equation*}
				\left\|\bcY \times_{k=1}^d P_{\U_k^0} - \bcX   \right\|_{\F} \leq \left\| (\bcY-\bcX) \times_{1} P_{\U^0_1} \times \cdots \times_{d} P_{\U^0_d}\right\|_{\F} + \sum_{k=1}^d \left\|\U_{k\perp}^{0\top} \mathcal{M}_{k}(\bcX)\right\|_\F.
			\end{equation*} 
		\end{Lemma}
		
		\begin{Definition} Suppose $\scT$ is a subset of a normed space. A set $\cN^{\scT}_\epsilon \subseteq \scT$ is called a $\epsilon$-net of $\scT$ with respect to the norm $\|\cdot\|$ if for each $\bcT \in \scT$, there exists $\bcT_0 \in \cN_\epsilon^\scT$ with $\| \bcT - \bcT_0 \| \leq \epsilon$. The minimal cardinality of an $\epsilon$-net of $\scT$ with respect to the norm $\|\cdot\|$ is denoted by $\cN(\scT, \|\cdot\|, \epsilon )$ and is called the covering number of $\scT$ at scale $\epsilon$.
		\end{Definition}
		Next, we give a sharp bound on the covering number for the set of low Tucker rank tensors with respect to the Frobenius norm and it improves upon \cite[Lemma 2]{rauhut2017low}.
		
		\begin{Lemma}[Covering Number for Low Tucker Rank Tensors] \label{lm: covering-number} Let $0< \epsilon < 1$, the covering number of $\scT_{\bp,\br}:= \{ \bcT \in \bbR^{p_1 \times \cdots \times p_d}: \tuckerrank(\bcT) \leq \br:= (r_1,\ldots, r_d), \|\bcT\|_\F \leq 1 \}$ with respect to the Frobenious norm satisfies $\cN( \scT_{\bp,\br}, \|\cdot\|_\F, \epsilon ) \leq \left(  \frac{3(d+1)}{\epsilon} \right)^{ \prod_{i=1}^d r_i } \prod_{i=1}^d \left(  \frac{c_0(d+1)}{\epsilon} \right)^{ (p_i - r_i) r_i } $ for some absolute constant $c_0 > 0$.
		\end{Lemma}
		{\noindent \bf Proof of Lemma \ref{lm: covering-number}.} Given any $\bcT \in \scT_{\bp, \br}$, $\bcT$ has the Tucker decomposition $\bcT = \bcS \times_{i=1}^d \U_i$ for some $\bcS$ satisfying $\|\bcS\|_\F \leq 1$ and $\U_i \in \bbO_{p_i,r_i}$. We will construct a $\epsilon$-net for $\scT_{\bp,\br}$ by covering $\{\U_i \}_{i=1}^d$ and the set of $r_1 \times \cdots \times r_d$ tensors with Frobenius norm at most $1$.
		
		By \cite[Lemma 7]{zhang2018tensor}, we can construct a $\epsilon/(d+1)$-net ($0<\epsilon/(d+1) < 1$) $\{ \bcS^{(1)},\cdots, \bcS^{(N_{\bcS})} \}$ for $\{ \bcS' \in \bbR^{r_1 \times \cdots \times r_{d}}: \|\bcS'\|_\F \leq 1 \}$ such that $$\sup_{\bcS': \|\bcS'\|_\F \leq 1 } \min_{i \leq N_{\bcS}} \|\bcS' - \bcS^{(i)}\|_\F \leq \epsilon/(d+1) $$
		with $N_{\bcS} \leq (\frac{3(d+1)}{\epsilon})^{ \prod_{i=1}^{d} r_i } $.
		
		At the same time, by \cite[Proposition 8]{szarek1982nets}, for each $k = 1,\ldots, d$, we can construct a $\epsilon/(d+1)$-net $\{ \U_k^{(1)},\ldots, \U_k^{(N_k)} \}$ on the Grassmann manifold of $r_k$-dimensional subspaces in $\bbR^{p_k}$ with the metric $d(\U_1, \U_2) = \| \U_1 \U_1^\top - \U_2 \U_2^\top \|_\F$ such that 
		\begin{equation*}
			\sup_{\U_k \in \bbO_{p_k, r_k}} \min_{i \leq N_k} d( \U_k, \U_k^{(i)} ) \leq \epsilon/(d+1)
		\end{equation*} with $N_k \leq (\frac{c_0(d+1)}{\epsilon})^{r_k(p_k - r_k)} $ for some absolute constant $c_0 > 0$.
		
		Given any fixed $\bcT \in \scT_{\bp,\br}$ with Tucker decomposition $\bcS \times_{i=1}^d \U_i$, we can find $\U_k^{(i_k)} $ in the corresponding $\epsilon/(d+1)$-net such that $d(\U_k^{(i_k)}, \U_k) \leq \epsilon/(d+1)$. Let $\O_k = \argmin_{\O \in \bbO_{r_k}} \| \U_k \O - \U_k^{(i_k)} \|_\F$. By \cite[Lemma 1]{cai2018rate}, we have $\| \U_k \O_k - \U_k^{(i_k)} \|_\F \leq  d(\U_k^{(i_k)}, \U_k) \leq \epsilon/(d+1)$. Denote $\widebar{\bcS} = \bcS \times_{k=1}^{d} \O_k^\top$ and let $\bcS^{(i_0)}$ be the one in the core tensor $\epsilon/(d+1)$-net such that $\|\bcS^{(i_0)} - \widebar{\bcS}\|_\F \leq \epsilon/(d+1)$. Thus
		\begin{equation*} 
			\begin{split}
				& \| \bcT - \bcS^{(i_0)} \times_{k=1}^{d} \U^{(i_k)}_k  \|_\F \\
				=&\left\| \widebar{\bcS} \times_{k=1}^{d} \U_k \O_k- \bcS^{(i_0)} \times_{k=1}^{d} \U^{(i_k)}_k \right \|_\F \\
				=& \left\| (\widebar{\bcS} - \bcS^{(i_0)}) \times_{k=1}^{d} \U_k \O_k + \sum_{k=1}^{d} \bcS^{(i_0)} \times_{j < k}  \U_j^{(i_j)}  \times_k ( \U_k \O_k- \U^{(i_k)}_k) \times_{j > k} \U_j \O_j \right \|_\F\\
				\leq & \|(\widebar{\bcS} - \bcS^{(i_0)}) \times_{k=1}^{d} \U_k \O_k\|_\F + \sum_{k=1}^d \|\bcS^{(i_0)} \times_{j < k}  \U_j^{(i_j)}  \times_k ( \U_k \O_k- \U^{(i_k)}_k) \times_{j > k} \U_j \O_j\|_\F\\
				\leq & (d+1) \epsilon/(d+1) = \epsilon.
			\end{split}
		\end{equation*} This finishes the proof of this lemma. \quad $\blacksquare$
		
		The following Lemma \ref{lm:SVD-projection} quantifies the projection error under the perturbation model.
		\begin{Lemma}(A perturbation projection error bound \citep[Theorem 2]{luo2021schatten})\label{lm:SVD-projection}
			Suppose $\B = \A + \Z$ for some rank-$r$ matrix $\A$ and perturbation matrix $\Z$. Denote the top rank $r$ truncated SVD of $\B$ as $\widehat{\U} \widehat{\bSigma} \widehat{\V}^\top$. Then for any $q \in [1, \infty]$,
			\begin{equation*}
				\max\left\{\|P_{\widehat{\U}_\perp}\A\|_q,\|\A P_{\widehat{\V}_\perp}\|_q \right\} \leq 2 \|\Z_{\max(r)}\|_q.
			\end{equation*} Here $\|\cdot\|_q$ denotes the matrix Schatten-$q$ norm.
		\end{Lemma}

	\end{sloppypar}
	
\end{document}